\documentclass[preprint]{elsarticle}

\usepackage[top=1.2in, bottom=1.2in, left=1.5in, right=1.5in]{geometry}
\usepackage{soul}
\usepackage{color}

\usepackage[most]{tcolorbox}
\usepackage{lineno}
\usepackage[allfiguresdraft]{draftfigure}
\usepackage{graphicx}
\usepackage[percent]{overpic}
\usepackage{epstopdf}
\DeclareGraphicsExtensions{.eps}
\usepackage{empheq}
\usepackage{scalerel,amssymb}
\usepackage{stmaryrd}
\usepackage{array}
\usepackage{amsmath}
\usepackage{leftidx}
\usepackage{amsfonts}
\usepackage{multicol}
\usepackage{mathrsfs}
\usepackage{tensor}
\usepackage{subfig}
\usepackage{hhline}
\usepackage{upgreek}
\usepackage{cancel}
\usepackage{ulem}
\usepackage{multirow}
\usepackage{setspace}

\usepackage{framed}
\usepackage{nomencl}
\makenomenclature

\usepackage[pagebackref=true,
            colorlinks=true,
            bookmarks=true,
           ]{hyperref}
           
\def\onedot{$\mathsurround0pt\ldotp$}
\def\cddot{
	\mathbin{\vcenter{\baselineskip.67ex
			\hbox{\onedot}\hbox{\onedot}}%
}}

\DeclareMathAlphabet\mathbfcal{OMS}{cmsy}{b}{n}
\newcommand{\vecgrad}{\vec{\nabla}}

\newcommand{\bs}{\boldsymbol}

\newcommand{\stt}[1]{\text{\scriptsize #1}}
\newcommand{\tightoverset}[2]{\mathop{#2}\limits^{\vbox to -.5ex{\kern-0.75ex\hbox{$#1$}\vss}}}

\newlength{\dheight}
\newcommand{\doubletilde}[1]{\smash{%
	\settoheight{\dheight}{\ensuremath{\tilde{#1}}}%
	\addtolength{\dheight}{-0.25ex}%
	\tilde{\vphantom{\rule{1pt}{\dheight}}%
		\smash{\tilde{\bs{#1}}}}\,\;}}		
\newcommand{\tildevec}[1]{\smash{%
	\settoheight{\dheight}{\ensuremath{\vec{#1}}}%
	\addtolength{\dheight}{-0.25ex}%
	\tilde{\vphantom{\rule{1pt}{\dheight}}%
		\smash{\vec{#1}}}\,\;}}
\newcommand{\hatvec}[1]{\smash{%
	\settoheight{\dheight}{\ensuremath{\vec{#1}}}%
	\addtolength{\dheight}{-0.25ex}%
	\hat{\vphantom{\rule{1pt}{\dheight}}%
		\smash{\vec{\bs{#1}}}}\,\;}}

\modulolinenumbers[5]

\journal{Elsevier}

\begin{document}

\allowdisplaybreaks[4]

\begin{frontmatter}
\title{Mesostructural origins of the anisotropic compressive properties of low-density closed-cell foams: A deeper understanding}

\author[1]{\corref{cor}L. Liu}\ead{lei.liu@chalmers.se}
\author[1]{F. Liu}
\author[2]{D. Zenkert}
\author[2]{M. {\AA}kermo}
\author[1]{M. Fagerstr\"{o}m}
\address[1]{Department of Industrial and Materials Science, Chalmers University of Technology, SE-41296 Gothenburg, Sweden}
\address[2]{Department of Engineering Mechanics, KTH Royal Institute of Technology, SE-10044 Stockholm, Sweden}
\cortext[cor]{Corresponding author.}

\begin{abstract}
Many closed-cell foams exhibit an elongated cell shape in the foam rise direction, resulting in anisotropic compressive properties, e.g.\ modulus and strength. Nevertheless, the underlying deformation mechanisms and how cell shape anisotropy induces this mechanical anisotropy are not yet fully understood, in particular for the foams with a high cell face fraction and low relative density. Moreover, the impacts of mesostructural stochastics are often overlooked.

This contribution conducts a systematic numerical study on the anisotropic compressive behaviour of \textcolor{red}{low-density closed-cell foams (with a relative density $<0.15$)}, which accounts for cell shape anisotropy, cell structure and different mesostructural stochastics. Representative volume elements (RVE) of foam mesostructures are modelled, with cell walls described as Reissner-Mindlin shells in a finite rotation setting. A mixed stress-strain driven homogenization scheme is introduced, which allows for enforcing an overall uniaxial stress state. Uniaxial compressive loadings in different global directions are applied.

Quantitative analysis of the cell wall deformation behaviour confirms the dominant role of membrane deformation in the initial elastic region, while the bending contribution gets important only \textcolor{red}{after buckling, followed by membrane yielding}. Based on the identified deformation mechanisms, analytical models are developed that relates mechanical anisotropy to cell shape anisotropy. It is found that cell shape anisotropy translates into the anisotropy of compressive properties through three pathways, cell load-bearing area fraction, cell wall buckling strength and cell wall inclination angle. Besides, the resulting mechanical anisotropy is strongly affected by the cell shape anisotropy stochastics while almost insensitive to the cell size and cell wall thickness stochastics. \textcolor{red}{The present findings provide deeper insights into the relationships between the anisotropic compressive properties and mesostructures of low-density closed-cell foams}. 
\end{abstract}

\begin{keyword}
Closed-cell foams \sep anisotropic compressive properties \sep cell shape anisotropy \sep stochastic variations \sep Laguerre tessellation \sep strain energy partitioning
\end{keyword}

\end{frontmatter}



\setlength{\nomlabelwidth}{2.05cm}
\setlength{\nomitemsep}{-0.025cm}
\begin{table*}[ht!] 
\color{red}  
\begin{framed}
\nomenclature{$E, \nu, \sigma_\stt{y}$}{Base material Young's modulus, Poisson's ratio and yield stress}

\nomenclature{$E_\stt{w}, \sigma_\stt{c,w}, \sigma_\stt{y,w}$}{Cell wall membrane modulus, buckling strength and yield strength}
\nomenclature{$k_\stt{c}$}{Cell wall buckling coefficient}
\nomenclature{$L_\stt{w}, B_\stt{w}, t$}{Cell wall length, width and thickness}
\nomenclature{$\theta_\stt{w}$}{Cell wall inclination angle}
\nomenclature{$\mathcal{R}_\stt{w}$}{Cell wall aspect ratio}
\nomenclature{$\mathcal{R}_\stt{f},\mathcal{R}_\stt{c},\mathcal{R}_\theta$}{Cell load-bearing area fraction ratio, cell wall buckling strength ratio and cell wall inclination angle ratio}

\nomenclature{$\mathcal{I}_\stt{w},\mathcal{B}_\stt{w}$}{Cell wall strain energy partitioning indicator and buckling detector}
\nomenclature{$\mathcal{J}_\stt{w},\mathcal{Y}_\stt{w}$}{Cell wall membrane plasticity indicator and yielding detector}

\nomenclature{$L_\stt{v}$}{Cell dimension in the global direction}
\nomenclature{$d_\stt{v}$}{Cell equivalent diameter}
\nomenclature{$\mathcal{R}_\stt{v}$}{Cell shape anisotropy}
\nomenclature{$\mathcal{R}_\stt{v}^E,\mathcal{R}_\stt{v}^\sigma$}{Cell compressive modulus anisotropy and strength anisotropy}

\nomenclature{$\hat{E}, \hat{\nu}, \hat{\sigma}_\stt{y}$}{RVE compressive modulus, Poisson's ratio and yield strength}

\nomenclature{$L$}{RVE dimension in the global direction}
\nomenclature{$\mathcal{R}$}{RVE cell shape anisotropy}
\nomenclature{$\mathcal{R}^E,\mathcal{R}^\sigma$}{RVE compressive modulus anisotropy and strength anisotropy}

\nomenclature{$\mathbf{F}$}{Mesoscale deformation gradient tensor}
\nomenclature{$\hat{\mathbf{F}}$}{Macroscale deformation gradient tensor}
\nomenclature{$\doubletilde{\mathbf{H}}^\stt{c}, \vec{G}^\stt{c}, \mathbf{K}^\stt{c}$}{Mesoscale shell membrane strain tensor, transverse shear strain vector and bending curvature tensor}
\nomenclature{$\mathbf{L},\mathbf{K}$}{Mesoscale shell cross-sectional deformation gradient tensor and its through-thickness gradient}
\nomenclature{$\mathbf{N},\mathbf{M}$}{Mesoscale shell stress resultant tensor and couple-stress resultant tensor}
\nomenclature{$\doubletilde{\mathbf{N}}^\stt{c}, \vec{V}^\stt{c}, \mathbf{M}^\stt{c}$}{Mesoscale shell membrane stress resultant tensor, transverse shear stress resultant vector and bending moment tensor}
\nomenclature{$\mathbf{P}$}{Mesoscale first Piola-Kirchhoff stress tensor}
\nomenclature{$\hat{\mathbf{P}}$}{Macroscale first Piola-Kirchhoff stress tensor}

\printnomenclature

\end{framed}
\color{black}
\end{table*}

\section{Introduction}\label{sec:introduction}
\begin{figure}[ht]
\centering
\begin{overpic}[draft=false,width=0.71\textwidth]{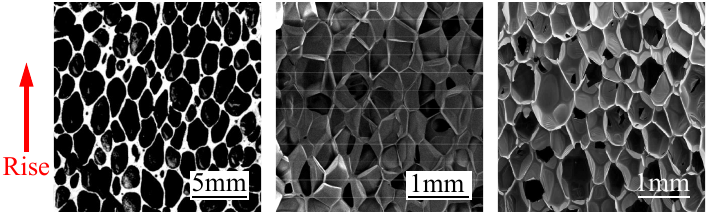}
\put(20.5,31.25){\color{black}\footnotesize \textbf{(a)}}
\put(51,31.25){\color{black}\footnotesize \textbf{(b)}}
\put(82.5,31.25){\color{black}\footnotesize \textbf{(c)}}
\end{overpic}
\caption{Examples of closed-cell foam mesostructures made from different base materials: (a) aluminium, (b) polyvinylklorid (PVC) and (c) polyisocyanurate (PIR) foams. The red arrow indicates the foam rise direction. Reproduced from \cite{Mu2010}, \cite{Zhou2023a} and \cite{Andersons2016}, respectively, with permission from Elsevier.}\label{fig:examples_meso_foam}
\end{figure}

Closed-cell foams are widely utilized in modern engineering applications due to their appealing specific mechanical properties with respect to low density, e.g.\ high stiffness and strength, and great energy absorption capacity \cite{Smith2012, Sun2018, Rahimidehgolan2023}. These properties are attributed to the underlying mesostructure, which consists of a large number of cells isolated by thin cell walls (see Figure~\ref{fig:examples_meso_foam}). During the foaming process, cells elongate in the foam rise direction, resulting in an anisotropic cell shape \cite{Mu2010, Zhou2023a}. Cell walls are usually thicker around the edges and thinner close to the face centers \cite{Jang2015, Tang2022}, commonly described by the cell edge/face material partitioning \cite{Gibson1997_ch5}. \textcolor{red}{For some closed-cell foams, cell walls may have apparent initial curvature with wriggles and corrugations, and even be missing \cite{Andrews1999, Jeon2005, PerezTamarit2019}}. Moreover, many mesostructural features, e.g.\ relative density, cell shape, cell size and cell wall thickness, are highly variable \cite{Jang2015, Ghazi2020a, Zhou2023a}. All the above lead to a broad spectrum of mechanical properties.

Given the exploitation of closed-cell foams for load-bearing applications, the compressive behaviour is often of interest \cite{Smith2012, Sun2018, Rahimidehgolan2023}. For most elasto-plastic foams, the compressive stress-strain response can be divided into three regions: elasticity, plateau and densification \cite{Gibson1997_ch5}. The first region is governed by the elastic membrane (or stretching) and bending deformations of cell walls. As the load increases, cell walls start to buckle elastically or collapse plastically. The elastic buckling and plastic collapse are localized failure modes, which occur first in the weakest cell walls and gradually propagate through the entire mesostructure, resulting in a plateau region with the compressive stress almost constant. Besides, the cell wall elastic buckling is the leading failure mode for low-density foams, while plastic collapse \textcolor{red}{is} the corresponding failure mode for high-density foams \cite{Tan2005, Zenkert2009, Michailidis2011, Kidd2012, Koohbor2018, Duan2019}. \textcolor{red}{The critical transition relative density depends on the base material properties and is also influenced by mesostructural features \cite{Michailidis2011, Kidd2012, Kader2017, Duan2019}}. Finally, upon densification the severely deformed cell walls come to contact and interact, leading to a rapid increase of the compressive stiffness.

To guide the closed-cell foam design, numerous studies have been conducted on the structure-property relationships. Among different mesostructural features, relative density is recognized as the most important in determining the compressive modulus and strength \cite{Gibson1997_ch5}. The relationships between these compressive properties and relative density have been well established. They can be expressed through power functions regardless of specific mechanisms \cite{Gibson1997_ch5}. Other features are implicitly accounted for by a set of constants of proportionality, usually identified from the experimental data. These relationships have demonstrated great success for a variety of foams (see e.g.\ aluminium \cite{Benouali2005, Bafti2013, Cheng2018}, polyvinylklorid (PVC) \cite{Saha2005, Tang2022, Zhou2023b}, aluminium composite \cite{Mondal2009, Guo2015}, ceramic \cite{Kim2005} and carbon \cite{Celzard2010} foams). 

Many closed-cell foams exhibit apparent anisotropic properties under compression (see e.g.\ aluminium \cite{Deshpande2001, EdwinRaj2009, Mu2010, Linul2018}, PVC \cite{Liu2020, Tang2022, Zhou2023b} and polyurethane (PU) \cite{Hamilton2013, MarviMashhadi2018, Li2019} foams). For instance, the compressive modulus and strength in the foam rise direction (see Figure~\ref{fig:examples_meso_foam}), are noticeably higher than the transverse direction. This mechanical anisotropy has been understood to primarily originate from cell shape anisotropy (see e.g.\ \cite{Deshpande2001, EdwinRaj2009, Mu2010, Hamilton2013, Linul2018, MarviMashhadi2018, Li2019, Liu2020, Tang2022, Zhou2023b}), while base material anisotropy plays a secondary role \cite{Linul2013}. These facts motivate detailed investigations on the impacts of cell shape anisotropy, arguably the second most important mesostructural feature for tailoring the compressive properties.

Compared with relative density, precisely controlling cell shape anisotropy is hardly possible in experiments, and thus micromechanical modelling is often employed. By idealizing a foam mesostructure as rectangular parallelepiped cell, Gibson and Ashby \cite{Gibson1997_ch6} pioneeringly proposed a semi-analytical model to predict the anisotropy of compressive properties in terms of cell shape anisotropy. \textcolor{red}{In this model, it is assumed that the cell edge bending accompanied by the face tension along the direction perpendicular to the compressive loading, and the cell wall plastic collapse are the dominant deformation and failure modes, respectively. Accordingly, the effective compressive properties in different global directions can be expressed in terms of the base material properties, cell wall thickness and cell sizes, followed by the mechanical anisotropy expressions}. Later on, Gong et al.\ \cite{Gong2005}, Sullivan et al.\ \cite{Sullivan2008} and Andersons et al.\ \cite{Andersons2016} improved the Gibson-Ashby model by introducing Kelvin cell, which could more accurately represent a foam mesostructure\footnote{These Kelvin cell-based analytical models are in principle developed for open-cell foams despite being applied for closed-cell foams in many studies.}. The idealized cell-based analytical models have been widely applied for realistic foams, showing capabilities to capture the general trends in the experimental data (see e.g.\ \cite{Gibson1997_ch6, Gong2005, Sullivan2008, EspadasEscalante2015, Andersons2016, Doyle2019, Liu2020, Zhou2023b}). However, the predictive deviations vary significantly from one case to another (sometimes $>100\%$), and are commonly regarded to arise from different uncertainties in the real foam mesostructures and experiments. Limited attention is paid to the mechanistic assumptions that have been introduced, which conflict with a few detailed experimental observations. For example, the cell wall elastic buckling, rather than plastic collapse, dominates the failure of many low-density foams (see e.g.\ aluminium \cite{Michailidis2011, Kader2017}, PVC \cite{Poapongsakorn2011, Kidd2012, Luong2013, Concas2019}, PU \cite{Koohbor2018, Bolintineanu2021} and polymethacrylimide (PMI) \cite{Zenkert2009, Chai2020} foams). \textit{This asks for a deeper understanding of the underlying mechanisms as well as the impacts of mesostructural features.}

To investigate the foam deformation behaviour in detail, finite element (FE) micromechanical modelling has been extensively performed. First, numerical models based on the idealized cell structures are developed (see e.g.\ rectangular \cite{Santosa1998}, Kelvin \cite{Simone1998a, Simone1998b, Grenestedt1998, Grenestedt2000, DeGiorgi2010, Sadek2013, Chen2018, Duan2019, Shakibanezhad2022} and Weaire–Phelan \cite{Chen2018, Shakibanezhad2022} cells), allowing for a systematic study of different mesostructural features and mechanisms. For example, Simone and Gibson \cite{Simone1998a} reported that both the compressive modulus and strength did not vary significantly against the cell edge/face material partitioning, suggesting that closed-cell foams deformed primarily by the cell wall stretching. \textcolor{red}{Grenestedt and Bassinet \cite{Grenestedt1998, Grenestedt2000} found that the cell wall curvature and thickness stochastics only weakly affected the compressive modulus, likely because the cell wall membrane deformation was largely involved the initial elastic region. Follow-up studies by Simone and Gibson \cite{Simone1998b} showed that compared with the cell wall curvature, the corrugations resulted in more pronounced reduction on the compressive modulus and strength. This may be because the cell wall corrugations promote bending deformation more effectively}. Chen et al.\ \cite{Chen2018} and Duan et al.\ \cite{Duan2019} confirmed that the failure of low-density foams was triggered by the cell wall elastic buckling, in alignment with experimental observations \cite{Zenkert2009, Michailidis2011, Kidd2012, Koohbor2018, Duan2019, Chai2020}. To the authors' best knowledge, the idealized cell-based numerical models are rarely employed to investigate the impacts of cell shape anisotropy (see one study in \cite{Sadek2013}, where no detailed mechanistic discussion is given).

Along with the advancement of computer tomography (CT) techniques, CT-based numerical models have also been developed (see e.g.\ \cite{Caty2008, Jeon2010, Sulong2015, Natesaiyer2015, Sun2017, Chen2017b, Talebi2019, Ghazi2020a}). These models provide a high-fidelity tool to study the underlying mechanisms. Sun et al.\ \cite{Sun2017} showed that the minimal ratio of the cell wall thickness to cell size determined the weakest region, where the first collapse (or crush) band formed under compression. Similar results were reported by Chen et al.\ \cite{Chen2017b} and Ghazi et al.\ \cite{Ghazi2020a}, that the larger and thinner cell walls tended to buckle earlier, followed by plastic deformation, eventually developing into the collapse bands. Most CT-based models are discretized by turning the voxels into cubic elements (see e.g.\ \cite{Natesaiyer2015, Sun2017}) or tetrahedral elements after geometric reconstruction (see e.g.\ \cite{Jeon2010, Sulong2015, Chen2017b}), leading to high computational costs. Therefore, CT-based models discretized by shell elements have been proposed, exhibiting excellent computational efficiency while preserving accuracy \cite{Caty2008}. Nevertheless, due to the inflexibility of manipulating the geometrical configurations, CT-based models are rarely used to systematically investigate the impacts of mesostructural features.

To fairly approximate the real foam mesostructures and meanwhile preserve the flexibility of manipulation, tessellation-based numerical models have received the most attention. Using Voronoi tessellation techniques, cell shape irregularity and randomness can be included. Song et al.\ \cite{Song2010} compared the results obtained using the tessellation-based and idealized cell-based models, showing that the compressive strength decreased along with cell shape irregularity. Further studies by Shi et al.\ \cite{Shi2018} and Vengatachalam et al.\ \cite{Vengatachalam2019} revealed that this strength reduction was attributed to the emergence of weak regions induced by cell shape irregularity. In contrast, the compressive modulus receives limited influence from cell shape irregularity, indicating that the cell wall membrane deformation dominates the initial elastic region \cite{Shi2018} (see also \cite{Grenestedt2000}). Roberts and Garboczi \cite{Roberts2001}, and K\"{o}ll and Hallstr\"{o}m \cite{Koll2016} studied the impacts of cell edge/face material partitioning on the compressive modulus, and observed substantial mismatch between the numerical data and Gibson-Ashby model predictions. It is pointed out that the cell wall membrane contribution is non-negligible. 

More recently, Laguerre tessellation-based models have been developed, which enable to incorporate the cell size stochastics. Chen et al.\ \cite{Chen2015, Chen2017a} showed that both the compressive modulus and strength decreased as the cell size and cell wall thickness stochastics increased, which again could be explained using the weakest link principle (see also \cite{Shi2018, Vengatachalam2019}). Compared with the compressive modulus, the strength is more sensitive to these mesostructural stochastics. \textcolor{red}{Marvi-Mashhadi et al.\ \cite{MarviMashhadi2020} showed that the entrapped gas inside cells generally stiffened the compressive response. Yet, this effect is nearly invisible in the elastic and early plateau regions (see also e.g.\ \cite{Sun2015, Zhang2015}), and thus becomes secondary for the compressive modulus and strength}. By elongating the original tessellation structures, the anisotropic foam mesostructural models, can be generated (see e.g.\ \cite{Su2018, Gebhart2019, MarviMashhadi2018, MarviMashhadi2020, Su2022, HossingerKalteis2022, Zhou2023b, Ding2023}). Gahlen and Stommel et al.\ \cite{Gahlen2022a, Gahlen2022b} further improved them to prescribe cell shape anisotropy stochastics. With the cell shape anisotropy control, the Laguerre tessellation-based models have shown great success to reproduce the anisotropic compressive stress-strain curves, even in quantitative agreement with the experimental data. 

\textcolor{red}{In addition, more general techniques based on inclusion packings have been proposed \cite{Sonon2015, Ghazi2019, Ghazi2020b}. With a control on the cell elongation, and cell wall thickness (linked to cell size) and curvature, these techniques can create the foam mesostructural models with arbitrary-shaped cells, providing a flexible representation of realistic foams. Ghazi et al.\ \cite{Ghazi2019, Ghazi2020b} showed that the compressive properties were strongly affected by the cell wall thickness stochastics, while less sensitive to the initial curvature and presence of missing cell walls.}

Nevertheless, very few of the above numerical studies provide quantitative analysis of the cell wall deformation behaviour and elaborate how cell shape anisotropy leads to mechanical anisotropy. Attempts have been made in \cite{MarviMashhadi2018, MarviMashhadi2020, Ding2023}, on PU foams where $>80\%$ of base materials are occupied by cell edges\footnote{These foams can be modelled in practice as open-cell foams since the cell face contribution is negligible.}. Through detailed analysis of the cell edge forces, it is revealed that the compressive load applied in the foam rise direction is initially carried by the cell edge axial deformation. The compressive load applied in the transverse direction is carried by both axial and bending deformations, leading to less stiff response and buckling at a lower applied stress. These findings are consistent with the experimental observations in \cite{Li2019} and rationalize the anisotropic compressive properties of PU foams. Yet, the obtained insights may not be representative for many foams with a high cell face fraction (likely $>$ 0.8, see e.g.\ aluminium \cite{Deshpande2001, EdwinRaj2009, Mu2010, Linul2018}, PVC \cite{Liu2020, Tang2022, Zhou2023a} and PMI \cite{Chai2020, Huo2022} foams).

To summarize, it is believed that the anisotropic compressive properties of closed-cell foams mainly originate from cell shape anisotropy. Analytical models have been proposed in the literature which relate mechanical anisotropy to cell shape anisotropy. However, the introduced mechanistic assumptions may not be valid for the foams with a high cell face fraction and low relative density. In these cases, the cell face contribution gets crucial and the cell wall elastic buckling becomes the leading failure mode. Extensive numerical studies have further suggested that the cell wall membrane deformation dominates the initial elastic region, which, nevertheless, are lacking confirmation, especially for anisotropic foams. Recently, attempts have been made to unravel the anisotropic compressive properties through quantitative analysis of the cell wall behaviour. However, these studies focus on the foams with a low cell face fraction and thus may not be representative for many other foams. \textit{More importantly, the detailed relationships between mechanical anisotropy and cell shape anisotropy remain unclear.} 

In addition, the intrinsic randomness of mesostructural features, \textcolor{red}{especially the cell wall thickness}, have been found to largely affect the compressive properties and may also impact mechanical anisotropy. \textit{These variations are usually overlooked when attempting to untangle the anisotropic compressive properties.} Accordingly, the present paper aims at addressing the following interconnected questions:
\begin{enumerate}
\item What are the key deformation mechanisms governing the anisotropic compressive behaviour of closed-cell foams with a high cell face fraction and low relative density?
\item How does cell shape anisotropy translate into the anisotropy of compressive properties?
\item Is this mechanical anisotropy influenced by the mesostructural stochastics?
\end{enumerate}
To the end, a systematic numerical study on the anisotropic compressive behaviour is conducted, which takes into account cell shape anisotropy, cell structure and the stochastic variations of different mesostructural features. Representative volume elements (RVE) of foam mesostructures are modelled\footnote{Due to the large randomness of mesostructural features, RVE should be interpreted as statistical volume element (SVE) \cite{OstojaStarzewski2006}. We retain the term RVE for consistency with convention in the foam community.}, where cell walls are described as Reissner-Mindlin shells \cite{Reissner1961} in a finite rotation setting \cite{Campello2003}. A mixed stress-strain driven homogenization scheme (see e.g.\ \cite{vanDijk2016, Saadat2023, Larsson2023}) is adopted to formulate the RVE problem such that an overall uniaxial stress state can be enforced. \textcolor{red}{Besides, to quantify the cell wall deformation behaviour, a strain energy partitioning indicator followed by a buckling detector, and a membrane plasticity indicator followed by a yielding detector, are proposed}.

Rectangular parallelepiped cell structures with different shape anisotropy are first modelled. Second, Kelvin cell structures are modelled which further account for the cell wall inclination angle. Third, foam mesostructures generated using Laguerre tessellation techniques are modelled which incorporate the stochastic variations of cell size, cell wall thickness and cell shape anisotropy. Based on the numerical analyses of the two idealized cell-based models, analytical models are derived which relates the anisotropy of compressive properties to cell shape anisotropy. The mechanical anisotropy determined from the tessellation-based numerical models is compared against the present analytical model predictions, to identify the influence of mesostructural stochastics. Besides, two widely used analytical models \cite{Gibson1997_ch6, Sullivan2008} are assessed, to investigate the importance of imposing appropriate mechanistic assumptions. \textcolor{red}{This work will provide a deeper understanding on the relationships between the anisotropic compressive properties and mesostructures for low-density closed-cell foams}.

The paper is organized as follows. In \textcolor{blue}{Section}~\ref{sec:problem_description}, the foam mesostructural RVE problem is formulated. In \textcolor{blue}{Section}~\ref{sec:quan_method_cell_wall}, a quantification method for the cell wall deformation behaviour is introduced. In \textcolor{blue}{Section}~\ref{sec:num_model_method}, RVE simulation setup is elaborated. In \textcolor{blue}{Section}~\ref{sec:analy_ideal_meso_model}, numerical results of the idealized cell-based models (rectangular parallelepiped and Kelvin cells) are analysed. Based on the insights obtained, analytical models that describe the relationships between mechanical anisotropy and cell shape anisotropy are developed in \textcolor{blue}{Section}~\ref{sec:relation_anisotropy}. In \textcolor{blue}{Section}~\ref{sec:analy_tess_meso_model}, numerical results of the tessellation-based models are analysed. The extracted anisotropy of compressive properties is compared against the analytical model predictions as well as experimental data, followed by discussions on the impacts of mesostructural stochastics. The main conclusions are summarized in \textcolor{blue}{Section}~\ref{sec:conclusions}.

Scalars, vectors, second-order tensors and forth-order tensors in this paper are denoted as e.g.\ $a$, $\vec{a}$, $\mathbf{A}$ and $\mathbb{A}$, respectively. The length of a vector (Euclidean norm) is denoted by $|| \vec{\bullet} ||$. The transpose and inverse of a second-order tensor are denoted by $(\bullet)^\stt{T}$ and $(\bullet)^\stt{-1}$, respectively. The macroscale effective quantities are denoted by $(\hat{\bullet})$.

\section{RVE problem description}\label{sec:problem_description}
\begin{figure}[ht]
\centering
\begin{overpic}[draft=false,width=0.95\textwidth]{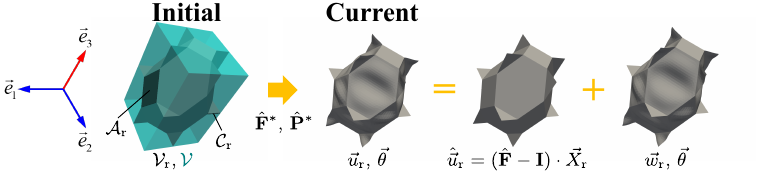}
\put(35.5,20.375){\color{black}\footnotesize \textbf{(a)}}
\put(66.75,20.375){\color{black}\footnotesize \textbf{(b)}}
\put(86.625,20.375){\color{black}\footnotesize \textbf{(c)}}
\end{overpic}
\caption{A foam mesostructural RVE with cell walls described as shell continuum: (a) initial to current configurations after imposing the macroscale stress $\hat{\mathbf{P}}$ and deformation gradient $\hat{\mathbf{F}}$ in a mixed manner. $(\bullet)^*$ indicates a quantity with its components partially prescribed. The space-filling volume domain including voids are indicated by the green shadows; decomposition of the mesoscale mid-surface displacement field $\vec{u}_\stt{r}$ into the (b) trend field $\hatvec{u}_\stt{r}$ and (c) fluctuation field $\vec{w}_\stt{r}$. The rotation angle field $\vec{\theta}$ is not visible.}\label{fig:rve_deform_configuration_decomp}
\end{figure}

Consider a foam mesostructural RVE, with its space-filling volume domain (shaded green) indicated by $\mathcal{V}$ and base material volume domain by $\mathcal{V}_\stt{r}$, see Figure~\ref{fig:rve_deform_configuration_decomp}(a). Only a single cell structure is illustrated for simplicity. Reissner–Mindlin shell description \cite{Reissner1961} is adopted for the cell walls. The shell reference mid-surface is denoted by $\mathcal{A}_\stt{r}$ and external boundary contour by $\mathcal{C}_\stt{r}$. The global coordinate system $\{\vec{e}_1, \vec{e}_2, \vec{e}_3\}$ is chosen such that $\vec{e}_3$ is parallel to the foam rise direction, and $\vec{e}_1$ and $\vec{e}_2$ the two transverse directions.

\subsection{Shell kinematics and stress resultants}
To capture geometrically nonlinear behaviour of cell walls, a finite rotation shell formulation \cite{Campello2003} is adopted. Let introduce a curvilinear coordinate system $\{\vec{e}_1^{\stt{\,c}},\vec{e}_2^{\stt{\,c}},\vec{e}_3^{\stt{\,c}}\}$, such that the plane $\{\vec{e}_1^{\stt{\,c}},\vec{e}_2^{\stt{\,c}}\}$ is tangent to the mid-surface, while $\vec{e}_3^{\stt{\,c}}$ is normal to this tangent plane. Position vectors of any material point in the initial and current configurations are given by
\begin{subequations}\label{eq:positions_initial_current}
	\begin{align}
	&\vec{X}=\vec{X}_\stt{r}+\eta{\vec{D}},\quad \eta\in\mathcal{H}\\
	&\vec{x}=\vec{x}_\stt{r}+\eta{\vec{d}},\quad \eta\in\mathcal{H}
	\end{align}
\end{subequations}
\color{red}
where $\vec{X}_\stt{r}$ and $\vec{x}_\stt{r}$ define a point on the mid-surface before and after deformation, respectively; $\vec{D}$ denotes the director in the initial configuration, i.e.\ $\vec{D}=\vec{e}_3^{\stt{\,c}}$; $\vec{d}$ denotes the director in the current configuration; $\eta$ is the through-thickness coordinate, belonging to the thickness domain $\mathcal{H}=[-\frac{t}{2},\frac{t}{2}]$, with $t$ the thickness which can vary from one location of $\vec{X}_\stt{r}$ to another.

Following the Reissner–Mindlin theory that the director remains straight after deformation but not necessarily perpendicular to the deformed mid-surface (enabling to account for the transverse shear effect), $\vec{d}$ is related to $\vec{D}$ through
\begin{equation}\label{eq:directors_initial_current}
\vec{d}=\mathbf{R}(\vec{\theta})\cdot\vec{D},
\end{equation}
where $\mathbf{R}$ denotes a rotation tensor in terms of the Euler rotation angle vector $\vec{\theta}$ (e.g.\ according to the Euler–Rodrigues formula \cite{Campello2003}).
\color{black}

Subtracting eq.~(\ref{eq:positions_initial_current}a) from eq.~(\ref{eq:positions_initial_current}b) gives the displacement field as
\begin{equation}\label{eq:displacement}
\vec{u}=\vec{u}_\stt{r}+(\mathbf{R}(\vec{\theta})-\mathbf{I})\cdot\eta\vec{D},
\end{equation}
implying that the finite rotation shell kinematics can be fully parametrized using $\vec{u}_\stt{r}$ and $\vec{\theta}$.

The deformation gradient tensor follows from eq.~\eqref{eq:positions_initial_current} as (see \cite{Coenen2010, Liu2021} for details)
\begin{equation}\label{eq:deformation_gradient}
\mathbf{F}=(\vecgrad_{0}\vec{x})^\stt{T} =\left(\frac{\partial\vec{x}_\stt{r}}{\partial\vec{X}_\stt{r}}\right)^\stt{T}+\eta\left(\frac{\partial\vec{d}}{\partial\vec{X}_\stt{r}}\right)^\stt{T}+\vec{d}\otimes\vec{D},
\end{equation}
where $\vecgrad_{0}={\partial(\bullet)}/{\partial\vec{X}_\stt{r}}+{\partial(\bullet)}/{\partial(\eta\vec{D})}$ denotes the gradient operator with respect to $\vec{X}$. Substituting eq.~\eqref{eq:directors_initial_current} into eq.~\eqref{eq:deformation_gradient} yields
\begin{equation}\label{eq:deformation_gradient_2}
\mathbf{F}=\mathbf{L}+\eta\mathbf{K},
\end{equation}
with
\begin{subequations}\label{eq:membrane_bending}
	\begin{align}
	&\mathbf{L}=(\tildevec{\nabla}_0\otimes\vec{x}_\stt{r})^\stt{T}+\mathbf{R}(\vec{\theta})\cdot\vec{D}\otimes\vec{D},\\
	&\mathbf{K}=\bs{\varGamma}(\vec{\theta})\cdot\left(\tildevec{\nabla}_0\otimes(\vec{\theta}\times\vec{D})\right)^\stt{T},
	\end{align}
\end{subequations}
where $\tildevec{\nabla}_0={\partial(\bullet)}/{\partial\vec{X}_\stt{r}}$ denotes the gradient operator with respect to $\vec{X}_\stt{r}$; \textcolor{red}{$\bs{\varGamma}$ denotes a rotation curvature tensor in terms of $\vec{\theta}$ (see \cite{Campello2003} for details)}.

For the convenience of constitutive model formulation, a back-rotated configuration is introduced by eliminating $\mathbf{R}$. The back-rotated counterpart of $\mathbf{F}$ is given by
\begin{equation}\label{eq:back_rot_deformation_gradient}
\mathbf{F}^\stt{c}=\mathbf{R}^\stt{T}\cdot\mathbf{F}=\mathbf{I}+\mathbf{H}^\stt{c}+\eta\mathbf{K}^\stt{c},
\end{equation}
where $\mathbf{H}^\stt{c}=\mathbf{R}^\stt{T}\cdot\mathbf{L}-\mathbf{I}$ and $\mathbf{K}^\stt{c}=\mathbf{R}^\stt{T}\cdot\mathbf{K}$ represent the cross-sectional generalized strain and bending curvature, respectively
\begin{subequations}\label{eq:disp_grad_bend_curv}
	\begin{align}
	&\mathbf{H}^\stt{c}=	
	\mathbf{R}^\stt{T}(\vec{\theta})\cdot(\tildevec{\nabla}_0\otimes\vec{x}_\stt{r})^\stt{T}-\doubletilde{\mathbf{I}},\\
	&\mathbf{K}^\stt{c}=\bs{\varGamma}^\stt{T}(\vec{\theta})\cdot\left(\tildevec{\nabla}_0\otimes(\vec{\theta}\times\vec{D})\right)^\stt{T},
	\end{align}
\end{subequations}
with $\doubletilde{\mathbf{I}}=(\tildevec{\nabla}_0\otimes\vec{X}_\stt{r})^\stt{T}=\mathbf{I}-\vec{D}\otimes\vec{D}$. Here, equation~\eqref{eq:membrane_bending} has been substituted and the property $\bs{\varGamma}^\stt{T}=\mathbf{R}^\stt{T}\cdot\bs{\varGamma}$ has been applied. The membrane strain tensor $\doubletilde{\mathbf{H}}^\stt{c}$ and transverse shear strain vector $\vec{G}^\stt{c}$ can be identified from $\mathbf{H}^\stt{c}$ as
\begin{subequations}\label{eq:memb_disp_grad_trans_shear_strain}
	\begin{align}
	&\doubletilde{\mathbf{H}}^\stt{c}=\doubletilde{\mathbf{I}}\cdot\mathbf{H}^\stt{c}, \\
	&\vec{G}^\stt{c}=\vec{D}\cdot\mathbf{H}^\stt{c}.
	\end{align}
\end{subequations}

\textcolor{red}{Denoting the first Piola–Kirchhoff stress tensor by $\mathbf{P}$ and thus its back-rotated counterpart $\mathbf{P}^\stt{c}=\mathbf{R}^\stt{T}\cdot\mathbf{P}$, the resultants conjugate to $\doubletilde{\mathbf{H}}^\stt{c}$, $\vec{G}^\stt{c}$ and $\mathbf{K}^\stt{c}$ are defined as}
\begin{subequations}\label{eq:membrane_trans_stress_bend_mom}
	\begin{align}
	&\doubletilde{\mathbf{N}}^\stt{c}
=\doubletilde{\mathbf{I}}\cdot\int_{\mathcal{H}}\mathbf{P}^\stt{c}\, {\rm d}\eta, \\
	&\vec{V}^\stt{c}
=\vec{D}\cdot\int_{\mathcal{H}}\mathbf{P}^\stt{c}\, {\rm d}\eta, \\
	&\mathbf{M}^\stt{c}=\int_{\mathcal{H}}\eta\mathbf{P}^\stt{c}\, {\rm d}\eta, 
	\end{align}
\end{subequations}
representing the membrane stress resultant tensor, transverse shear stress resultant vector and bending moment tensor, respectively. The constitutive relations between $\{\doubletilde{\mathbf{N}}^\stt{c}, \vec{V}^\stt{c}, \mathbf{M}^\stt{c}\}$ and $\{\doubletilde{\mathbf{H}}^\stt{c}, \vec{G}^\stt{c}, \mathbf{K}^\stt{c}\}$ will be specified in \textcolor{blue}{Section}~\ref{sec:num_model_method}. 

Similarly, the resultants conjugate to $\mathbf{L}$ and $\mathbf{K}$ are defined as
\begin{subequations}\label{eq:stress_resultants}
	\begin{align}
	&\mathbf{N}=\int_{\mathcal{H}}\mathbf{P}\, {\rm d}\eta, \\
	&\mathbf{M}=\int_{\mathcal{H}}\eta\mathbf{P}\, {\rm d}\eta. 
	\end{align}
\end{subequations}
\textcolor{red}{Applying the principle $\int_{\mathcal{H}}\delta\mathbf{F}^\stt{T}\cddot\mathbf{P}\, {\rm d}\eta=\int_{\mathcal{H}}(\delta\mathbf{F}^\stt{c})^\stt{T}\cddot\mathbf{P}^\stt{c}\, {\rm d}\eta$ and accounting for the plane stress state in the thickness direction}, $\mathbf{N}$ and $\mathbf{M}$ can be expressed in terms of $\doubletilde{\mathbf{N}}^\stt{c}$, $\vec{V}^\stt{c}$ and $\mathbf{M}^\stt{c}$ as
\begin{subequations}\label{eq:stress_resultants_2}
	\begin{align}
	&\mathbf{N}=\mathbf{R}\cdot(\doubletilde{\mathbf{N}}^\stt{c}+\vec{D}\otimes\vec{V}^\stt{c}), \\
	&\mathbf{M}=\mathbf{R}\cdot\mathbf{M}^\stt{c}. 
	\end{align}
\end{subequations}

\subsection{Weak form of the balance equations}
The mid-surface displacement $\vec{u}_\stt{r}$ and rotation angle $\vec{\theta}$ are commonly adopted as the primary field variables to be solved for a finite rotation shell problem. In absence of body force, the weak form can be stated as: find $\{\vec{u}_\stt{r},\vec{\theta}\}$ such that for all admissible variations $\{\delta\vec{u}_\stt{r},\delta\vec{\theta}\}$, the balance holds between the internal and external virtual works:
\begin{equation}\label{eq:weak_form}
	\int_{\mathcal{A}_\stt{r}}\delta{\mathbf{L}}^\stt{T}(\vec{u}_\stt{r},\vec{\theta})\cddot{\mathbf{N}}\,{\rm{d}}A
    + \int_{\mathcal{A}_\stt{r}} \delta{\mathbf{K}}^\stt{T}(\vec{\theta})\cddot{\mathbf{M}}\,{\rm{d}}A 
	= \int_{\mathcal{C}_{\stt{r}}}\delta\vec{u}_\stt{r}\cdot\vec{n}\,{\rm{d}}C
    + \int_{\mathcal{C}_{\stt{r}}}\delta\vec{\theta}\cdot\left(\bs{\varGamma}^\stt{T}(\vec{\theta})\cdot\vec{m}\right)\,{\rm{d}}C,\quad\forall\{\delta\vec{u}_\stt{r},\delta\vec{\theta}\},
\end{equation}
with the traction resultant and bending moment on the external boundary $\mathcal{C}_\stt{r}$ given by
\begin{equation}
\vec{n}=\mathbf{N}\cdot\vec{N}_\stt{r},\quad \vec{m}=\mathbf{M}\cdot\vec{N}_\stt{r}\quad  {\text{on }}\mathcal{C}_{\rm{r}}.
\end{equation}
Here $\vec{N}_\stt{r}$ is the outward normal to $\mathcal{C}_{\rm{r}}$. \textcolor{red}{The boundary conditions required to complete the problem~\eqref{eq:weak_form} will be subsequently determined through downscaling.}

\subsection{Standard strain driven formulation}
Following the classical homogenization (see e.g.\ \cite{Kouznetsova2001, Miehe2002}), applying the macroscale effective deformation gradient tensor $\hat{\mathbf{F}}$ to the foam mesostructural RVE (see Figure~\ref{fig:rve_deform_configuration_decomp}(a)) yields the mesoscale relative position vector field, over the space-filling volume domain $\mathcal{V}$
\begin{equation}\label{eq:meso_relative_position}
\Delta\vec{x}=\hat{\mathbf{F}}\cdot\Delta\vec{X}+\Delta\vec{w},\quad \vec{X}\in\mathcal{V}
\end{equation}
with $\Delta\vec{X}=\vec{X}-\vec{X}_\stt{o}$, $\Delta\vec{x}=\vec{x}-\vec{x}_\stt{o}$ and  $\Delta\vec{w}=\vec{w}-\vec{w}_\stt{o}$, where $\vec{X}_\stt{o}$ denotes the initial position vector of a reference origin point; the current position vector of this point is denoted by $\vec{x}_\stt{o}$; $\vec{w}$ reflects the fluctuations induced by heterogeneities. 

\textcolor{red}{The classical downscaling requires that the volume average of the mesoscale deformation gradient ${\mathbf{F}}$ is equated to $\hat{\mathbf{F}}$}:
\begin{equation}\label{eq:downscaling}
\hat{\mathbf{F}}=\frac{1}{V}\int_{\mathcal{V}}\mathbf{F}\,{\rm{d}}V,
\end{equation}
with $V$ being the RVE space-filling volume (see Figure~\ref{fig:rve_deform_configuration_decomp}(a)). Substituting eqs.~\eqref{eq:deformation_gradient} and~\eqref{eq:meso_relative_position} into eq.~\eqref{eq:downscaling}, followed by applying divergence theorem, gives the constraint in terms of the fluctuations $\Delta\vec{w}$
\begin{equation}\label{eq:constraints_flu}
\int_{\partial\mathcal{V}/\partial\mathcal{V}_\stt{r}} \Delta\vec{w}\otimes\vec{N}\,{\rm{d}}\partial V
+
\int_{\partial\mathcal{V}_\stt{r}} \Delta\vec{w}\otimes\vec{N}\,{\rm{d}}\partial V=\mathbf{0}.
\end{equation}
where the total external surface domain $\partial\mathcal{V}$ has been split into the void $\partial\mathcal{V}/\partial\mathcal{V}_\stt{r}$ and base material $\partial\mathcal{V}_\stt{r}$ parts. \textcolor{red}{Notice that $\Delta\vec{w}$ is not available on $\partial\mathcal{V}/\partial\mathcal{V}_\stt{r}$}. Therefore, the following choice is made to satisfy constraint~\eqref{eq:constraints_flu}
\begin{subequations}\label{eq:constraints_flu_decomp}
	\begin{align}
	&\int_{\partial\mathcal{V}/\partial\mathcal{V}_\stt{r}} \Delta\vec{w}\otimes\vec{N}\,{\rm{d}}\partial V=\mathbf{0}, \\
	&\int_{\partial\mathcal{V}_\stt{r}} \Delta\vec{w}\otimes\vec{N}\,{\rm{d}}\partial V=\mathbf{0}.
	\end{align}
\end{subequations}
The first constraint can be fulfilled by appropriate choice of $\Delta\vec{w}$ on $\partial\mathcal{V}/\partial\mathcal{V}_\stt{r}$, which in practice does not affect the RVE solution. 

Constraint~(\ref{eq:constraints_flu_decomp}b) is next elaborated. For the convenience of derivations, let choose the reference origin point $\vec{X}_\stt{o}$ to be located within the mid-surface of cell walls. Moreover, the fluctuation at $\vec{X}_\stt{o}$ will be constrained, i.e.\ $\vec{w}_\stt{o}=\vec{0}$, to eliminate the rigid-body translation. 

Consistent with the shell kinematics~\eqref{eq:positions_initial_current}, $\vec{w}$ at any material point can be expressed as
\begin{equation}\label{eq:meso_flu_decomp}
\vec{w}=\vec{w}_\stt{r}+\left(\mathbf{R}(\vec{\theta})-\hat{\mathbf{R}}\right)\cdot\eta\vec{D},\quad \eta\in\mathcal{H}
\end{equation}
where $\vec{w}_\stt{r}$ collects the mid-surface displacement fluctuations, while the other term collects the director fluctuations; $\hat{\mathbf{R}}$ is the rotation part of $\hat{\mathbf{F}}$. Substituting eq.~\eqref{eq:meso_flu_decomp} into  eq.~(\ref{eq:constraints_flu_decomp}b) leads to the constraint in terms of $\vec{w}_\stt{r}$ and $\vec{\theta}$:
\begin{equation}\label{eq:constraints_flu_1}
\int_{\mathcal{C}_\stt{r}}\int_{\mathcal{H}}\left(\vec{w}_\stt{r}+\left(\mathbf{R}(\vec{\theta})-\hat{\mathbf{R}}\right)\cdot\eta\vec{D}\right)\otimes\vec{N}_\stt{r}\,{\rm{d}}\eta{\rm{d}}C=\mathbf{0},
\end{equation}
\textcolor{red}{implying a coupling between $\vec{w}_\stt{r}$ and $\vec{\theta}$ on the external boundaries in general}. Further taking into account the thickness domain $\mathcal{H}=[-\frac{t}{2},\frac{t}{2}]$, constraint~\eqref{eq:constraints_flu_1} can be simplified into
\begin{equation}\label{eq:constraints_flu_2}
\int_{\mathcal{C}_\stt{r}}t\vec{w}_\stt{r}\otimes\vec{N}_\stt{r}\,{\rm{d}}C=\mathbf{0},
\end{equation}
\textcolor{red}{implying that the constraints on $\vec{\theta}$ are not compulsory for the symmetric shell formulation.}

The standard strain driven homogenization scheme assumes that $\hat{\mathbf{F}}$ is fully known, and thus the mesoscale mid-surface displacement fluctuation field $\vec{w}_\stt{r}$ can be replaced by
\begin{equation}\label{eq:meso_mid_surf_flu}
\vec{w}_\stt{r}=\vec{u}_\stt{r}-(\hat{\mathbf{F}}-\mathbf{I})\cdot\vec{X}_\stt{r}.
\end{equation}
Substituting eq.~\eqref{eq:meso_mid_surf_flu} into eq.~\eqref{eq:constraints_flu_2} yields minimal kinematic boundary conditions in terms of the mesoscale mid-surface displacements $\vec{u}_\stt{r}$
\begin{equation}\label{eq:cond_min_kin}
\int_{\mathcal{C}_\stt{r}}t\vec{u}_\stt{r}
\otimes\vec{N}_\stt{r}\,{\rm{d}}C=(\hat{\mathbf{F}}-\mathbf{I})\cdot\int_{\mathcal{C}_\stt{r}}t\vec{X}_\stt{r}
\otimes\vec{N}_\stt{r}\,{\rm{d}}C.
\end{equation}
The other common choice to fulfil constraint~\eqref{eq:constraints_flu_2} is fully prescribed boundary conditions, which can be obtained by enforcing $\vec{w}_\stt{r}=\vec{0}$ on the external boundary $\mathcal{C}_{\rm{r}}$
\begin{equation}\label{eq:cond_fully_pres}
\vec{u}_\stt{r}=(\hat{\mathbf{F}}-\mathbf{I})\cdot\vec{X}_\stt{r} \quad  {\text{on }}\mathcal{C}_{\rm{r}}
\end{equation}
\textcolor{red}{For an RVE with geometrical periodicity, periodic boundary conditions is frequently adopted}
\begin{equation}\label{eq:cond_periodic}
\vec{u}_\stt{r}^{-}-\vec{u}_\stt{r}^{+}=(\hat{\mathbf{F}}-\mathbf{I})\cdot(\vec{X}_\stt{r}^{-}-\vec{X}_\stt{r}^{+}) \quad  {\text{on }}\mathcal{C}_{\rm{r}}
\end{equation}
Here the superscripts ``$-/+$" denote the opposite boundary pair. \textcolor{red}{Notice that extra conditions $\vec{\theta}^{-}=\vec{\theta}^{+}$ are introduced in many foam mesostructural RVE studies (see e.g.\ \cite{Chen2018, Ghazi2020b, Gahlen2022a}, where no detailed derivations are given). This choice intuitively follows from the periodic boundary conditions for the solid continuum problem}.

Choosing one of the conditions~\eqref{eq:cond_min_kin},~\eqref{eq:cond_fully_pres} and~\eqref{eq:cond_periodic} completes the RVE problem~\eqref{eq:weak_form}. After solving the RVE problem, the macroscale effective stress $\hat{\mathbf{P}}$ can be derived using the Hill-Mandel condition \cite{Hill1963}, as the volume average of the mesoscale stress $\mathbf{P}$
\begin{equation}\label{eq:upscaling}
\hat{\mathbf{P}}=\frac{1}{V}\int_{\mathcal{V}}\mathbf{P}\,{\rm{d}}V.
\end{equation}
Accounting for $\mathbf{P}=\mathbf{0}$ in the void part $\mathcal{V}/\mathcal{V}_\stt{r}$, substituting eq.~(\ref{eq:stress_resultants}a) into eq.~\eqref{eq:upscaling} gives $\hat{\mathbf{P}}$ in terms of the mesoscale resultants $\mathbf{N}$:
\begin{equation}\label{eq:eff_stress}
\hat{\mathbf{P}}=\frac{1}{V}\int_{\mathcal{A}_\stt{r}}\mathbf{N}\,{\rm{d}}A
=\frac{1}{V}\int_{\mathcal{C}_\stt{r}}\vec{n}\otimes\vec{X}_\stt{r}\,{\rm{d}}C.
\end{equation}
Here, divergence theorem has been applied to obtain the second equality. Clearly, the mesoscale resultant $\mathbf{M}$ (see eq.~(\ref{eq:stress_resultants}b)) does not contribute to $\hat{\mathbf{P}}$.

\subsection{Mixed stress-strain driven formulation}
The RVE response under a specific stress state, e.g.\ overall stress free or uniaxial stress (to be considered in this study) state, is of interest in some cases. This is, however, not straightforward to \textcolor{red}{(precisely) enforce} with the strain driven formulation. Therefore, several mixed stress-strain driven formulations for the solid continuum problem that enable to prescribe $\hat{\mathbf{P}}$ or $(\hat{\mathbf{F}},\hat{\mathbf{P}})$ in a mixed manner, have been proposed in the literature (see e.g.\ \cite{vanDijk2016, Saadat2023, Larsson2023}). Following the same spirit, the early introduced strain driven formulation for the shell problem is extended to a mixed stress-strain driven formulation. \textcolor{red}{To the authors' best knowledge, a mixed stress-strain driven formulation for the shell problem has not been reported in the literature}.

As sketched in Figure~\ref{fig:rve_deform_configuration_decomp} (see also eq.~\eqref{eq:meso_mid_surf_flu}), the mesoscale mid-surface displacement field $\vec{u}_\stt{r}$ consists of two contributions:
\begin{equation}\label{eq:meso_mid_surf_disp_decomp}
\vec{u}_\stt{r}=\hatvec{u}_\stt{r}+\vec{w}_\stt{r}=(\hat{\mathbf{F}}-\mathbf{I})\cdot\vec{X}_\stt{r}+\vec{w}_\stt{r},
\end{equation}
where the first term represents a linear, trend field induced by $\hat{\mathbf{F}}$ while the second term $\vec{w}_\stt{r}$ the fluctuation field. Considering $\hat{\mathbf{F}}$, $\vec{w}_\stt{r}$ and $\vec{\theta}$ as the primary field variables instead, substituting eq.~\eqref{eq:meso_mid_surf_disp_decomp} into eq.~\eqref{eq:weak_form} allows to reformulate the RVE problem as:
\begin{equation}\label{eq:weak_form_re}
	\begin{split}
	\int_{\mathcal{A}_\stt{r}}\delta{\mathbf{L}}^\stt{T}(\hat{\mathbf{F}},\vec{w}_\stt{r},\vec{\theta})\cddot{\mathbf{N}}\,{\rm{d}}A
    + \int_{\mathcal{A}_\stt{r}} \delta{\mathbf{K}}^\stt{T}(\vec{\theta})\cddot{\mathbf{M}}\,{\rm{d}}A \\
	= V\delta\hat{\mathbf{F}}^\stt{T}\cddot\hat{\mathbf{P}}+
	\int_{\mathcal{C}_{\stt{r}}}\delta\vec{w}_\stt{r}\cdot\vec{n}\,{\rm{d}}C
    + \int_{\mathcal{C}_{\stt{r}}}\delta\vec{\theta}\cdot\left(\bs{\varGamma}^\stt{T}(\vec{\theta})\cdot\vec{m}\right)\,{\rm{d}}C,\quad\forall\{\delta\hat{\mathbf{F}},\delta\vec{w}_\stt{r},\delta\vec{\theta}\}.
	\end{split}
\end{equation}
Here, equation~\eqref{eq:eff_stress} has been applied to the right-hand side as such $\hat{\mathbf{P}}$ appears like an ``external force", which can be fully or partially prescribed.

Accordingly, the RVE boundary conditions~\eqref{eq:cond_min_kin},~\eqref{eq:cond_fully_pres} and~\eqref{eq:cond_periodic} derived from the strain driven homogenization shall be reformulated in terms of $\vec{w}_\stt{r}$. Constraint~\eqref{eq:constraints_flu_2} can be directly adopted as minimal kinematic boundary conditions. Fully prescribed boundary conditions and periodic boundary conditions become $\vec{w}_\stt{r}=\vec{0}$ and $\vec{w}_\stt{r}^{-}=\vec{w}_\stt{r}^{+}$, respectively.

\section{Quantification method of the cell wall behaviour}\label{sec:quan_method_cell_wall} 
To directly identify the deformation mechanisms, a quantification method for the cell wall deformation behaviour is developed. As discussed in \textcolor{blue}{Section}~\ref{sec:introduction}, cell wall stretching and bending are the two most commonly observed deformation modes. Accordingly, one choice is to measure the partition of strain energy in cell walls into membrane and bending modes (see also e.g.\ \cite{Alkhader2009, Ding2024}). A cell wall-wise (marked by the subscript ``w") strain energy partitioning indicator $\mathcal{I}_\stt{w}$ is defined through:
\begin{equation}\label{eq:strain_energy_par_ind}
\mathcal{I}_\stt{w}=\frac{W_\stt{b}-W_\stt{m}}{W_\stt{b}+W_\stt{m}},
\end{equation}
with $W_\stt{m}$ and $W_\stt{b}$ being the membrane and bending energy, respectively. $\mathcal{I}_\stt{w}=-1$ indicates a pure membrane deformation mode while $\mathcal{I}_\stt{w}=1$ a pure bending mode. Making use of the shell kinetic quantities (see eqs.~\eqref{eq:disp_grad_bend_curv}-\eqref{eq:membrane_trans_stress_bend_mom}), $W_\stt{m}$ and $W_\stt{b}$ are evaluated as
\begin{subequations}\label{eq:membrane_bending_strain_energy}
	\begin{align}
	&W_\stt{m}=\int_{\mathcal{A}_\stt{w}}\left(\int(\doubletilde{\mathbf{N}}^\stt{c})^\stt{T}\cddot{\rm d}\doubletilde{\mathbf{H}}^\stt{c}\right)\, {\rm d}A, \\
	&W_\stt{b}=\int_{\mathcal{A}_\stt{w}}\left(\int(\mathbf{M}^\stt{c})^\stt{T}\cddot{\rm d}\mathbf{K}^\stt{c}\right)\, {\rm d}A, 
	\end{align}
\end{subequations}
with $\mathcal{A}_\stt{w}$ being the mid-surface of a probed cell wall. Notice that the proposed indicator~\eqref{eq:strain_energy_par_ind} is applicable to an arbitrary constitutive model choice, in contrast to those in the literature where small strain and isotropic elasticity are assumed \cite{Alkhader2009, Ding2024}.

\begin{figure}[ht]
\centering
\begin{overpic}[draft=false,width=0.9\textwidth]{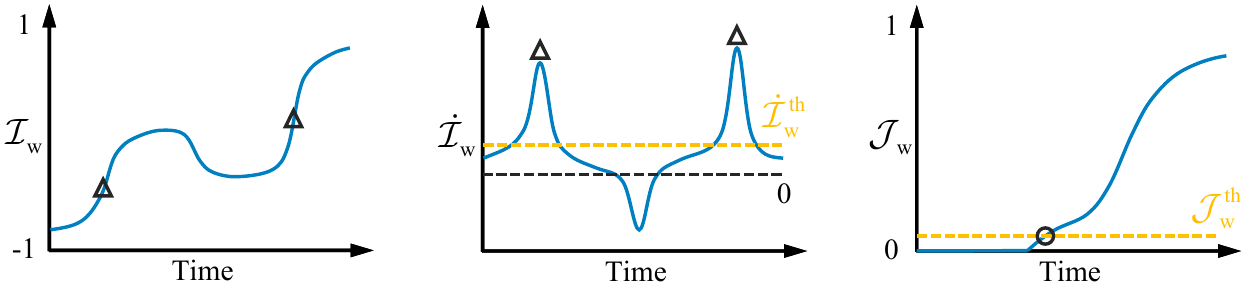}
\put(14.25,21){\color{black}\footnotesize \textbf{(a)}}
\put(49.25,21){\color{black}\footnotesize \textbf{(b)}}
\put(84.25,21){\color{black}\footnotesize \textbf{(c)}}
\end{overpic}
\caption{\textcolor{red}{Sketches of the (a) strain energy partitioning indicator $\mathcal{I}_\stt{w}$, (b) strain energy partitioning indicator rate $\dot{\mathcal{I}}_\stt{w}$ and (c) membrane plasticity indicator $\mathcal{J}_\stt{w}$ versus time, for a probed cell wall. The buckling and membrane yielding events are indicated by the black triangles and circles, respectively.}}\label{fig:indicator_time}
\end{figure}

For low-density foams, the cell wall elastic buckling is the main failure mode (see \textcolor{blue}{Section}~\ref{sec:introduction}). Therefore, the occurrence of buckling is tracked. As buckling is often accompanied by a sharp transition from the membrane to bending mode, a sudden increase on the indicator profile $\mathcal{I}_\stt{w}$ is expected, as sketched in Figure~\ref{fig:indicator_time}(a). This sudden increase may be characterized as a positive peak on the rate profile $\dot{\mathcal{I}}_\stt{w}$, as sketched in Figure~\ref{fig:indicator_time}(b). A cell wall-wise buckling detector $\mathcal{B}_\stt{w}$ is proposed, which is defined in a time-wise way:
\begin{equation}\label{eq:BK_detector}
\mathcal{B}_\stt{w}^n=
\begin{cases}
1 & {\rm if}\,\, \dot{\mathcal{I}}_\stt{w}^n>{\rm max}\{\dot{\mathcal{I}}_\stt{w}^{n-1},\dot{\mathcal{I}}_\stt{w}^{n+1}\}\,\,{\rm and}\,\,\dot{\mathcal{I}}_\stt{w}^n>\dot{\mathcal{I}}_\stt{w}^\stt{th} \\
0 & {\rm else}
\end{cases},
\end{equation} 
where $n$ indicates the probed time step, and a threshold $\dot{\mathcal{I}}_\stt{w}^\stt{th}$ has been introduced to extract remarkable peaks only. The threshold value is proposed as the time average of positive $\dot{\mathcal{I}}_\stt{w}$. A cell wall is considered to be buckled as long as $\mathcal{B}_\stt{w}=1$ appears once.

\color{red}
It should be emphasized that different from beams, plates (or shells) after buckling remain capable to take additional load to a large extent, before membrane yielding occurs \cite{Timoshenko1961}. Accordingly, it becomes relevant to identify when the cell wall enters the plastic regime and yields. Considering that the membrane energy is dominating under compression, one choice is to measure the cell wall portion where the membrane plastic deformation has started. A membrane plasticity indicator $\mathcal{J}_\stt{w}$ is defined through:
\begin{equation}\label{eq:membrane_plas_ind}
\mathcal{J}_\stt{w}=\frac{A_\stt{p}}{A_\stt{w}},
\end{equation}
with $A_\stt{p}$ and $A_\stt{w}$ being the plastic and total areas of $\mathcal{A}_\stt{w}$, respectively. Notice that we do not model plasticity explicitly (to be discussed in \textcolor{blue}{Section}~\ref{sec:num_model_method}). Instead, von Mises plasticity criterion is assumed such that $A_\stt{p}$ can be evaluated by post-processing $\doubletilde{\mathbf{N}}^\stt{c}$
\begin{equation}\label{eq:plas_area}
A_\stt{p}=\int_{\mathcal{A}_\stt{w}}\mathcal{F}\stt{p}\left(\frac{\doubletilde{\mathbf{N}}^\stt{c}}{t}\right)\, {\rm d}A,
\end{equation}
where $\mathcal{F}_\stt{p}$ denotes the plasticity flag, i.e.\ $\mathcal{F}_\stt{p}=0$ before plasticity and $\mathcal{F}_\stt{p}=1$ after plasticity. 

Furthermore, by introducing a small threshold $\mathcal{J}_\stt{w}^\stt{th}$, as sketched in Figure~\ref{fig:indicator_time}(c), a membrane yielding detector $\mathcal{Y}_\stt{w}$ (in the post-buckling regime) is proposed:
\begin{equation}\label{eq:Y_detector}
\mathcal{Y}_\stt{w}^n=
\begin{cases}
1 & {\rm if}\,\,\mathcal{J}_\stt{w}^n>\mathcal{J}_\stt{w}^\stt{th} \\
0 & {\rm else}
\end{cases}.
\end{equation} 
The threshold value is taken to be 0.01. A cell wall is considered to have reached yielding and thereby reached its load-bearing capacity, as long as $\mathcal{Y}_\stt{w}=1$ appears once.
\color{black}

\section{RVE numerical simulation setup}\label{sec:num_model_method}

\subsection{Geometrical model configurations}
To systematically investigate the influence of different foam mesostructural features, RVE models with different levels of complexity are considered. \textcolor{red}{Since closed-cell foams with a high cell face fraction and low relative density are focused on, constant thickness is assumed for each cell wall (see also e.g.\ \cite{Chen2015, Vengatachalam2019, Zhou2023b}).}

\textcolor{red}{A preliminary study on the impacts of cell wall curvature has been particularly conducted, which shows a more complex relationship between the cell wall deformation behaviour and initial curvature level than what has been reported previously (see e.g.\ \cite{Grenestedt1998, Simone1998b, Ghazi2019}). Nevertheless, these impacts shall be negligible at least for polymer foams where most cell walls have a small normalized curvature ($< 0.1$). More details can be found in~\ref{sec:influence_cell_wall_curv}. Therefore, for the sake of simplicity, and without too much loss of generality for the focused foams, we model each cell wall as a flat plate. Besides, the presence of missing cell walls is not explicitly modelled. This effect can be implicitly accounted for by assigning low base material properties (see e.g.\ \cite{Ghazi2019}), or a small thickness and thus partially included by the cell wall thickness stochastics.}

For each cell (marked by the subscript ``v"), the shape of which is (approximately) transversely isotropic, shape anisotropy is defined as $\mathcal{R}_\stt{v}=(\mathcal{R}_\stt{v,31}\mathcal{R}_\stt{v,32})^{\frac{1}{2}}$, where the two ratios are given by $\mathcal{R}_\stt{v,31}=L_\stt{v,3}/L_\stt{v,1}$ and $\mathcal{R}_\stt{v,32}=L_\stt{v,3}/L_\stt{v,2}$. Here, $L_{{\rm v},i}$ denotes the cell dimension in the global direction $\vec{e}_i$. Besides, an equivalent diameter $d_\stt{v}=\left(\frac{6}{\pi}V_\stt{v}\right)^{\frac{1}{3}}$ is introduced, with $V_\stt{v}$ being the cell volume. The overall cell shape anisotropy for an RVE model consisting of multiple cells is evaluated as ${\mathcal{R}}=\frac{1}{N_{\rm v}}\sum\mathcal{R}_\stt{v}$, where $N_\stt{v}$ denotes the total number of cells. Finally, the overall relative density is evaluated as $\rho/\rho_\stt{r}=V_\stt{r}/V$, where $\rho$ and $\rho_\stt{r}$ represent the RVE and base material mass densities, respectively; $V=\sum V_\stt{v}=L_1L_2L_3$ the total volume of cells (or RVE space-filling volume), with $L_i$ the RVE dimension in the global direction $\vec{e}_i$; $V_\stt{r}=\sum tA_\stt{w}$ the total volume of cell walls, i.e.\ base material volume, $t$ the thickness and $A_\stt{w}$ the mid-surface area.

\begin{figure}[ht]
\centering
\begin{overpic}[draft=false,width=0.79\textwidth]{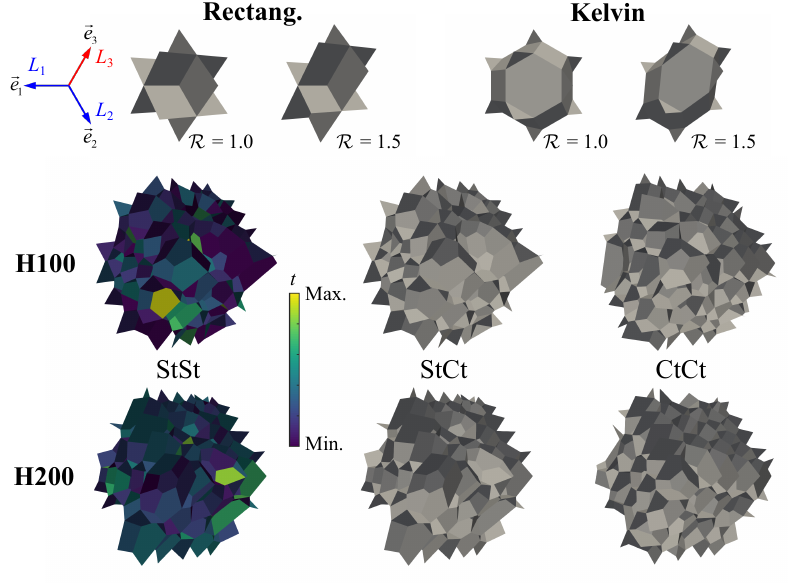}
\put(20.5,72.125){\color{black}\footnotesize \textbf{(a)}}
\put(67,72.125){\color{black}\footnotesize \textbf{(b)}}
\put(-3,40.25){\color{black}\footnotesize \textbf{(c)}}
\put(-3,13.25){\color{black}\footnotesize \textbf{(d)}}
\end{overpic}
\caption{Geometrical model configurations of different foam mesostructural RVE: (a-b) idealized cell-based models for two shape anisotropy. All cell walls are assigned with a constant thickness; (c-d) tessellation-based models. ``StSt" accounts for the stochastic variations of both cell size and cell wall thickness. ``StCt" accounts for the cell size stochastics while assigning a constant thickness to all cell walls. ``CtCt" assigns a constant equivalent diameter to all cells and a constant thickness to all cell walls.}\label{fig:geo_model_configuration}
\end{figure}
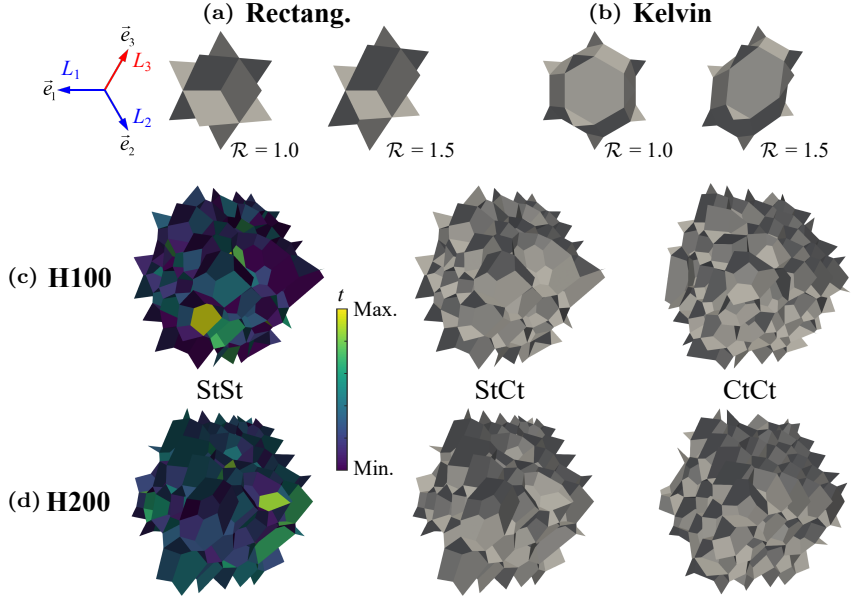

The idealized cell-based models are first introduced, including rectangular parallelepiped and Kelvin cells, which have been widely employed in the literature (see e.g.\ \cite{Santosa1998, Sadek2013}). The use of the rectangular parallelepiped cells allows for investigating the impacts purely by cell shape anisotropy, while Kelvin cells further take the cell wall inclination angle into account and better approximate the real foam mesostructures. Geometrical model configurations for the two shape anisotropy, $\mathcal{R}=1.0$ and $\mathcal{R}=1.5$, are shown in Figures~\ref{fig:geo_model_configuration}(a-b) as examples, with the geometrical parameters adopted for $\mathcal{R}=1.0$ listed in Table~\ref{tab:ideal_cell_geo_parameters}. 

By assuming periodicity, each RVE model represents a perfectly repeatable, infinite foam mesostructure. Model configurations with $\mathcal{R}$ from $1.0$ to $2.0$ are considered. $L_i$ is scaled according to the prescribed $\mathcal{R}$, with $V$ and $t$ preserved, i.e.\ $L_1=L_2=V^\frac{1}{3}\mathcal{R}^{\stt{-}\frac{1}{3}}$ and $L_3=V^\frac{1}{3}\mathcal{R}^{\frac{2}{3}}$. Note, that $\rho/\rho_\stt{r}$ would slightly increase as $\mathcal{R}$ increases. For rectangular parallelepiped cells, $\rho/\rho_\stt{r}$ varies from 0.0750 to 0.0788, while for Kelvin cells, \textcolor{red}{from 0.0753 to 0.0815}. The resulting $\rho/\rho_\stt{r}$ are representative for low-density foams (see e.g.\ \cite{Kader2017, Concas2019, Chai2020}).

\renewcommand\arraystretch{1.1}
\begin{table}[ht]
\centering
\caption{Geometrical parameters of the reference idealized cell-based models ($\mathcal{R}=1.0$) in Figures~\ref{fig:geo_model_configuration}(a-b).}
\label{tab:ideal_cell_geo_parameters}
\begin{tabular}{ll|ll}
\hline
\textbf{Parameter}          & \textbf{Symbol} & \textbf{Rectang.} & \textbf{Kelvin} \\ \hline
RVE dimension 1     & $L_1$ & 0.4 [mm] & 0.4 [mm] \\
RVE dimension 2     & $L_2$ & 0.4 [mm] & 0.4 [mm] \\
RVE dimension 3     & $L_3$ & 0.4 [mm] & 0.4 [mm] \\
Cell wall thickness & $t$ & 0.01 [mm] & \textcolor{red}{0.009} [mm] \\
Relative density    & $\rho/\rho_\stt{r}$ & 0.0750 [-] & \textcolor{red}{0.0753} [-] \\
\hline
\end{tabular}
\end{table}
\renewcommand\arraystretch{1}

Next, tessellation-based models are introduced which incorporate the stochastic variations of cell size, cell wall thickness and cell shape anisotropy observed in real foam mesostructures. Diab Divinycell foam H100 and H200 are considered, given their representativeness as low-density foams and availability of mesostructural characterization data in the literature (see e.g.\ \cite{Tang2022, Zhou2023a, Skeens2024}). The corresponding RVE geometrical model configurations are generated using Laguerre tessellation techniques, supported in the open-source package Neper \cite{Quey2018}, in accordance to the detailed experimental measurements \cite{Zhou2023a}. 

For each tessellation-based RVE model, the overall cell shape anisotropy $\mathcal{R}$ is prescribed \textcolor{red}{by uniformly transforming the original tessellation model (see also e.g.\ \cite{MarviMashhadi2018, Zhou2023b, Ding2023})}. The resulting shape anisotropy $\mathcal{R}_\stt{v}$ of individual cells approximately follows a normal distribution (see Figures~\ref{fig:PDF_meso_foam}(a-b)), associated with the cell shape irregularity. The cell equivalent diameters $d_\stt{v}$ are assigned using a log-normal distribution \cite{Zhou2023a}, with the mean $\mu_d$ and standard deviation $\sigma_d$. The cell wall thickness $t$ are assigned using a gamma distribution \cite{Zhou2023a}, with the mean $\mu_t$ and standard deviation $\sigma_t$. The adopted geometrical parameters are listed in Table~\ref{tab:tess_geo_parameters}, with the generated model set ``StSt" shown in Figures~\ref{fig:geo_model_configuration}(c-d). It has been validated that the numerically realized distributions of different mesostructural features agree well with the prescribed distributions (see Figure~\ref{fig:PDF_meso_foam}). More details can be found in~\ref{sec:meso_stochastics}. The resulting $\rho/\rho_\stt{r}$ for ``StSt" H100 and H200 are 0.0806 and 0.1480, respectively, which are slightly higher than the nominal values (see Table~\ref{tab:tess_geo_parameters}). Nevertheless, such differences are within the variation of $\rho/\rho_\stt{r}$ in practice, $+15$/$-10\%$ \cite{Diab2023}.

To investigate the influence of different mesostructural stochastics, two extra model sets ``StCt" and ``CtCt" are introduced (see Figures~\ref{fig:geo_model_configuration}(c-d)), with the cell shape anisotropy stochastics close to ``StSt". ``StCt" is obtained by prescribing a constant $t$ on ``StSt", i.e.\ $\sigma_t=0$ [mm]. The resulting $\rho/\rho_\stt{r}$ for ``StCt" are nearly the same as ``StSt". ``CtCt" is obtained by further enforcing a constant $d_\stt{v}$ on ``StCt", i.e.\ $\sigma_d=0$ [mm]. Notice that the resulting $\rho/\rho_\stt{r}$ for ``CtCt" H100 and H200 are 0.0922 and 0.1668, apparently higher than the nominal values (see Table~\ref{tab:tess_geo_parameters}). This implies the importance of incorporating the cell size stochastics in order to fairly approximate realistic foams.

\renewcommand\arraystretch{1.1}
\begin{table}[ht]
\centering
\caption{Geometrical parameters of the tessellation-based model set ``StSt" in Figures~\ref{fig:geo_model_configuration}(c-d).}
\label{tab:tess_geo_parameters}
\begin{tabular}{ll|ll}
\hline
\textbf{Parameter}          & \textbf{Symbol} & \textbf{H100} & \textbf{H200} \\ \hline
RVE dimension 1     & $L_1$ & 1.50 [mm] & 1.45 [mm] \\
RVE dimension 2     & $L_2$ & 1.50 [mm] & 1.45 [mm] \\
RVE dimension 3     & $L_3$ & 1.50 [mm] & 1.45 [mm] \\
Cell shape anisotropy & $\mathcal{R}$ & 1.2 [-] & 1.4 [-] \\
Cell equivalent diameter & $(\mu_d,\sigma_d)$ & $(0.35,0.10)$ [mm] & $(0.34,0.09)$  [mm] \\
Cell wall thickness & $(\mu_t,\sigma_t)$ & $(0.0115,0.0059)$ [mm] & $(0.0200,0.0067)$  [mm] \\
Nominal relative density    & $\rho/\rho_\stt{r}$ & 0.0714 [-] & 0.1429 [-] \\
\hline
\end{tabular}
\end{table}
\renewcommand\arraystretch{1}

Each RVE model is discretized in the open-source mesh generator Gmsh \cite{Geuzaine2009}, by triangular shell elements with the mid-surface displacement fluctuation $\vec{w}_\stt{r}$ and rotation $\vec{\theta}$ as degrees of freedom (DOF). To avoid shear locking, second-order Lagrange interpolation is adopted for $\vec{w}_\stt{r}$ while Crouzeix-Raviart interpolation for $\vec{\theta}$, as suggested in \cite{Campello2003}. For the idealized cell-based models, a fine mesh with averaged element size $\sim$ 0.01 [mm] is used to resolve the local deformation pattern in sufficient detail. For the tessellation-based models, a relatively large element size $\sim$ 0.03 [mm] is chosen to balance the computational accuracy and cost. A mesh sensitivity check has been performed for each tessellation-based model, and confirmed that the cell wall-wise and macroscale effective responses are both converged for the adopted discretization. \textcolor{red}{Besides, it has been confirmed that the chosen RVE dimensions (see Table~\ref{tab:tess_geo_parameters}) are sufficiently large (with $\sim$ 120 cells) to deliver the converged effective responses, even when different random realizations are considered}. More details can be found in~\ref{sec:influence_RVE_size_real}.

\subsection{Material model}
This study focuses on the anisotropy of compressive modulus and strength for low-density foams, where failure is mainly triggered by the cell wall elastic buckling (see \textcolor{blue}{Section}~\ref{sec:introduction}). \textcolor{red}{Therefore, plasticity is disregarded in the material modelling but remains involved in the cell wall behaviour analysis and foam yield strength analysis, through post-processing instead (see \textcolor{blue}{Section}~\ref{sec:quan_method_cell_wall}).}

A finite-strain isotropic elasticity is used to describe the base material behaviour. The bulk elasticity follows the linear relation $\mathbf{S}=\mathbb{C}\cddot\mathbf{E}$, where $\mathbf{S}=\mathbf{F}^\stt{-1}\cdot\mathbf{P}$ denotes the second Piola-Kirchhoff stress tensor, $\mathbf{E}=\frac{1}{2}(\mathbf{F}^\stt{T}\cdot\mathbf{F}-\mathbf{I})$ the Green–Lagrange strain tensor, and $\mathbb{C}$ the forth-order elasticity tensor fully determined by Young's modulus $E$ and Poisson's ratio $\nu$. Substituting the bulk elasticity into eqs.~\eqref{eq:disp_grad_bend_curv}-\eqref{eq:membrane_trans_stress_bend_mom}, followed by neglecting the higher-order terms, leads to the constitutive relations consistent with the finite rotation shell formulation:
\begin{subequations}\label{eq:shell_cons_relations}
	\begin{align}
	&\doubletilde{\mathbf{N}}^\stt{c}=\doubletilde{\mathbf{L}}^\stt{c}\cdot\mathbb{D}_\stt{m}\cddot\tfrac{1}{2}\left((\doubletilde{\mathbf{L}}^\stt{c})^\stt{T}\cdot\doubletilde{\mathbf{L}}^\stt{c}-\doubletilde{\mathbf{I}}\right), \\
	&\vec{V}^\stt{c}=\mathbf{D}_\stt{t}\cdot\vec{G}^\stt{c}, \\
	&\mathbf{M}^\stt{c}=\doubletilde{\mathbf{L}}^\stt{c}\cdot\mathbb{D}_\stt{b}\cddot\tfrac{1}{2}\left((\doubletilde{\mathbf{L}}^\stt{c})^\stt{T}\cdot{\mathbf{K}}^\stt{c}+
({\mathbf{K}}^\stt{c})^\stt{T}\cdot\doubletilde{\mathbf{L}}^\stt{c}\right),
	\end{align}
\end{subequations}
with $\doubletilde{\mathbf{L}}^\stt{c}=\doubletilde{\mathbf{I}}+\doubletilde{\mathbf{H}}^\stt{c}$. Here, $\mathbb{D}_\stt{m}=t\mathbb{C}^\stt{ps}$, $\mathbf{D}_\stt{t}=\frac{\kappa tE}{2(1+\nu)}\doubletilde{\mathbf{I}}$ and $\mathbb{D}_\stt{b}=\frac{t^3}{12}\mathbb{C}^\stt{ps}$ represent the generalized stiffness tensors for the membrane, transverse shear and bending modes, respectively; $\mathbb{C}^\stt{ps}$ denotes a degraded $\mathbb{C}$ by assuming a plane stress state in the thickness direction $\vec{e}_3^{\stt{\,c}}$; $\kappa=\frac{5}{6}$ is the shear correction factor. Besides, an element-wise fictitious stiffness $t^3E$ is added to stabilize the drilling rotation. It has been confirmed that the resulting drilling strain energy is sufficiently small.

The material model described above is implemented using the open-source code generator MFront \cite{Helfer2020}. In consistency with Divinycell foams, PVC is adopted as the base material with its parameters $E=2.7$ [GPa], $\nu=0.38$ and $\sigma_\stt{y}=62$ [MPa] taken from \cite{Zhou2023b}. \textcolor{red}{The yield stress $\sigma_\stt{y}$ is required for the use of the membrane plasticity indicator~\eqref{eq:membrane_plas_ind}. It has been confirmed that for the considered foams (with a relative density $<0.15$), the compressive properties of interest, i.e.\ modulus and strength, can be reasonably well determined using the elastic numerical results, although the experimentally observed plateau region cannot be captured. More details can be found in~\ref{sec:influence_plas}}.

\subsection{Boundary conditions and loads}
For the idealized cell-based models, to fulfil the periodicity assumptions, periodic boundary conditions are enforced on both $\vec{w}_\stt{r}$ and $\vec{\theta}$, i.e.\ $\vec{w}_\stt{r}^{-}=\vec{w}_\stt{r}^{+}$ and $\vec{\theta}^{-}=\vec{\theta}^{+}$ (see also \textcolor{blue}{Section}~\ref{sec:problem_description}). Since no periodicity is present for the tessellation-based models and minimal kinematic boundary conditions is usually too weak (see e.g.\ \cite{Miehe2002}), fully prescribed boundary conditions are enforced on $\vec{w}_\stt{r}$, i.e.\ $\vec{w}_\stt{r}=\vec{0}$ (see also \textcolor{blue}{Section}~\ref{sec:problem_description}). This choice is appropriate as long as the RVE dimensions are large enough, which is the case in the present study \textcolor{red}{(see~\ref{sec:influence_RVE_size_real} for details)}. 

To investigate the anisotropic compressive behaviour, uniaxial compressive loadings in different global directions are applied, by prescribing the macroscale effective quantities $(\hat{\mathbf{F}},\hat{\mathbf{P}})$ in a mixed manner. Loading case $\vec{e}_1$ is specified below as an example:
\begin{equation}\label{eq:uniaxial_compressive_load}
\hat{\mathbf{F}} =
\begin{bmatrix}
\lambda & 0          & 0 \\
\times       & \times    & 0 \\
\times       & \times         & \times
\end{bmatrix},\,\,
\hat{\mathbf{P}}=
\begin{bmatrix}
\times& \times        & \times \\
0       & 0 & \times \\
0       & 0         & 0
\end{bmatrix}.
\end{equation}
Here, $\hat{\mathbf{F}}$ and $\hat{\mathbf{P}}$ have been projected to the global basis $\vec{e}_i$; $\lambda$ denotes the prescribed stretch ratio and ``$\times$" the unprescribed (unknown) components. Notice that half of the non-diagonal components in $\hat{\mathbf{F}}$ have also been prescribed to eliminate the rigid-body rotation. 

The complete RVE problem~\eqref{eq:weak_form_re} is implemented by coupling the open-source computing platform FEniCS \cite{Logg2012, Bleyer2018} with MFront \cite{Helfer2020}, where $\hat{\mathbf{F}}$ is treated as global DOF in addition to local DOF $\vec{w}_\stt{r}$ and $\vec{\theta}$.

\section{Analyses of the idealized cell-based models}\label{sec:analy_ideal_meso_model}
Numerical results of the idealized cell-based models introduced in \textcolor{blue}{Section}~\ref{sec:num_model_method} will be analysed in this section, to exemplify the anisotropic compressive behaviour and deformation mechanisms of closed-cell foams, as well as how different compressive properties evolve with cell shape anisotropy, in absence of any mesostructural stochastics.

\subsection{Rectangular parallelepiped cell}
\begin{figure}[ht]
\centering
\begin{overpic}[draft=false,width=0.9\textwidth]{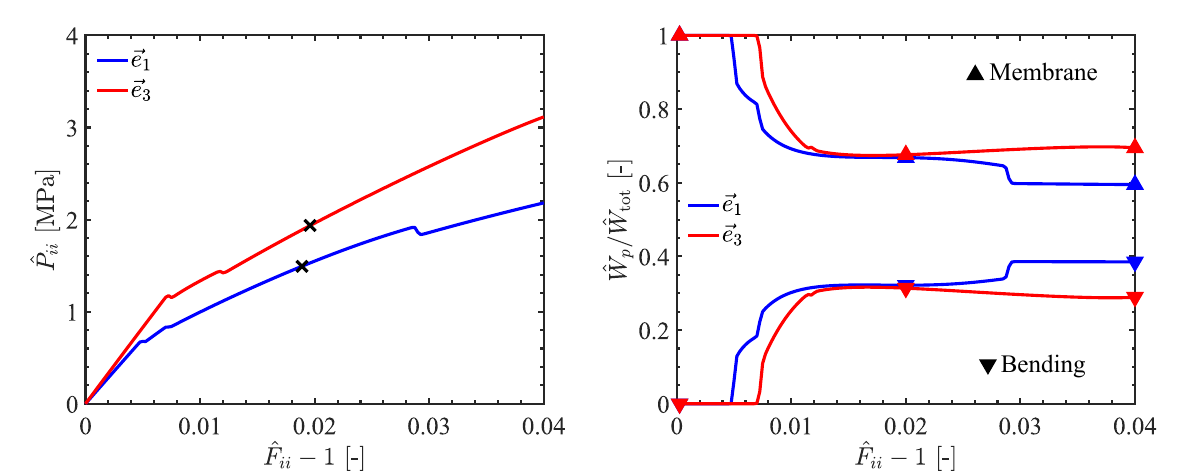}
\put(24.5,37.5){\color{black}\footnotesize \textbf{(a)}}
\put(73.625,37.5){\color{black}\footnotesize \textbf{(b)}}
\end{overpic}
\caption{Effective responses of the rectangular parallelepiped cell-based model with $\mathcal{R}=1.5$, under uniaxial compression in the transverse ($\vec{e}_1$) and foam rise ($\vec{e}_3$) directions: (a) stress and (b) strain energy fraction versus applied strain. The yield points are indicated in (a) by the black crosses.}\label{fig:eff_response_cube_cell_05}
\end{figure}

The macroscale effective stresses $\hat{\mathbf{P}}$ and strain energy fractions $\hat{W}_p/\hat{W}_\stt{tot}$ of different deformation modes, for shape anisotropy $\mathcal{R}=1.5$, are plotted as functions of the applied strain in Figure~\ref{fig:eff_response_cube_cell_05}. The effective responses for loading case $\vec{e}_2$ are identical to those for $\vec{e}_1$ and thus not displayed. Figure~\ref{fig:eff_response_cube_cell_05}(a) shows that for each loading case, the stress first increases linearly, \textcolor{red}{followed by multiple times of stress drops}. As expected, the stress under compression in the foam rise direction ($\vec{e}_3$) is apparently higher than the transverse direction ($\vec{e}_1$), indicating an anisotropic compressive behaviour. Figure~\ref{fig:eff_response_cube_cell_05}(b) shows that for each loading direction, the cell structure experiences a nearly pure membrane deformation mode ($\hat{W}_\stt{m}/\hat{W}_\stt{tot}\sim 1$) in the initial elastic region, while the bending mode gets pronounced after the first stress drop.

\begin{figure}[ht]
\centering
\begin{overpic}[draft=false,width=0.9\textwidth]{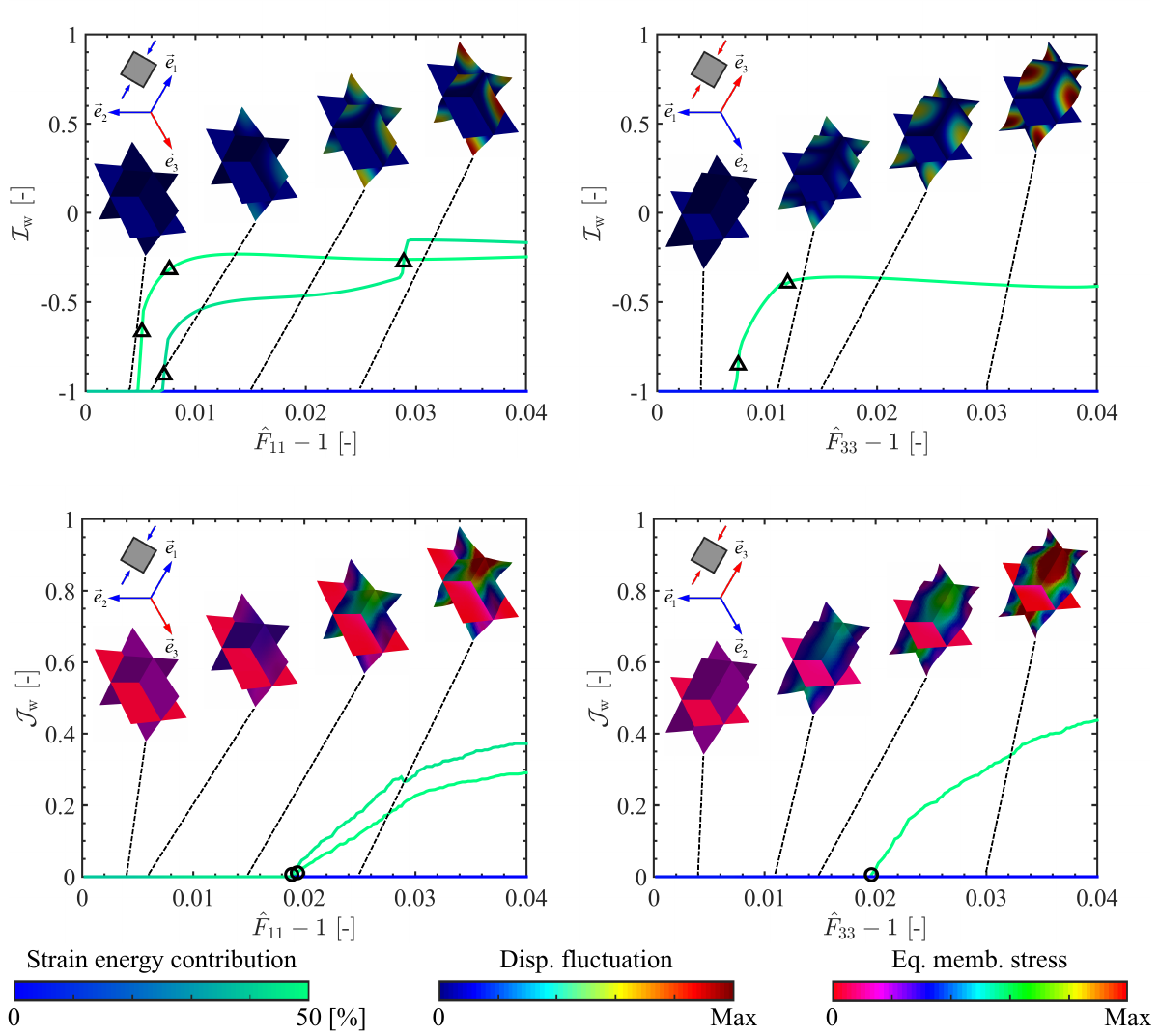}
\put(24.5,87.625){\color{black}\footnotesize \textbf{(a)}}
\put(73.625,87.625){\color{black}\footnotesize \textbf{(b)}}
\put(24.5,45.75){\color{black}\footnotesize \textbf{(c)}}
\put(73.625,45.75){\color{black}\footnotesize \textbf{(d)}}
\end{overpic}
\caption{\textcolor{red}{Cell wall strain energy partitioning indicators, membrane plasticity indicators, and deformed configurations at different stages of the rectangular parallelepiped cell-based model with $\mathcal{R}=1.5$, under uniaxial compression in the (a, c) transverse and (b, d) foam rise directions. The two green curves in (b, d) are overlapping. The buckling points in (a-b) and yield points in (c-d) are indicated by the black triangles and circles, respectively. Each loading direction is represented by a pair of opposite arrows}.}\label{fig:cell_wall_cube_cell_05}
\end{figure}

\textcolor{red}{To rationalize the observations in Figure~\ref{fig:eff_response_cube_cell_05}, individual cell wall deformation behaviours are quantified using the method introduced in \textcolor{blue}{Section}~\ref{sec:quan_method_cell_wall}, with the results reported in Figure~\ref{fig:cell_wall_cube_cell_05}. The strain energy partitioning indicator $\mathcal{I}_\stt{w}$ (using eq.~\eqref{eq:strain_energy_par_ind}) and membrane plasticity indicators $\mathcal{J}_\stt{w}$ (using eq.~\eqref{eq:membrane_plas_ind}) of each cell wall are plotted as functions of the applied strain, and colored by its strain energy contribution percentage. The deformation configurations are colored by the displacement fluctuations for better visualization of buckling in Figures~\ref{fig:cell_wall_cube_cell_05}(a-b), and by the equivalent membrane stresses (with the maximum being the base material yield stress) in Figures~\ref{fig:cell_wall_cube_cell_05}(c-d), respectively}.

It can be seen from Figures~\ref{fig:cell_wall_cube_cell_05}(a-b) that only the two cell walls parallel with the loading direction accommodate strain energy. These cell walls deform first by the membrane mode ($\mathcal{I}_\stt{w}\sim -1$) and then switch to a mixed membrane-bending mode rapidly after buckling (black triangles). In particular, under compression in the transverse direction, the larger cell wall buckles earlier than the other one. The cell wall buckling points correspond well with the stress drop points in Figure~\ref{fig:eff_response_cube_cell_05}(a).

\textcolor{red}{Figures~\ref{fig:cell_wall_cube_cell_05}(c-d) show that for each loading case, yielding (black circles) happens later than buckling, indicating that the elastic buckling triggers failure. As expected, each buckled cell wall undergoes a stress redistribution. The compressive load becomes mainly carried by the cell wall portion close to cell edges, implying a load-bearing efficiency reduction. In the transverse direction, a lower load-bearing efficiency can be observed for the larger cell wall, compared with the other one}. The results for other shape anisotropy have confirmed similar behaviour to those presented for $\mathcal{R}=1.5$ and are thus not shown here.

\begin{figure}[ht]
\centering
\begin{overpic}[draft=false,width=0.9\textwidth]{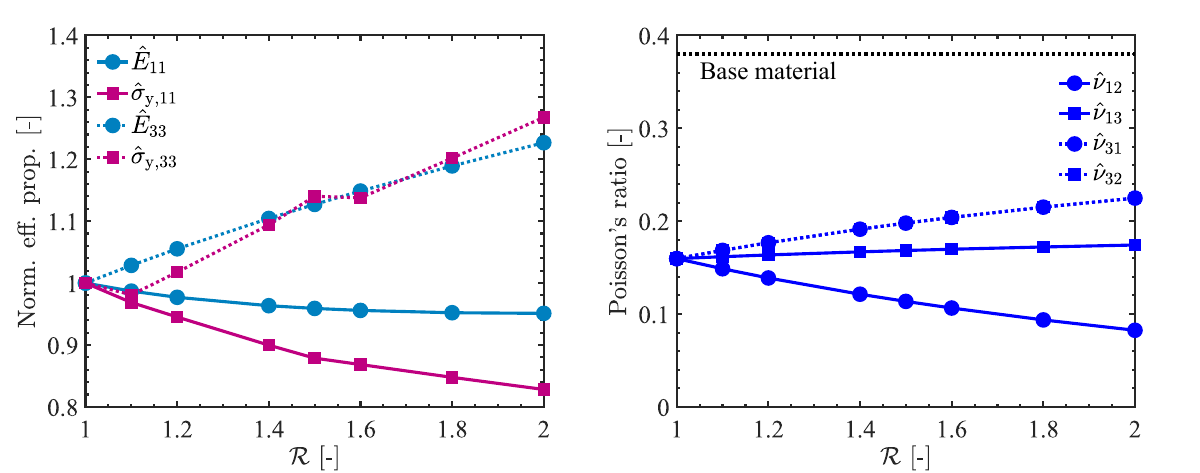}
\put(24.5,37.5){\color{black}\footnotesize \textbf{(a)}}
\put(73.625,37.5){\color{black}\footnotesize \textbf{(b)}}
\end{overpic}
\caption{Effective compressive properties of the rectangular parallelepiped cell-based models with different shape anisotropy $\mathcal{R}$: (a) modulus and yield strength, and (b) Poisson's ratio. Results in (a) have been normalized with respect to those at $\mathcal{R}=1.0$.}\label{fig:eff_props_cube_cell}
\end{figure}

Next, the effective compressive properties are extracted and plotted against varying shape anisotropy $\mathcal{R}$ in Figure~\ref{fig:eff_props_cube_cell}. For each loading direction $\vec{e}_i$, the compressive modulus and Poisson's ratio are evaluated as $\hat{E}_{ii}=\Delta \hat{P}_{ii}/\Delta \hat{F}_{ii}$ and $\hat{\nu}_{ij}=-\Delta \hat{F}_{jj}/\Delta \hat{F}_{ii}$, respectively, in the initial elastic region. The yield strength is determined by $\hat{\sigma}_{\stt{y},ii}=\hat{P}_{ii}$ \textcolor{red}{at the first cell wall yield point (see Figures~\ref{fig:cell_wall_cube_cell_05}(c-d) and also indicated in Figure~\ref{fig:cell_wall_cube_cell_05}(a))}. Regarding elastic properties, Figure~\ref{fig:eff_props_cube_cell}(a) shows that $\hat{E}_{33}$ increases while $\hat{E}_{11}$ decreases with increasing $\mathcal{R}$, indicating an increasing modulus ratio $\hat{E}_{33}/\hat{E}_{11}$. Figure~\ref{fig:eff_props_cube_cell}(b) shows that $\hat{\nu}_{31}=\hat{\nu}_{32}$ for a given $\mathcal{R}$, as expected for a transverse isotropy. $\hat{\nu}_{13}$ and $\hat{\nu}_{31}$ both increase along with $\mathcal{R}$ while $\hat{\nu}_{12}$ decreases. Compared with $\hat{\nu}_{13}$, $\hat{\nu}_{31}$ presents a stronger dependency on $\mathcal{R}$. This can be understood through the well-known relation $\hat{\nu}_{13}/\hat{E}_{11}=\hat{\nu}_{31}/\hat{E}_{33}$.

Back to Figure~\ref{fig:eff_props_cube_cell}(a), the compressive strength $\hat{\sigma}_{\stt{y},33}$ tends to increase \textcolor{red}{while $\hat{\sigma}_{\stt{y},11}$ decreases with increasing $\mathcal{R}$}. The resulting strength ratio $\hat{\sigma}_{\stt{y},33}/\hat{\sigma}_{\stt{y},11}$ implies an overall rising trend. Moreover, a stronger dependency on $\mathcal{R}$ is observed for $\hat{\sigma}_{\stt{y},33}/\hat{\sigma}_{\stt{y},11}$, compared with $\hat{E}_{33}/\hat{E}_{11}$.

\begin{figure}[ht]
\centering
\begin{overpic}[draft=false,width=0.925\textwidth]{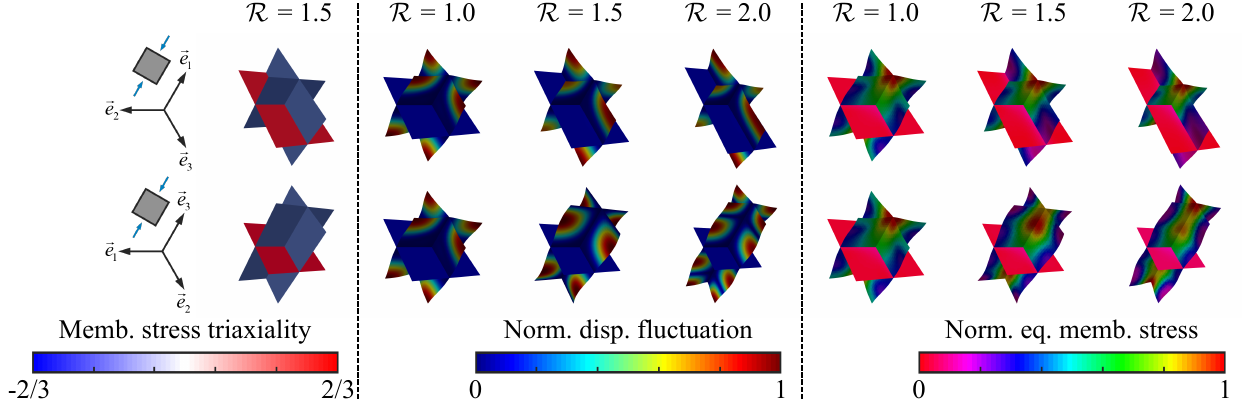}
\put(2.5,23){\color{black}\footnotesize \textbf{(a)}}
\put(2.5,11.625){\color{black}\footnotesize \textbf{(b)}}
\end{overpic}
\caption{\textcolor{red}{Initial membrane stress triaxiality distributions, and buckling patterns and membrane stress patterns at the yield points of the rectangular parallelepiped cell-based models with different shape anisotropy $\mathcal{R}$, under uniaxial compression in the (a) transverse and (b) foam rise directions. The practical triaxiality values range from $-0.48$ to $0.67$ in (a), and from $-0.43$ and $0.67$ in (b), respectively. The loading direction is represented by a pair of opposite arrows}.}\label{fig:stress_triaxiality_buckling_mode_cube_cell}
\end{figure}

\textcolor{red}{In addition, the initial stress states, and buckling patterns and membrane stress patterns upon yielding are analysed}. Since the cell wall membrane deformation dominates the initial elastic region, the membrane stress triaxiality is evaluated to better understand the pre-buckling stress states. The results for several $\mathcal{R}$ are reported in Figure~\ref{fig:stress_triaxiality_buckling_mode_cube_cell}. It is found that the initial stress states of the two cell walls for $\mathcal{R}=1.5$, parallel with the loading direction, are close to uniaxial compression (theoretical triaxiality $-1/3$). Similar stress triaxiality distributions have been confirmed for other shape anisotropy and are thus not presented. 

\textcolor{red}{The buckling patterns and membrane stress patterns are visualized through the deformed configurations colored by the normalized displacement fluctuations and equivalent membrane stresses, respectively. Note, that the present buckling patterns more or less deviate from the theoretical buckling modes since the geometrical nonlinear effect is involved. It can be seen that the buckling pattern and membrane stress pattern generally depend on the loading direction and $\mathcal{R}$. In particular, under compression in the foam rise direction, these patterns become more wavy (from one to three half waves) with increasing $\mathcal{R}$}.

\subsection{Kelvin cell}
\begin{figure}[ht]
\centering
\begin{overpic}[draft=false,width=0.9\textwidth]{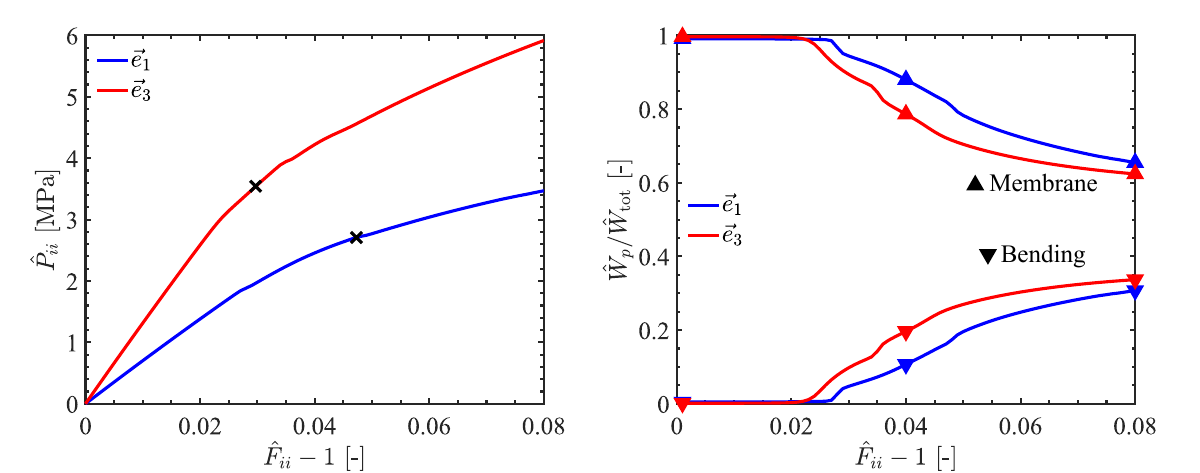}
\put(24.5,37.5){\color{black}\footnotesize \textbf{(a)}}
\put(73.625,37.5){\color{black}\footnotesize \textbf{(b)}}
\end{overpic}
\caption{Effective responses of the Kelvin cell-based model with $\mathcal{R}=1.5$, under uniaxial compression in the transverse ($\vec{e}_1$) and foam rise ($\vec{e}_3$) directions: (a) stress and (b) strain energy fraction versus applied strain.}\label{fig:eff_response_Kelvin_cell_05}
\end{figure}

The macroscale effective responses for shape anisotropy $\mathcal{R}=1.5$ are reported in Figure~\ref{fig:eff_response_Kelvin_cell_05}. For each loading case, the initial elastic region is followed by multiple times of stiffness reduction rather than apparent stress drops observed for the rectangular parallelepiped cell (see Figure~\ref{fig:eff_response_cube_cell_05}(a)). The compressive stress in the foam rise direction is well above that in the transverse direction, indicating an anisotropic compressive behaviour. The strain energy fraction profiles (see Figure~\ref{fig:eff_response_Kelvin_cell_05}(b)) demonstrate that for each loading direction, the membrane deformation mode ($\hat{W}_\stt{m}/\hat{W}_\stt{tot}\sim 1$) governs the initial elastic region and the bending mode becomes important only after the stiffness reduction. The corresponding strain energy redistribution proceeds gradually instead of in a sudden manner observed for the rectangular parallelepiped cell (see Figure~\ref{fig:eff_response_cube_cell_05}(b)).

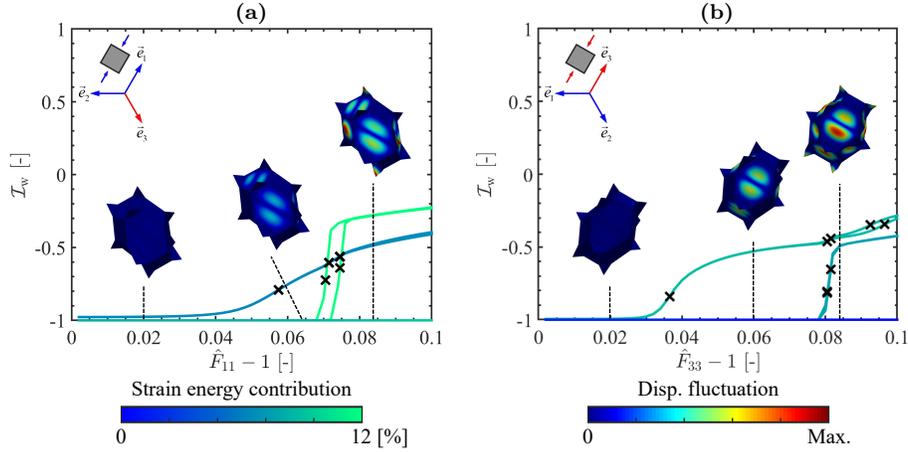
\begin{figure}[ht]
\centering
\begin{overpic}[draft=false,width=0.9\textwidth]{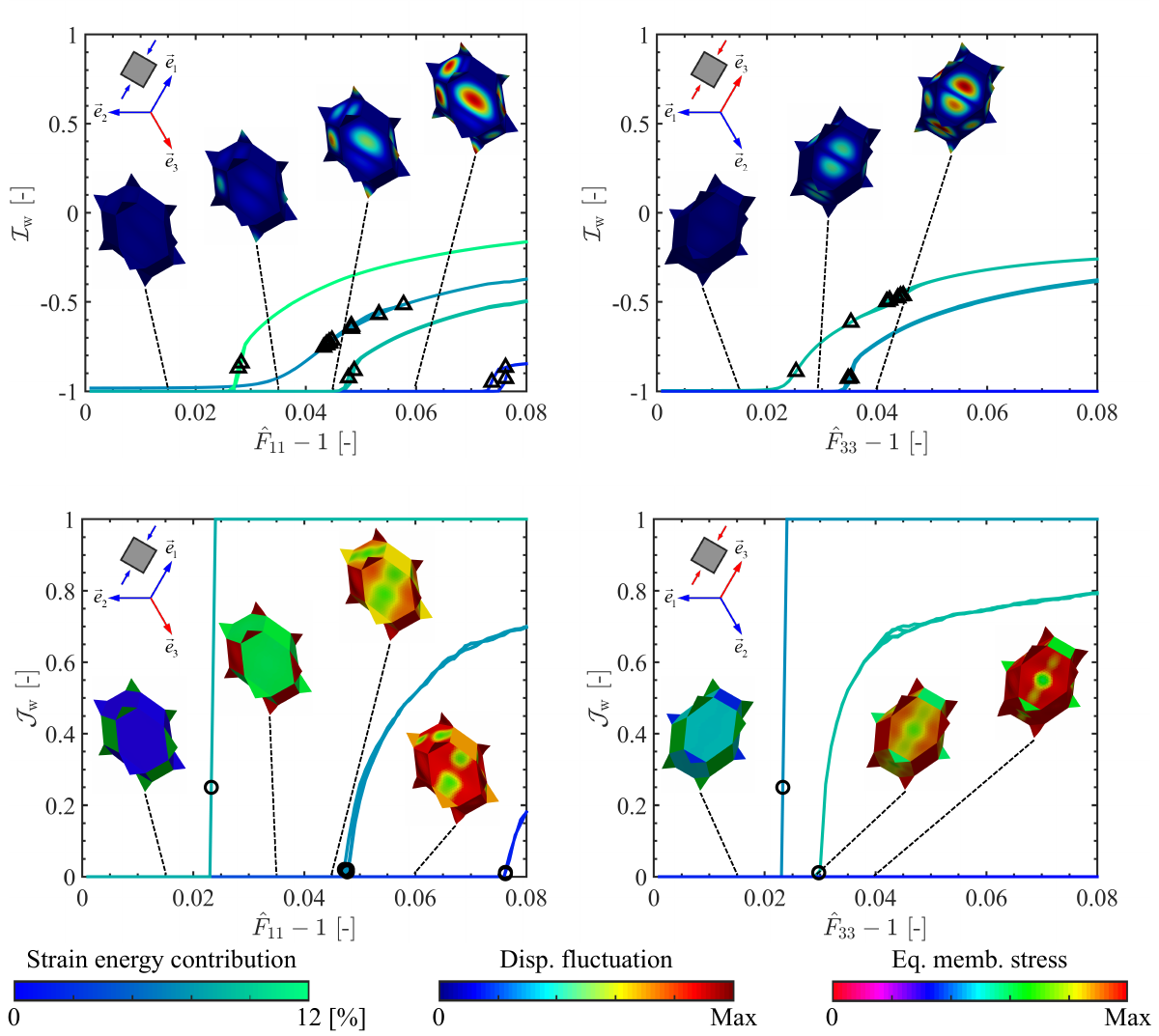}
\put(24.5,87.625){\color{black}\footnotesize \textbf{(a)}}
\put(73.625,87.625){\color{black}\footnotesize \textbf{(b)}}
\put(24.5,45.75){\color{black}\footnotesize \textbf{(c)}}
\put(73.625,45.75){\color{black}\footnotesize \textbf{(d)}}
\end{overpic}
\caption{\textcolor{red}{Cell wall strain energy partitioning indicators, membrane plasticity indicators, and deformed configurations at different stages of the Kelvin cell-based model with $\mathcal{R}=1.5$, under uniaxial compression in the (a, c) transverse and (b, d) foam rise directions.}}\label{fig:cell_wall_Kelvin_cell_05}
\end{figure}

\textcolor{red}{The cell wall-wise strain energy partitioning indicators and membrane plasticity indicators are reported in Figure~\ref{fig:cell_wall_Kelvin_cell_05}}. It can be seen from Figures~\ref{fig:cell_wall_Kelvin_cell_05}(a-b) that the strain energy is mainly accommodated by the eight hexagonal cell walls inclined about the loading direction and the three quadrilateral ones parallel with the loading direction. These cell walls deform first by the membrane mode ($\mathcal{I}_\stt{w}\sim -1$), followed by transition towards a mixed membrane-bending mode after buckling (black triangles). The inclined cell walls present a gradual deformation mode transition in contrast to the parallel ones. This can be explained by that the inclination angle increases with increasing applied strain, leading to a reduction of the load portion projected in the cell wall plane. Besides, the large, inclined cell walls tend to buckle earlier than the small, parallel ones, especially under compression in the foam rise direction. The cell wall buckling points correspond well with the stiffness reduction points in Figure~\ref{fig:eff_response_Kelvin_cell_05}(a).

\textcolor{red}{Figures~\ref{fig:cell_wall_Kelvin_cell_05}(c-d) show that for each loading case, the inclined cell walls yield (black circles) after buckling, while the parallel ones yield already before buckling and thus would fail by plastic collapse instead. Focusing on the inclined cell walls which carry the most compressive load, a similar stress redistribution as observed in Figures~\ref{fig:cell_wall_cube_cell_05}(c-d) is confirmed. However, the redistributed stresses in these cell walls are much more uniform compared with those for the rectangular parallelepiped cell (see Figure~\ref{fig:cell_wall_cube_cell_05}(c-d)), implying a higher load-bearing efficiency for the Kelvin cell (in the post-buckling regime)}. Similar observations have been confirmed for other shape anisotropy.

\begin{figure}[ht]
\centering
\begin{overpic}[draft=false,width=0.9\textwidth]{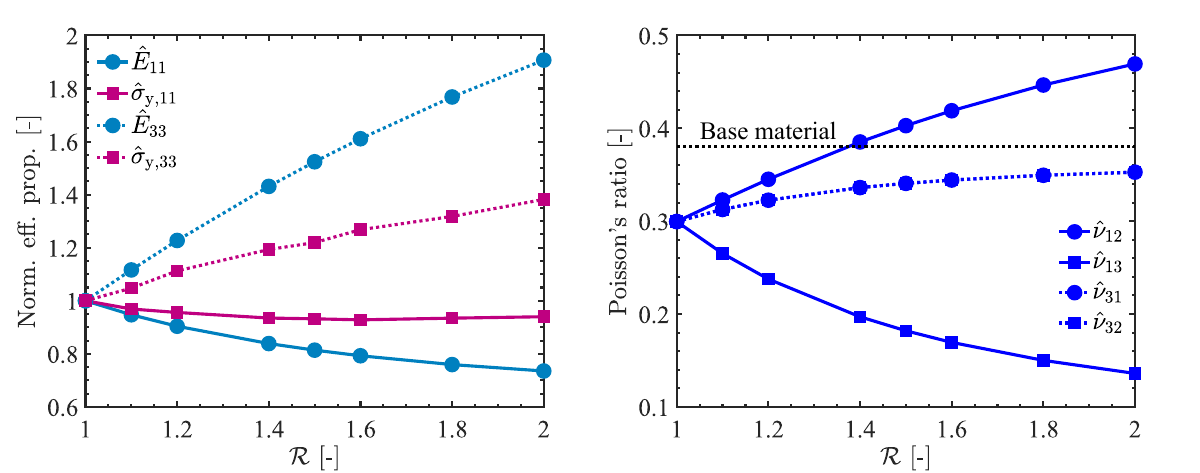}
\put(24.5,37.5){\color{black}\footnotesize \textbf{(a)}}
\put(73.625,37.5){\color{black}\footnotesize \textbf{(b)}}
\end{overpic}
\caption{Effective compressive properties of the Kelvin cell-based models with different shape anisotropy $\mathcal{R}$: (a) modulus and yield strength, and (b) Poisson's ratio.}\label{fig:eff_props_Kelvin_cell}
\end{figure}

The effective compressive properties for different shape anisotropy $\mathcal{R}$ are reported in Figure~\ref{fig:eff_props_Kelvin_cell}. \textcolor{red}{Here, the yield strength is determined at the first inclined cell wall yield point (see Figures~\ref{fig:cell_wall_Kelvin_cell_05}(c-d) and also indicated in Figure~\ref{fig:cell_wall_Kelvin_cell_05}(a)))}. As $\mathcal{R}$ increases, the compressive modulus $\hat{E}_{33}$ increases while $\hat{E}_{11}$ decreases (see Figure~\ref{fig:eff_props_Kelvin_cell}(a)), implying a rapid increase of the modulus ratio $\hat{E}_{33}/\hat{E}_{11}$. The compressive Poisson's ratios $\hat{\nu}_{12}$ and $\hat{\nu}_{31}$ (identical to $\hat{\nu}_{32}$) both increase along with $\mathcal{R}$ while $\hat{\nu}_{13}$ decreases (see Figure~\ref{fig:eff_props_Kelvin_cell}(b)).

The compressive strength $\hat{\sigma}_{\stt{y},33}$ increases along with $\mathcal{R}$, while $\hat{\sigma}_{\stt{y},11}$ remains almost constant, leading to a slow increase of the strength ratio $\hat{\sigma}_{\stt{y},33}/\hat{\sigma}_{\stt{y},11}$ (see Figure~\ref{fig:eff_props_Kelvin_cell}(a)). Compared with $\hat{E}_{33}/\hat{E}_{11}$, $\hat{\sigma}_{\stt{y},33}/\hat{\sigma}_{\stt{y},11}$ is much less sensitive to $\mathcal{R}$. 

\begin{figure}[ht]
\centering
\begin{overpic}[draft=false,width=0.925\textwidth]{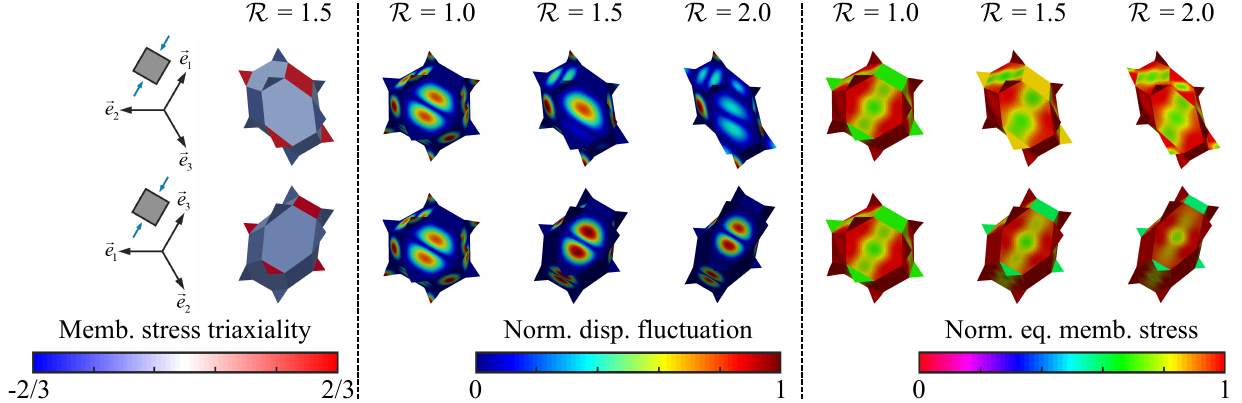}
\put(2.5,23){\color{black}\footnotesize \textbf{(a)}}
\put(2.5,11.625){\color{black}\footnotesize \textbf{(b)}}
\end{overpic}
\caption{\textcolor{red}{Initial membrane stress triaxiality distributions, and buckling patterns and membrane stress patterns at the yield points of the Kelvin cell-based models with different shape anisotropy $\mathcal{R}$, under uniaxial compression in the (a) transverse and (b) foam rise directions. The practical triaxiality values range from $-0.44$ to $0.64$ in (a), and from $-0.36$ and $0.67$ in (b), respectively}.}\label{fig:stress_triaxiality_buckling_mode_Kelvin_cell}
\end{figure}

\textcolor{red}{The initial membrane stress triaxiality distributions, buckling patterns and membrane stress patterns upon yielding for several $\mathcal{R}$ are reported in Figure~\ref{fig:stress_triaxiality_buckling_mode_Kelvin_cell}}. Nearly uniaxial compression stress states can be observed on the eight inclined and three parallel cell walls with respect to the loading direction. Similar stress states have also been confirmed for Kelvin cells with other shape anisotropy. 

For each $\mathcal{R}$, the bucking pattern and membrane stress pattern are dominated by those large cell walls despite being inclined, especially under compression in the foam rise direction.  In general, the bucking pattern and membrane stress pattern depend on the loading direction and $\mathcal{R}$. \textcolor{red}{It seems that these patterns for $\mathcal{R}=1.5$ in the transverse direction (see Figure~\ref{fig:stress_triaxiality_buckling_mode_Kelvin_cell}(a)) has not yet fully developed when yielding happens (see Figures~\ref{fig:cell_wall_Kelvin_cell_05}(a) and (c))}.

\subsection{Discussion}
Despite numerous simplifications, the idealized cell-based models have shown capabilities to qualitatively reproduce the anisotropic compressive behaviour of realistic foams. Through quantitative analysis of the cell wall deformation behaviour, a few preliminary insights into the deformation mechanisms are summarized:
\begin{enumerate}
\item The initial elastic region is primarily governed by the cell wall membrane deformation, regardless of the loading direction. 
\item The cell wall bending contribution becomes crucial only \textcolor{red}{after buckling, followed by membrane yielding}.
\end{enumerate}
Notice that the first insight has been raised to some extent elsewhere (see e.g.\ \cite{Simone1998a, Grenestedt2000, Shi2018}), however, limited to isotropic foams and lacking confirmation of the cell wall deformation. It will be shown in \textcolor{blue}{Section}~\ref{sec:analy_tess_meso_model} that these insights hold even when the mesostructural stochastics are taken into account. Accordingly, the anisotropic compressive properties of the foams with a high cell face fraction, as considered in this study, may not be simply explained by different deformation mechanisms for different loading directions, in contrast to those with a low cell face fraction (see e.g.\ \cite{MarviMashhadi2018, MarviMashhadi2020, Ding2023}). 

The Kelvin cells exhibit noticeably different anisotropy trends of compressive properties from the rectangular parallelepiped cells. Extensive experimental observations on Divinycell foams (see e.g.\ \cite{Liu2020, Tang2022, Zhou2023b}) have confirmed that, compared with the modulus anisotropy, the strength anisotropy is much less sensitive to shape anisotropy. This can be captured by the Kelvin cell but not by the rectangular parallelepiped cell. All the above suggest that the cell wall inclination angle has non-negligible impacts on the mechanical anisotropy of realistic foams, in addition to the primary role of cell shape anisotropy.

\textcolor{red}{In addition, it has been confirmed that the compressive strength anisotropy only varies $\sim \pm10\%$ for a broad range of relative densities (see~\ref{sec:influence_plas}), with respect to the chosen relative density 0.075. The modulus anisotropy is almost insensitive to the relative density. These results imply that the present findings are general for low-density foams}.

\section{Relationships between mechanical anisotropy and cell shape anisotropy}\label{sec:relation_anisotropy}
The insights obtained in \textcolor{blue}{Section}~\ref{sec:analy_ideal_meso_model} will guide derivations of analytical models in this section, to illustrate how cell shape anisotropy translates into mechanical anisotropy.

\subsection{Model development}
\subsubsection{Rectangular parallelepiped cell}
\begin{figure}[ht]
\centering
\begin{overpic}[draft=false,width=0.58\textwidth]{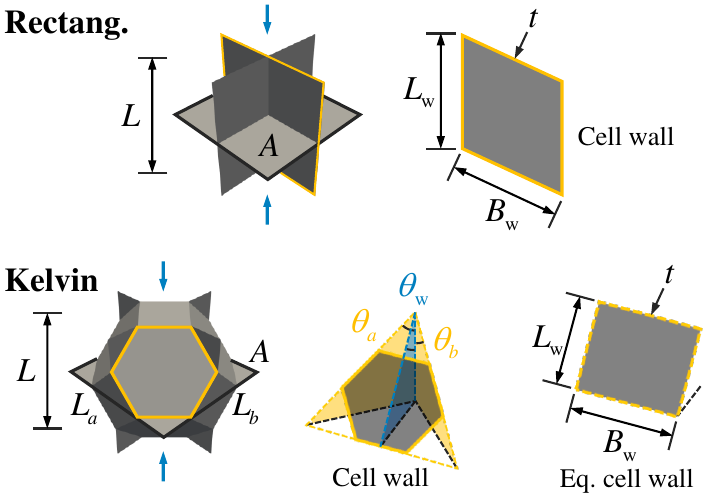}
\put(-6,66.375){\color{black}\footnotesize \textbf{(a)}}
\put(-6,30.25){\color{black}\footnotesize \textbf{(b)}}
\end{overpic}
\caption{Sketches of the two idealized cell structures and load-bearing cell walls under compression: (a) rectangular parallelepiped and (b) Kelvin cells. The loading direction is represented by a pair of opposite arrows.}\label{fig:geo_model_cell_wall}
\end{figure}

In a rectangular parallelepiped cell structure, the cell walls parallel with the loading direction constitute the primary load-bearing elements (see Figure~\ref{fig:cell_wall_cube_cell_05}), as sketched in Figure~\ref{fig:geo_model_cell_wall}(a). The cell length and cross-section area in the loading direction are indicated by $L$ and $A$, respectively. For instance, under compression in the foam rise direction $\vec{e}_3$ (see Figure~\ref{fig:geo_model_configuration}(a)), $L=L_3$ and $A=L_1L_2$. The cell wall length spans over the entire cell, i.e.\ $L_\stt{w}=L$, and the width is indicated by $B_\stt{w}$. 

The effective compressive modulus $\hat{E}$ can be expressed as:
\begin{equation}\label{eq:cube_cell_modulus}
\hat{E}=\frac{L}{A}\sum^N K_\stt{w},
\end{equation}
where $N$ is the number of load-bearing cell walls, and $K_\stt{w}$ is the cell wall-wise membrane stiffness given by
\begin{equation}\label{eq:cube_cell_wall_stiffness}
K_\stt{w}= E_\stt{w}\frac{t}{L}B_\stt{w},
\end{equation}
with $E_\stt{w}$ being the cell wall membrane modulus which accounts for the stress state effect. Assuming a uniaxial compression stress state (see Figure~\ref{fig:stress_triaxiality_buckling_mode_cube_cell}), $E_\stt{w}$ can be replaced by the material Young's modulus $E$.

The effective compressive strength in the loading direction can be expressed as:\color{red}
\begin{equation}\label{eq:cube_cell_yield_strength}
\hat{\sigma}_\stt{y}=\min\{\sigma_\stt{y,w}\}\frac{t}{A}\sum^N B_\stt{w},
\end{equation}
with $\sigma_\stt{y,w}$ being the cell wall-wise compressive strength. Note, that $\hat{\sigma}_\stt{y}$ is determined by the weakest cell wall, i.e.\ with the lowest aspect ratio $\mathcal{R}_\stt{w}=L_\stt{w}/B_\stt{w}$ (see Figure~\ref{fig:stress_triaxiality_buckling_mode_cube_cell}). Based on the effective width principle, $\sigma_\stt{y,w}$ of a rectangular plate can be approximated as $\sigma_\stt{y,w}=\sqrt{\sigma_\stt{c,w}\sigma_\stt{y}}$ \cite{Timoshenko1961}, with $\sigma_\stt{c,w}$ and $\sigma_\stt{y}$ being the plate buckling strength and material yield stress, respectively. Following the linear buckling theory \cite{Gerard1957}, $\sigma_\stt{c,w}$ is given by 
\color{black}
\begin{equation}\label{eq:cube_cell_wall_buckling_stress}
\sigma_\stt{c,w}\propto k_\stt{c}\left(\frac{t}{B_\stt{w}}\right)^2,
\end{equation}
where ``$\propto$" represents a proportional relationship; $k_\stt{c}$ is the buckling coefficient which depends on $\mathcal{R}_\stt{w}$. Although the theoretical solution of $k_\stt{c}$ is available only under certain boundary conditions, the functional form below provides a sufficient approximation in most cases (see also\ \cite{Gerard1957})
\begin{equation}\label{eq:plate_buckling_coefficient}
k_\stt{c}\propto \mathcal{K}_\stt{c}(\mathcal{R}_\stt{w})=1-k+k(\mathcal{R}_\stt{w})^p,
\end{equation}
where a normalized buckling coefficient function $\mathcal{K}_\stt{c}$ with $\mathcal{K}_\stt{c}(1)=1$ has been introduced; $k$ and $p$ are the two parameters which can be identified using numerical results.

Using the effective compressive properties $(\hat{E}_\stt{11},\hat{E}_\stt{33})$ and $(\hat{\sigma}_\stt{y,11},\hat{\sigma}_\stt{y,33})$ in the transverse ($\vec{e}_1$) and foam rise ($\vec{e}_3$) directions, which are given by eqs.~\eqref{eq:cube_cell_modulus} and~\eqref{eq:cube_cell_yield_strength}, the modulus anisotropy $\mathcal{R}^E$ and strength anisotropy $\mathcal{R}^\sigma$ for a rectangular parallelepiped cell, can be expressed as:
\begin{subequations}\label{eq:cube_cell_compressive_anisotropy}
	\begin{align}
	&\mathcal{R}^E=\frac{\hat{E}_\stt{33}}{\hat{E}_\stt{11}}=\mathcal{R}_\stt{f},\\
	&\textcolor{red}{\mathcal{R}^\sigma=\frac{\hat{\sigma}_\stt{y,33}}{\hat{\sigma}_\stt{y,11}}=\mathcal{R}_\stt{c}^{\frac{1}{2}}\mathcal{R}_\stt{f}},
	\end{align}
\end{subequations}
where $\mathcal{R}_\stt{f}$ denotes a cell load-bearing area fraction ratio and $\mathcal{R}_\stt{c}$ a cell wall buckling strength ratio, defined through
\begin{equation}\label{eq:cube_cell_ratios}
\mathcal{R}_\stt{f}=\frac{\{f_\stt{w}\}_3}{\{f_\stt{w}\}_1},\,\,\mathcal{R}_\stt{c}=\frac{\min\{\sigma_\stt{c,w}\}_3}{\min\{\sigma_\stt{c,w}\}_1},
\end{equation}
with $f_\stt{w}=\frac{t}{A}\sum B_\stt{w}$. \textcolor{red}{Note, that $\sigma_\stt{y,w}=\sqrt{\sigma_\stt{c,w}\sigma_\stt{y}}$ has been substituted in the derivation of $\mathcal{R}^\sigma$ and $\sigma_\stt{y}$ gets eliminated eventually}. It can be observed that $\mathcal{R}^E$ is purely determined by $\mathcal{R}_\stt{f}$, while $\mathcal{R}^\sigma$ depends on both $\mathcal{R}_\stt{c}$ and $\mathcal{R}_\stt{f}$.

Based on the geometrical relationships shown in Figure~\ref{fig:geo_model_cell_wall}(a) that allow for expressing $L$, $A$ and $B_\stt{w}$ in terms of $L_i$, followed by substituting $L_1=L_2=V^\frac{1}{3}\mathcal{R}^{\stt{-}\frac{1}{3}}$ and $L_3=V^\frac{1}{3}\mathcal{R}^{\frac{2}{3}}$ (see \textcolor{blue}{Section}~\ref{sec:num_model_method}) into eq.~\eqref{eq:cube_cell_ratios}, the two ratios $\mathcal{R}_\stt{f}$ and $\mathcal{R}_\stt{c}$ can be directly related to shape anisotropy $\mathcal{R}$ through:
\begin{equation}\label{eq:cube_cell_ratios_2}
\mathcal{R}_\stt{f}=\frac{2\mathcal{R}}{1+\mathcal{R}},\,\,\mathcal{R}_\stt{c}=\frac{\mathcal{K}_\stt{c}(\mathcal{R})}{\mathcal{K}_\stt{c}(\mathcal{R}^\stt{-1})}\mathcal{R}^2.
\end{equation}
Combining eqs.~\eqref{eq:cube_cell_compressive_anisotropy} and~\eqref{eq:cube_cell_ratios_2}, the mechanical anisotropy for a rectangular parallelepiped cell can be fully predicted with shape anisotropy as the input.

\subsubsection{Kelvin cell}
In a Kelvin cell structure, the cell walls inclined about the loading direction may be regarded as the primary load-bearing elements (see Figure~\ref{fig:cell_wall_Kelvin_cell_05}), as sketched in Figure~\ref{fig:geo_model_cell_wall}(b). The cell wall inclination angle is indicated by $\theta_\stt{w}$, which can be related to the two cell edge inclination angles $\theta_a$ and $\theta_b$ through
\begin{equation}\label{eq:cell_wall_inclination_angle}
\frac{1}{\tan^2\theta_\stt{w}}=\frac{1}{\tan^2\theta_a}+\frac{1}{\tan^2\theta_b},
\end{equation}
with $\theta_a$ and $\theta_b$ given by
\begin{equation}\label{eq:cell_edge_inclination_angles}
\tan\theta_a=\frac{L_a}{L},\,\,\tan\theta_b=\frac{L_b}{L},
\end{equation}
where $L_a$ and $L_b$ are the two cell global dimensions on the plane perpendicular to the loading direction. For instance, under compression in the transverse direction $\vec{e}_1$ (see Figure~\ref{fig:geo_model_configuration}(b)), $L_a=L_2$ and $L_b=L_3$. To continue the analytical model derivations, the hexagonal cell wall is then approximated using an equivalent rectangular cell wall with its length $L_\stt{w}=\frac{1}{2}L/\cos\theta_\stt{w}$ and width $B_\stt{w}=\frac{1}{2}\kappa_\stt{w}\sqrt{L_a^2+L_b^2}$; $\kappa_\stt{w}=\frac{\sqrt{3}}{2}$ is a correction factor such that the equivalent cell wall aspect ratio $\mathcal{R}_\stt{w}=L_\stt{w}/B_\stt{w}=1$ at $\mathcal{R}=1.0$. 

Taking into account $\theta_\stt{w}$, the effective compressive modulus can be expressed by modifying eq.~\eqref{eq:cube_cell_modulus} as:
\begin{equation}\label{eq:Kelvin_cell_modulus}
\hat{E}=\cos^2\theta_\stt{w}\frac{L}{A}\sum^{N/2} \frac{K_\stt{w}}{2},
\end{equation}
where $\cos^2\theta_\stt{w}$ quantifies the membrane stiffness portion of an inclined wall in the loading direction (see also \cite{Chai2020}), and the equivalent cell wall membrane stiffness $K_\stt{w}$ is given by substituting $L_\stt{w}=\frac{1}{2}L/\cos\theta_\stt{w}$ into eq.~\eqref{eq:cube_cell_wall_stiffness}
\begin{equation}\label{eq:Kelvin_cell_wall_stiffness}
K_\stt{w}= 2\cos\theta_\stt{w}E_\stt{w}\frac{t}{L}B_\stt{w}.
\end{equation}

Similarly, the effective compressive strength is expressed by modifying eq.~\eqref{eq:cube_cell_yield_strength} as:\color{red}
\begin{equation}\label{eq:Kelvin_cell_yield_strength}
\hat{\sigma}_\stt{y}=\cos\theta_\stt{w}\min\{\sigma_\stt{y,w}\}\frac{t}{A}\sum^N B_\stt{w},
\end{equation}\color{black}
where $\cos\theta_\stt{w}$ quantifies the membrane stress portion of an inclined wall in the loading direction (see also \cite{Chai2020}), \textcolor{red}{and $\sigma_\stt{y,w}=\sqrt{\sigma_\stt{c,w}\sigma_\stt{y}}$ with $\sigma_\stt{c,w}$ already given by eq.~\eqref{eq:cube_cell_wall_buckling_stress}}.

Using the effective compressive properties $(\hat{E}_\stt{11},\hat{E}_\stt{33})$ and $(\hat{\sigma}_\stt{y,11},\hat{\sigma}_\stt{y,33})$ in the transverse ($\vec{e}_1$) and foam rise ($\vec{e}_3$) directions, which are given by eqs.~\eqref{eq:Kelvin_cell_modulus} and~\eqref{eq:Kelvin_cell_yield_strength}, $\mathcal{R}^E$ and $\mathcal{R}^\sigma$ for a Kelvin cell, can be expressed as:
\begin{subequations}\label{eq:Kelvin_cell_compressive_anisotropy}
	\begin{align}
	&\mathcal{R}^E=\frac{\hat{E}_\stt{33}}{\hat{E}_\stt{11}}=\mathcal{R}_{\theta}^3\mathcal{R}_\stt{f},\\
	&\textcolor{red}{\mathcal{R}^\sigma=\frac{\hat{\sigma}_\stt{y,33}}{\hat{\sigma}_\stt{y,11}}=\mathcal{R}_{\theta}\mathcal{R}_\stt{c}^{\frac{1}{2}}\mathcal{R}_\stt{f}},
	\end{align}
\end{subequations}
where $\mathcal{R}_\stt{f}$ and $\mathcal{R}_\stt{c}$ are already defined in eq.~\eqref{eq:cube_cell_ratios}; the third ratio $\mathcal{R}_{\theta}$ has been introduced
\begin{equation}\label{eq:Kelvin_cell_ratio_theta}
\mathcal{R}_\theta=\frac{\cos\{\theta_\stt{w}\}_3}{\cos\{\theta_\stt{w}\}_1}.
\end{equation}
Here, $\{\theta_\stt{w}\}_i$ is the cell wall inclination angle with respect to the global direction $\vec{e}_i$ (see eq.~\eqref{eq:cell_wall_inclination_angle}). It can be seen that for the Kelvin cell, both $\mathcal{R}^E$ and $\mathcal{R}^\sigma$ are additionally influenced by $\theta_\stt{w}$ through $\mathcal{R}_{\theta}$, compared with the rectangular parallelepiped cell (see eq.~\eqref{eq:cube_cell_compressive_anisotropy}).

Expressing $\theta_\stt{w}$, $L$, $A$ and $B_\stt{w}$ in terms of $L_i$ according to the geometrical relationships shown in Figure~\ref{fig:geo_model_cell_wall}(b), in combination with $L_1=L_2=V^\frac{1}{3}\mathcal{R}^{\stt{-}\frac{1}{3}}$ and $L_3=V^\frac{1}{3}\mathcal{R}^{\frac{2}{3}}$ (see also \textcolor{blue}{Section}~\ref{sec:num_model_method}), the three ratios $\mathcal{R}_\theta$, $\mathcal{R}_\stt{f}$ and $\mathcal{R}_\stt{c}$ become:
\begin{equation}\label{eq:Kelvin_cell_ratios_2}
\mathcal{R}_\theta=\frac{\sqrt{2}\mathcal{R}}{(1+\mathcal{R}^2)^{\frac{1}{2}}},\,\,
\mathcal{R}_\stt{f}=\frac{\sqrt{2}\mathcal{R}}{(1+\mathcal{R}^2)^{\frac{1}{2}}},\,\,\mathcal{R}_\stt{c}=\frac{\mathcal{K}_\stt{c}\left(\tfrac{1}{\sqrt{3}}(1+2\mathcal{R}^2)^{\frac{1}{2}}\right)}{\mathcal{K}_\stt{c}\left(\tfrac{2}{\sqrt{3}}\frac{(1+2\mathcal{R}^2)^{\frac{1}{2}}}{1+\mathcal{R}^2}\right)}\frac{1+\mathcal{R}^2}{2}.
\end{equation}
Notice that $\mathcal{R}_\theta$ and $\mathcal{R}_\stt{f}$ essentially represent different factors despite the same expression. Combining eqs.~\eqref{eq:Kelvin_cell_compressive_anisotropy} and~\eqref{eq:Kelvin_cell_ratios_2}, the mechanical anisotropy for a Kelvin cell can be fully predicted with shape anisotropy as the input.

The above expressions~(\ref{eq:cube_cell_compressive_anisotropy}a) and~(\ref{eq:Kelvin_cell_compressive_anisotropy}a) derived for the compressive modulus anisotropy are parameter-free. The use of the compressive strength anisotropy expressions~(\ref{eq:cube_cell_compressive_anisotropy}b) and~(\ref{eq:Kelvin_cell_compressive_anisotropy}b) requires parameter identification for the normalized buckling coefficient function~\eqref{eq:plate_buckling_coefficient}, which depends on the specific geometry and boundary conditions of load-bearing cell walls \textcolor{red}{and is also influenced by the geometrical nonlinear effect}.

\subsection{Model assessment}
\subsubsection{Rectangular parallelepiped cell}
\begin{figure}[ht]
\centering
\begin{overpic}[draft=false,width=0.9\textwidth]{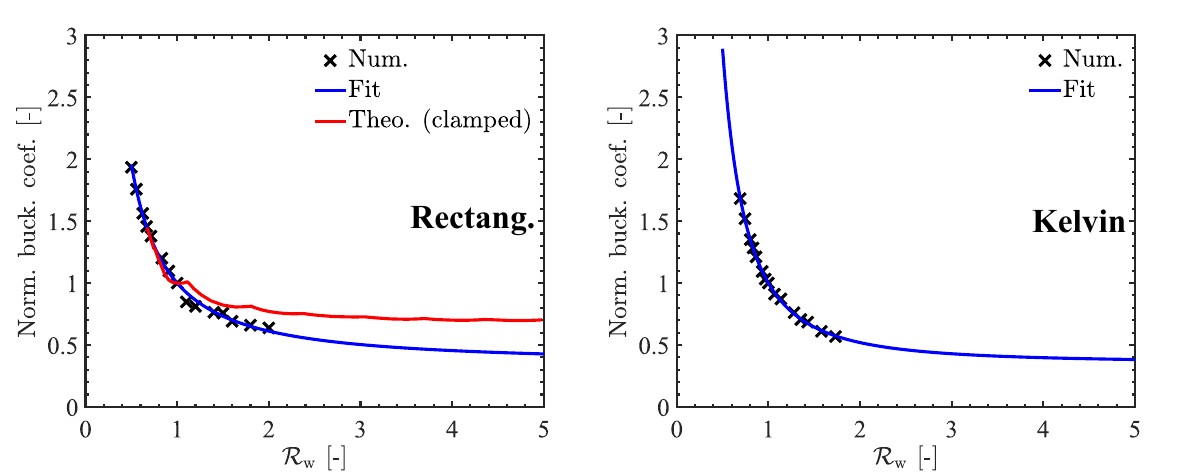}
\put(24.5,37.5){\color{black}\footnotesize \textbf{(a)}}
\put(73.625,37.5){\color{black}\footnotesize \textbf{(b)}}
\end{overpic}
\caption{Normalized buckling coefficients versus aspect ratios $R_\stt{w}$ of the weakest cell walls of the two idealized cell-based models: (a) rectangular parallelepiped and (b) Kelvin cells. Comparison between the numerical data, model fit and theoretical solution is shown.}\label{fig:norm_BK_coef}
\end{figure}

To identify the parameters $k$ and $p$ for the function $\mathcal{K}_\stt{c}(R_\stt{w})$ (see eq.~\eqref{eq:plate_buckling_coefficient}), the normalized buckling coefficients of the load-bearing cell walls are extracted (using eqs.~\eqref{eq:cube_cell_yield_strength} and~\eqref{eq:cube_cell_wall_buckling_stress}) from numerical simulations for all considered shape anisotropy values $\mathcal{R}$ and loading cases, and plotted in Figure~\ref{fig:norm_BK_coef}(a) against the aspect ratio $R_\stt{w}$. The model fit with \textcolor{red}{$k=0.6525$ and $p=-1.3033$} is shown in Figure~\ref{fig:norm_BK_coef}(a), which well reproduces all the numerical data. \textcolor{red}{The numerical data almost overlaps with the theoretical solution of a rectangular plate with fully clamped boundary conditions \cite{Gerard1957} when $R_\stt{w}<1$, while apparent deviations appear when $R_\stt{w}>1$. This may be attributed to the pronounced geometrical nonlinear effect for a large $R_\stt{w}$ (accompanied by a more wavy buckling pattern and see Figure~\ref{fig:stress_triaxiality_buckling_mode_cube_cell}), which is not accounted for in the theoretical solution}.

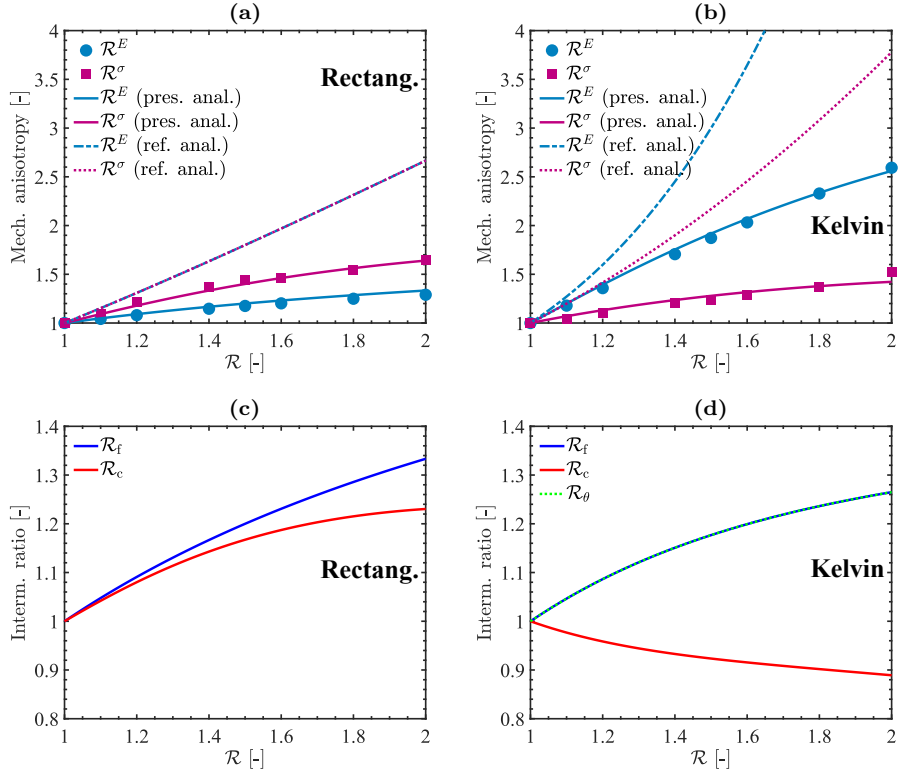
\begin{figure}[ht]
\centering
\begin{overpic}[draft=false,width=0.9\textwidth]{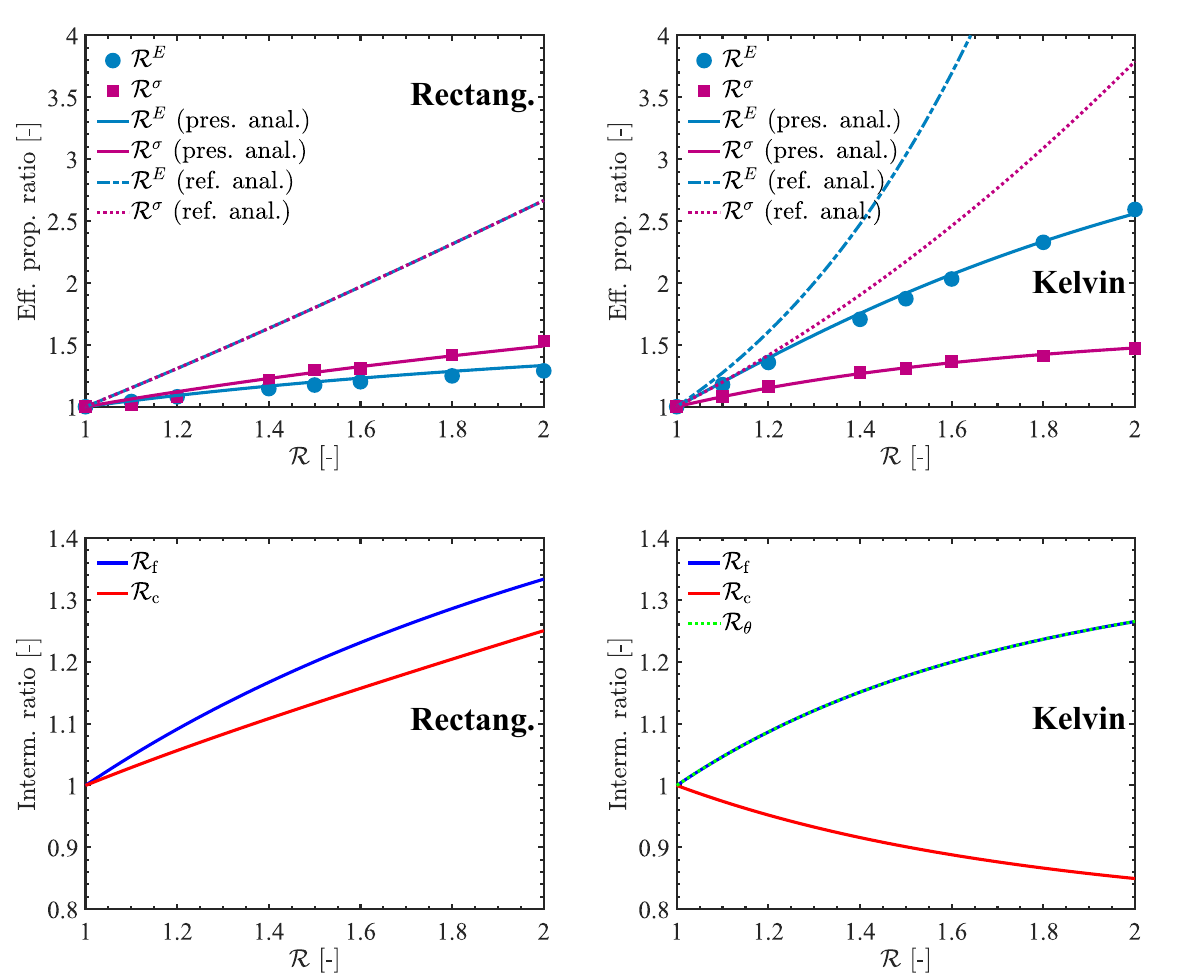}
\put(24.5,79.5){\color{black}\footnotesize \textbf{(a)}}
\put(73.625,79.5){\color{black}\footnotesize \textbf{(b)}}
\put(24.5,37.625){\color{black}\footnotesize \textbf{(c)}}
\put(73.625,37.625){\color{black}\footnotesize \textbf{(d)}}
\end{overpic}
\caption{(a-b) Mechanical anisotropy and (c-d) intermediate ratios versus shape anisotropy $\mathcal{R}$ of the two idealized cell-based models: (a, c) rectangular parallelepiped and (b, d) Kelvin cells. Comparison between the numerical data, and predictions by the present (\ref{eq:cube_cell_compressive_anisotropy},~\ref{eq:Kelvin_cell_compressive_anisotropy}) and reference \cite{Gibson1997_ch6,Sullivan2008} analytical models  is shown in (a-b).}\label{fig:eff_prop_ratios_cell}
\end{figure}

To demonstrate the predictive capabilities of the present analytical model, the mechanical anisotropy $\mathcal{R}^E$ and $\mathcal{R}^\sigma$ computed using eq.~\eqref{eq:cube_cell_compressive_anisotropy}, are plotted in Figure~\ref{fig:eff_prop_ratios_cell}(a) as functions of $\mathcal{R}$ and compared to those from numerical simulations. An excellent agreement can be observed between the numerical data and analytical model predictions. $\mathcal{R}^E$ and $\mathcal{R}^\sigma$ both increase along with $\mathcal{R}$, and $\mathcal{R}^\sigma$ develops faster than $\mathcal{R}^E$ (see also Figure~\ref{fig:eff_props_cube_cell}(a)). 

Moreover, the cell load-bearing area fraction ratio $\mathcal{R}_\stt{f}$ and cell wall buckling strength ratio $\mathcal{R}_\stt{c}$ computed using eq.~\eqref{eq:cube_cell_ratios_2}, are reported in Figure~\ref{fig:eff_prop_ratios_cell}(c). It can be seen that $\mathcal{R}_\stt{f}$ and $\mathcal{R}_\stt{c}$ both increase along with $\mathcal{R}$. According to eq.~\eqref{eq:cube_cell_compressive_anisotropy}, $\mathcal{R}^\sigma$ additionally depends on $\mathcal{R}_\stt{c}$ compared with $\mathcal{R}^E$, thus resulting in faster development of $\mathcal{R}^\sigma$ than $\mathcal{R}^E$, as reflected in Figure~\ref{fig:eff_prop_ratios_cell}(a).

In addition, the widely used Gibson-Ashby model \cite{Gibson1997_ch6} is assessed, which is also derived by assuming a rectangular parallelepiped cell structure. Detailed expressions of this reference analytical model can be found in~\ref{sec:ref_anal_models} and the corresponding model predictions are reported in Figure~\ref{fig:eff_prop_ratios_cell}(a). It can be seen that the Gibson-Ashby model overestimates the mechanical anisotropy, especially modulus anisotropy, \textcolor{red}{by $>200\%$}. Since the same geometrical assumptions have been adopted, these deviations can only be associated with the introduced mechanistic assumptions, i.e.\ load-bearing cell walls are subjected to a tensile stress state under compression and plastic collapse dominates the compressive failure, that are likely inappropriate at a high cell face fraction and low relative density.

\subsubsection{Kelvin cell}
To identify the parameters for the function $\mathcal{K}_\stt{c}(R_\stt{w})$, the normalized buckling coefficients extracted (using eqs.~\eqref{eq:Kelvin_cell_yield_strength} and~\eqref{eq:cube_cell_wall_buckling_stress}) from numerical simulations, are reported in Figure~\ref{fig:norm_BK_coef}(b). The model fit with \textcolor{red}{$k=0.6443 $ and $p=-1.9771$} for the function $\mathcal{K}_\stt{c}(R_\stt{w})$ is shown in Figure~\ref{fig:norm_BK_coef}(b), which again accurately reproduces all the numerical data.

The mechanical anisotropy $\mathcal{R}^E$ and $\mathcal{R}^\sigma$ computed using eq.~\eqref{eq:Kelvin_cell_compressive_anisotropy} are reported in Figure~\ref{fig:eff_prop_ratios_cell}(b), together with those obtained from numerical simulations. Once more, the analytical model predictions agree well with the numerical data. $\mathcal{R}^E$ and $\mathcal{R}^\sigma$ both increase along with $\mathcal{R}$, and $\mathcal{R}^E$ develops much faster than $\mathcal{R}^\sigma$  (see also Figure~\ref{fig:eff_props_Kelvin_cell}(a)), which is opposite to the trend observed for the rectangular parallelepiped cell (see Figure~\ref{fig:eff_prop_ratios_cell}(a)). Moreover, for a given $\mathcal{R}$, the Kelvin cell has a much higher $\mathcal{R}^E$ while slightly lower $\mathcal{R}^\sigma$, than the rectangular parallelepiped cell.

The intermediate ratios $\mathcal{R}_\stt{f}$, $\mathcal{R}_\stt{c}$ and $\mathcal{R}_\theta$ computed using eq.~\eqref{eq:Kelvin_cell_ratios_2} are reported in Figure~\ref{fig:eff_prop_ratios_cell}(d). $\mathcal{R}_\stt{f}$ and $\mathcal{R}_\theta$ increases along with $\mathcal{R}$, while $\mathcal{R}_\stt{c}$ slightly decreases. According to eq.~\eqref{eq:Kelvin_cell_compressive_anisotropy}, $\mathcal{R}^E$ has a cubic dependency on $\mathcal{R}_\theta$ while $\mathcal{R}^\sigma$ a linear dependency, leading to much faster development of $\mathcal{R}^E$ than $\mathcal{R}^\sigma$, as shown in Figure~\ref{fig:eff_prop_ratios_cell}(b). Interestingly, comparing Figures~\ref{fig:eff_prop_ratios_cell}(c) and (d) indicates that the curve $\mathcal{R}_\stt{f}(\mathcal{R})$ for the Kelvin cell is similar to that for the rectangular parallelepiped cell. The curves $\mathcal{R}_\stt{c}(\mathcal{R})$ for the two idealized cell structures are however completely different, which may be related to different shapes of their load-bearing cell walls. Due to the additional strong dependency on $\mathcal{R}_\theta$, the curve $\mathcal{R}^E(\mathcal{R})$ for the Kelvin cell (see Figure~\ref{fig:eff_prop_ratios_cell}(b)) is well above that for the rectangular parallelepiped cell (see Figure~\ref{fig:eff_prop_ratios_cell}(a)).

For the sake of assessment, the predictions by another widely used analytical model, Sullivan model \cite{Sullivan2008}, are reported in Figure~\ref{fig:eff_prop_ratios_cell}(b). Detailed expressions of this reference model are given in~\ref{sec:ref_anal_models}. Despite being applied for closed-cell foams in many studies, the Sullivan model is derived by assuming an open Kelvin cell structure, which is in principle inappropriate for closed-cell foams with a high cell face fraction. As expected, large deviations can be observed on the Sullivan model predictions, especially for the strength anisotropy, which is overestimated by $>200\%$.

\subsection{Discussion}
The present analytical models have shown capabilities to accurately reproduce the mechanical anisotropy obtained from the idealized cell-based numerical models. Detailed analysis on the impacts of cell shape anisotropy indicates that:
\begin{enumerate}
\item Cell shape anisotropy translates into mechanical anisotropy through three pathways, cell load-bearing area fraction, cell wall buckling strength and cell wall inclination angle.
\item The inclination angle plays an critical role in determining mechanical anisotropy, in particular modulus anisotropy.
\end{enumerate}
The specific relationships between mechanical anisotropy and cell shape anisotropy would vary from one case to another, depending on the competition among the three pathways above. \textcolor{red}{The base material yield stress is relevant for the compressive strengths in different global directions but does not contribute to the strength anisotropy}.

In addition, two widely used analytical models \cite{Gibson1997_ch6, Sullivan2008} exhibit large predictive deviations already for the idealized cell structures, where consistent geometrical assumptions are adopted. These deviations are believed to originate from the introduced mechanistic assumptions, which appear to be inappropriate for closed-cell foams with a high cell face fraction and low relative density. This explains why the predictive capabilities of these analytical models for realistic foams can sometimes be quite low (see e.g.\ \cite{EspadasEscalante2015, Doyle2019, Liu2020, Zhou2023b}).

\section{Analyses of the tessellation-based models}\label{sec:analy_tess_meso_model}
Numerical results of the tessellation-based models incorporating different mesostructural stochastics, will be analysed in this section.

\subsection{Deformation mechanisms}
\begin{figure}[ht]
\centering
\begin{overpic}[draft=false,width=0.9\textwidth]{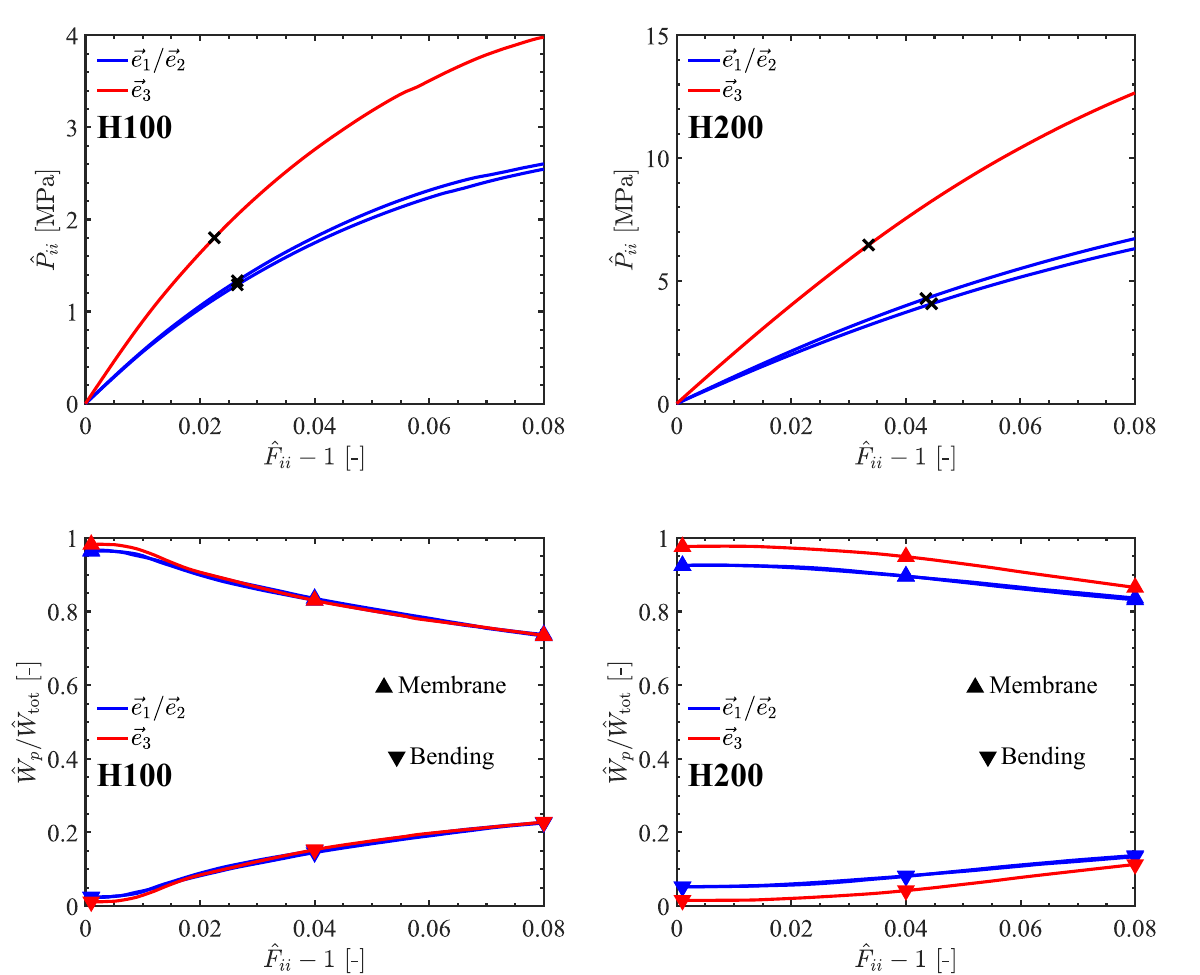}
\put(24.5,79.5){\color{black}\footnotesize \textbf{(a)}}
\put(73.625,79.5){\color{black}\footnotesize \textbf{(b)}}
\put(24.5,37.625){\color{black}\footnotesize \textbf{(c)}}
\put(73.625,37.625){\color{black}\footnotesize \textbf{(d)}}
\end{overpic}
\caption{(a-b) Effective stresses and (c-d) strain energy fractions versus applied strain of the tessellation-based model set ``StSt" for two Divinycell foam grades, under uniaxial compression in the transverse ($\vec{e}_1/\vec{e}_2$) and foam rise ($\vec{e}_3$) directions. The yield points are indicated in (a-b) by the black crosses.}\label{fig:eff_response_foam}
\end{figure}
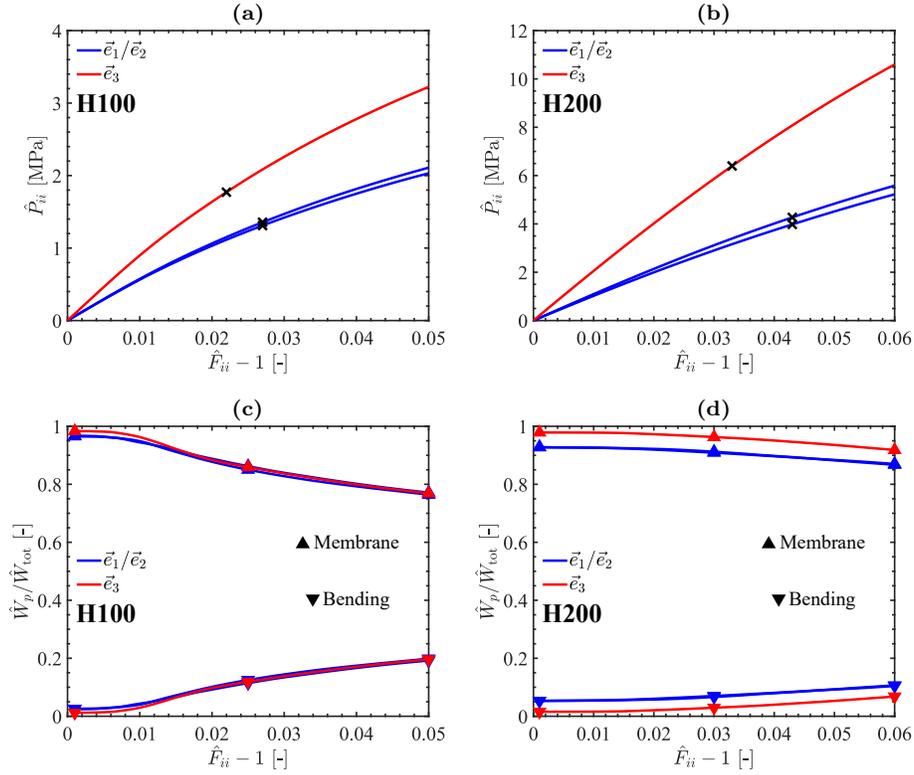

The macroscale effective responses of the model set ``StSt" for Divinycell foam H100 and H200, which account for the stochastic variations of cell size, cell wall thickness and cell shape anisotropy, are reported in Figure~\ref{fig:eff_response_foam}. As expected, the effective responses in two transverse directions ($\vec{e}_1/\vec{e}_2$) are quite close. Figure~\ref{fig:eff_response_foam}(a) shows that for each loading case of H100, the stress first increases linearly, followed by a continuous stiffness reduction. The compressive stress in the foam rise direction ($\vec{e}_3$) is higher than the transverse direction ($\vec{e}_1/\vec{e}_2$), indicating an anisotropic compressive behaviour. Figure~\ref{fig:eff_response_foam}(c) shows that for each loading direction, the initial elastic region is dominated by the membrane deformation mode ($\hat{W}_\stt{m}/\hat{W}_\stt{tot}> 0.95$), followed by a continuous increase of the bending contribution. Similar trends can be observed for H200 from Figures~\ref{fig:eff_response_foam}(b) and (d). Yet, the stiffness reduction is less pronounced, and the bending contribution increase rate is lower than H100.

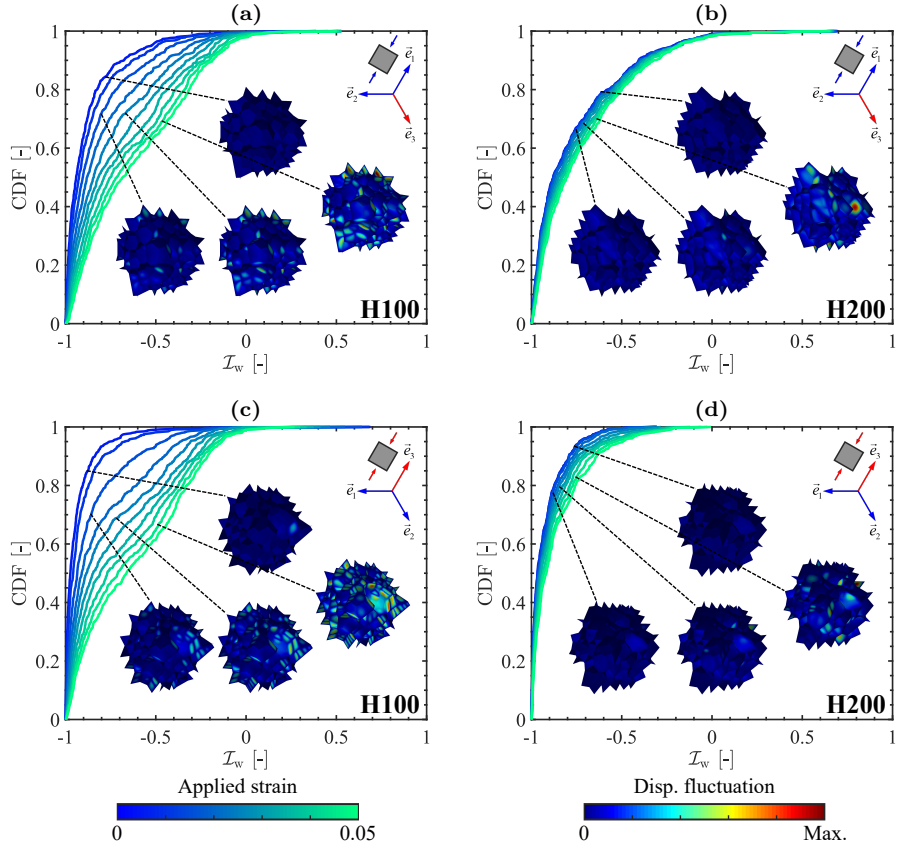
\begin{figure}[ht]
\centering
\begin{overpic}[draft=false,width=0.9\textwidth]{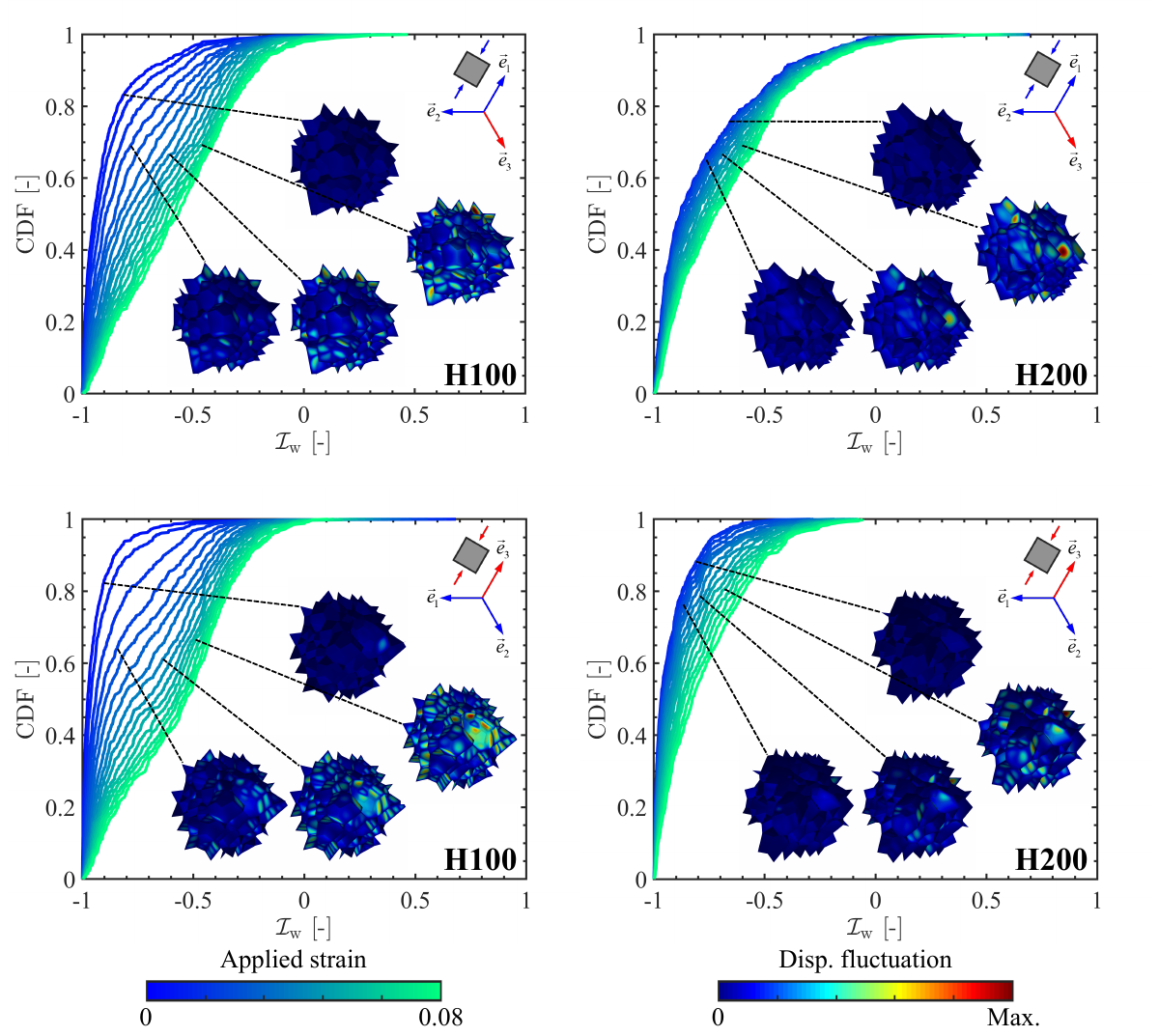}
\put(24.5,87.625){\color{black}\footnotesize \textbf{(a)}}
\put(73.625,87.625){\color{black}\footnotesize \textbf{(b)}}
\put(24.5,45.75){\color{black}\footnotesize \textbf{(c)}}
\put(73.625,45.75){\color{black}\footnotesize \textbf{(d)}}
\end{overpic}
\caption{Cumulative density functions (CDF) of the cell wall strain energy partitioning indicators and deformed configurations at different stages of the tessellation-based model set ``StSt" for two Divinycell foam grades, under uniaxial compression in the (a-b) transverse and (c-d) foam rise directions. The loading direction is represented by a pair of opposite arrows.}\label{fig:cell_wall_foam}
\end{figure}

To interpret the observations in Figure~\ref{fig:eff_response_foam}, cumulative density functions (CDF) of the cell wall strain energy partitioning indicators at different applied strains, are reported in Figure~\ref{fig:cell_wall_foam}.  Figures~\ref{fig:cell_wall_foam}(a) and (c) demonstrate that $>80\%$ and $>90\%$ of cell walls of H100 deform by a nearly pure membrane mode ($\mathcal{I}_\stt{w}<-0.8$) at the initial stage, under compression in the transverse and foam rise directions, respectively. As the loading proceeds, more and more cell walls buckle and switch to a mixed membrane-bending mode, \textcolor{red}{accompanied by a stress redistribution and load-carrying efficiency reduction (not shown)}. Attributed to the large variations of cell size and cell wall thickness (see Table~\ref{tab:tess_geo_parameters}), the buckling resistance greatly varies between individual cell walls, leading to sequential occurrence of buckling and thus gradual energy redistribution in Figure~\ref{fig:eff_response_foam}(c). This holds for both the transverse and foam rise directions. 

Figures~\ref{fig:cell_wall_foam}(b) and (d) for H200 demonstrate similar trends as those for H100. However, fewer cell walls buckle, resulting in a slower deformation mode transition compared with H100. This can be understood by the larger cell wall thickness of H200 (see Table~\ref{tab:tess_geo_parameters}), which gives rise to a higher buckling resistance. Numerical simulations on ``StCt" and ``CtCt" have delivered qualitatively similar results as ``StSt" and are thus omitted here.

The above analyses confirm that the deformation mechanisms identified using the idealized cell structures (see \textcolor{blue}{Section}~\ref{sec:analy_ideal_meso_model}) remain valid in the presence of mesostructural stochastics.

\begin{figure}[ht]
\centering
\begin{overpic}[draft=false,width=0.9\textwidth]{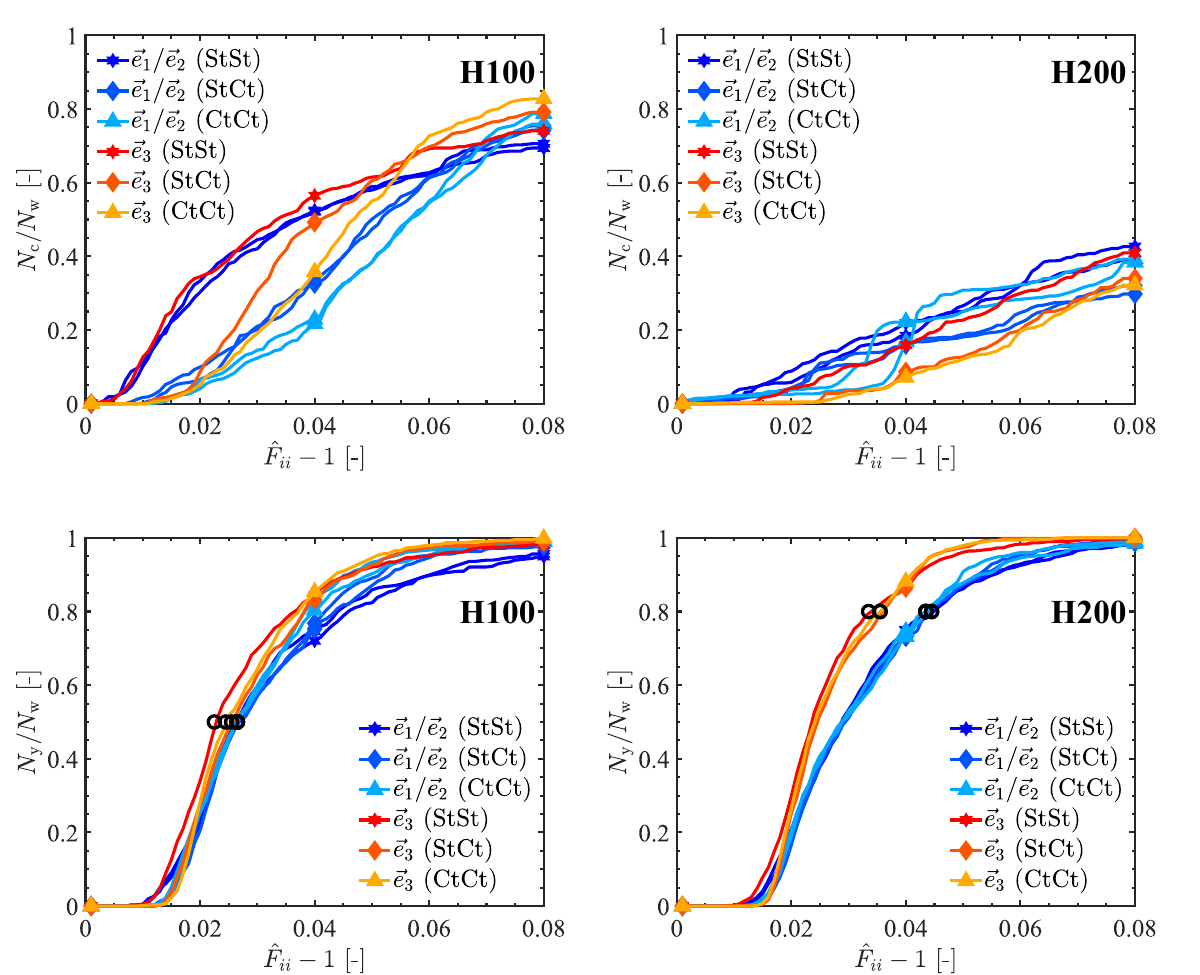}
\put(24.5,79.5){\color{black}\footnotesize \textbf{(a)}}
\put(73.625,79.5){\color{black}\footnotesize \textbf{(b)}}
\put(24.5,37.625){\color{black}\footnotesize \textbf{(c)}}
\put(73.625,37.625){\color{black}\footnotesize \textbf{(d)}}
\end{overpic}
\caption{\textcolor{red}{Fractions of (a-b) buckled and (c-d) yield cell walls versus applied strain of different tessellation-based model sets for two Divinycell foam grades, under uniaxial compression in the transverse ($\vec{e}_1/\vec{e}_2$) and foam rise ($\vec{e}_3$) directions. The yield points are indicated in (c-d) by the black circles}.}\label{fig:cell_wall_BK_Y_foam}
\end{figure}

\textcolor{red}{To reveal the impacts of mesostructural stochastics on the cell wall deformation behaviour, the fractions of buckled and yield cell walls of different model sets are provided in Figure~\ref{fig:cell_wall_BK_Y_foam}. Here, ${N}_\stt{c}$, ${N}_\stt{y}$ and ${N}_\stt{w}$ denote the numbers of buckled, yield and all cell walls, respectively. Figures~\ref{fig:cell_wall_BK_Y_foam}(a-b) show that the cell wall buckling events get promoted as more mesostructural stochastics are included (from ``CtCt" to ``StSt").  This can be explained through the weakest link principle (see also \cite{Shi2018, Vengatachalam2019}). Introducing the stochastic variations of more mesostructural features gives rise to the emergence of more weak regions. As expected, the cell wall buckling events for H200 are less active (see also Figure~\ref{fig:cell_wall_foam}) and accompanied by a lower growth rate of ${N}_\stt{c}/{N}_\stt{w}$, compared with H100}.

\textcolor{red}{Interestingly, Figures~\ref{fig:cell_wall_BK_Y_foam}(c-d) show that the cell wall yield events remain nearly unaffected by including more mesostructural stochastics (from ``CtCt" to ``StSt"). This is likely because the membrane deformation mode is dominating under compression. The cell wall yield events of H100 only start slightly earlier than H200, despite their highly different relative densities. Focusing on the model set ``StSt", attributed to the presence of many weak cell walls, the buckling events tend to occur before the yielding events, even for H200 which has a relative density (see Table~\ref{tab:tess_geo_parameters}) higher than the critical transition relative density \cite{Kidd2012}}.

\subsection{Effective properties}
\begin{figure}[ht]
\centering
\begin{overpic}[draft=false,width=0.909\textwidth]{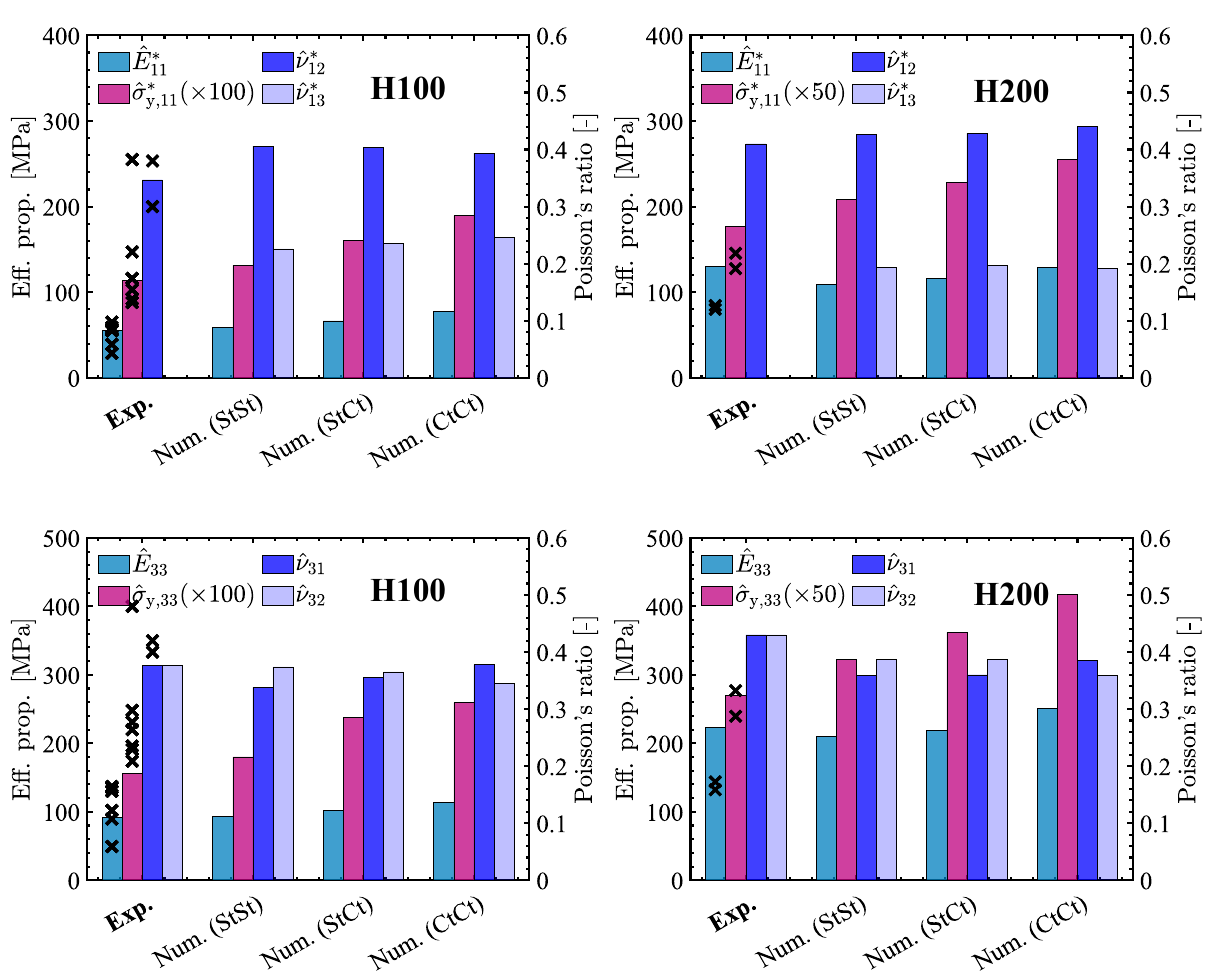}
\put(23.625,78.625){\color{black}\footnotesize \textbf{(a)}}
\put(73.375,78.625){\color{black}\footnotesize \textbf{(b)}}
\put(23.625,37.25){\color{black}\footnotesize \textbf{(c)}}
\put(73.375,37.25){\color{black}\footnotesize \textbf{(d)}}
\end{overpic}
\caption{Effective compressive properties of different tessellation-based model sets for two Divinycell foam grades. $(\bullet)^*$ in (a-b) indicates the average of quantities in two transverse directions ($\vec{e}_1/\vec{e}_2$). Experimental data collected from extensive literature are provided for reference. H100 data from \cite{Shafiq2015} and H200 data from \cite{Funari2021} are indicated by the leftmost bars. H100 data from \cite{Zhang2012, Taher2012, Chen2013, Liu2020, Funari2021, HooFatt2022, Tang2022, Tong2022} and H200 data from \cite{Liu2020, Magliaro2023} are indicated by the black crosses.}\label{fig:eff_props_foam}
\end{figure}

The effective compressive properties of different model sets are reported in Figure~\ref{fig:eff_props_foam}. \textcolor{red}{Here, the yield strength is determined using a yield criterion, i.e.\ when a sufficient number of cell walls yield. The critical ${N}_\stt{y}/{N}_\stt{w}$ of H100 and H200 are taken as 0.5 and 0.8 (see Figures~\ref{fig:cell_wall_BK_Y_foam}(c-d)), respectively, which are calibrated according to the strains at the the peak stress points of the experimental stress-strain curves in \cite{Shafiq2015} and \cite{Funari2021}. The resulting yield points for ``StSt" have been supplemented to Figures~\ref{fig:eff_response_foam}(a-b) as an example}.

Comparing the results of different model sets for H100 (see Figures~\ref{fig:eff_props_foam}(a) and (c)) shows that for each loading direction, as the stochastic variations of cell size and cell wall thickness are sequentially incorporated (from ``CtCt" to ``StSt"), the compressive modulus and strength both decrease, while the Poisson's ratios remain almost unchanged. Compared with the compressive modulus, the strength is more sensitive to these mesostructural stochastics. Similar trends hold for H200 (see Figures~\ref{fig:eff_props_foam}(b) and (d)). Given the higher relative density of H200 (see Table~\ref{tab:tess_geo_parameters}), the resulting compressive moduli and strengths are apparently higher than H100. The Poisson's ratios of H100 and H200 are quite similar despite their highly different mesostructures.

Nevertheless, one should be careful with interpreting the impacts of cell size stochastics. Based on empirical relationships between the compressive properties and relative density for closed-cell foams \cite{Gibson1997_ch5}, the net impacts of cell size stochastics may be secondary in practice. When analysing the model sets for each Divinycell foam grade, the resulting overall relative density $\rho/\rho_\stt{r}$ of ``CtCt" is $\sim 10\%$ higher than those of ``StCt" and ``StSt" (see \textcolor{blue}{Section}~\ref{sec:num_model_method}). At the same time, the compressive properties of ``CtCt" are $\sim 15\%$ higher than ``StCt" (see Figure~\ref{fig:eff_props_foam}). It is likely that this difference in compressive properties is mainly associated with the change in $\rho/\rho_\stt{r}$  rather than to the change in the cell size stochastics. To examine this inference, numerical simulations of ``CtCt" H100 and H200, with the cell wall thickness scaled such that the the resulting $\rho/\rho_\stt{r}$ are equal to those of ``StCt", are performed. It has been found that the resulting compressive moduli and strengths of the scaled ``CtCt" are quite close to those of ``StCt", with a relative difference $< 5\%$. This implies that the compressive properties receive secondary impacts from the cell size stochastics in practice, despite being still noticeable\footnote{More pronounced impacts by the cell size stochastics are claimed in other numerical studies \cite{Chen2015, Chen2017a}, which focus on Gurit Corecell foam M130, nearly isotropic. Mesostructural models with different cell size distributions while the overall relative density preserved, are considered. For a cell size distribution comparable to the present study, the resulting compressive modulus and strength are found to decrease by $\sim 5\%$ and $\sim 10\%$, respectively, compared with the case with a constant cell size. However, the overall cell equivalent diameter is not preserved and increases along with increasing cell size stochastic variations. This would already weaken especially the compressive strength. Therefore, we believe that the net impacts of cell size stochastics is less pronounced than what have been claimed in \cite{Chen2015, Chen2017a}.}.

The above observed decreasing trends of compressive properties with increasing mesostructural stochastics can be well linked to the observations in Figures~\ref{fig:cell_wall_BK_Y_foam}(a-b).

For the sake of reference, the experimental compressive properties of H100 \cite{Shafiq2015} and H200 \cite{Funari2021}, are provided in Figure~\ref{fig:eff_props_foam} (leftmost bars). $\hat{\nu}_{13}^*$ are not measured and $\hat{\nu}_{31}=\hat{\nu}_{32}$ is assumed in \cite{Shafiq2015, Funari2021}. As an indication, the experimental data from other literature are also provided in Figure~\ref{fig:eff_props_foam} (black crosses), although these studies are lacking either well-defined strain measurements or complete stress-strain curves under uniaxial compression. A remarkable inconsistency between the experimental data reported in different literature can be observed. This places a clear need of more attention to the experimental aspects, e.g.\ test method, sample shape, sample size and determination of compressive properties.

In the following, the numerical model predictions are compared with the experimental data from \cite{Shafiq2015, Funari2021} only, given their reliability and relevance. The model set ``StSt", with all the cell size, cell wall thickness and cell shape anisotropy stochastics incorporated, seems to deliver the closest predictions with respect to the experimental data. In particular for the compressive moduli and Poisson's ratios, an excellent agreement between the experimental data and ``StSt" predictions can be observed. Relatively large deviations appear on the strengths, which are overestimated by $\sim 15\%$, which shall be attributed to the disregarded plasticity in the present material modelling.

\subsection{Mechanical anisotropy}
\begin{figure}[ht]
\centering
\begin{overpic}[draft=false,width=0.9\textwidth]{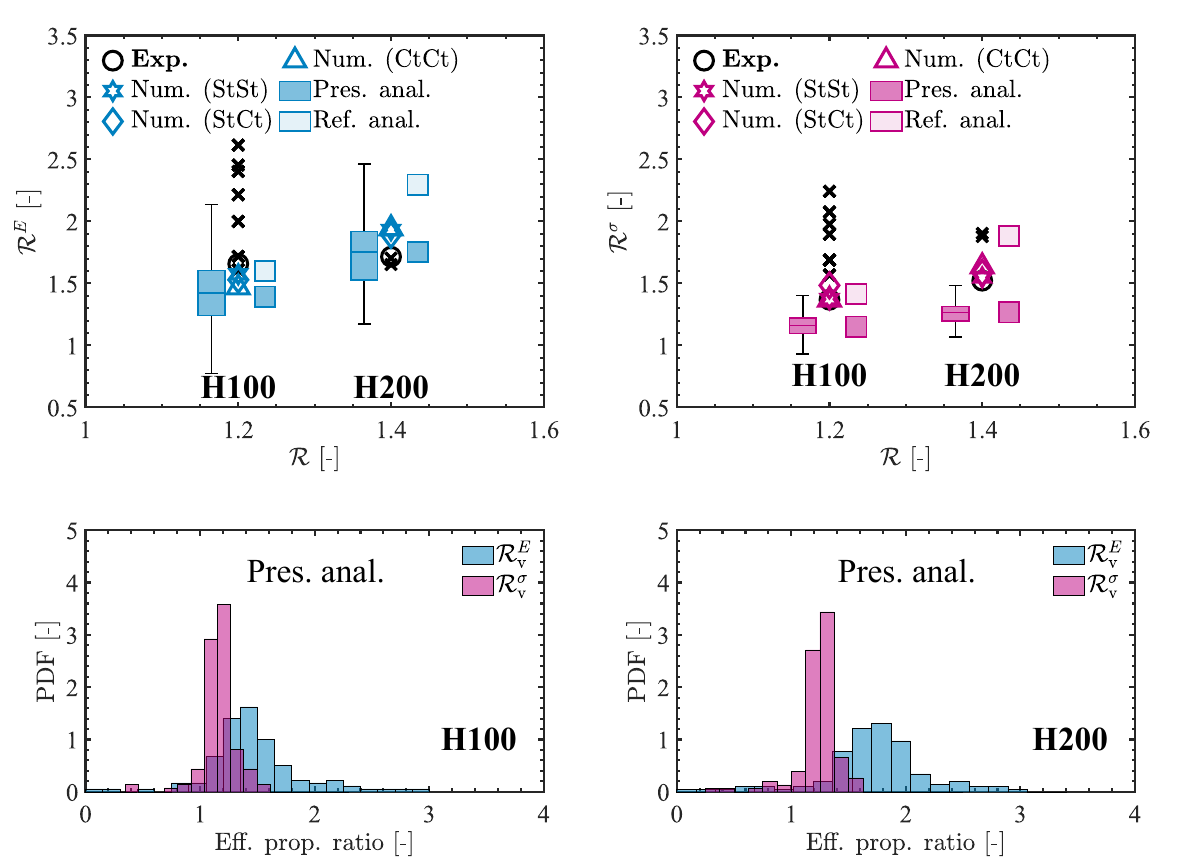}
\put(24.5,69.75){\color{black}\footnotesize \textbf{(a)}}
\put(73.625,69.75){\color{black}\footnotesize \textbf{(b)}}
\put(24.5,28.5){\color{black}\footnotesize \textbf{(c)}}
\put(73.625,28.5){\color{black}\footnotesize \textbf{(d)}}
\end{overpic}
\caption{Mechanical anisotropy predicted by different tessellation-based numerical models, and the idealized cell-based present~\eqref{eq:Kelvin_cell_compressive_anisotropy} and reference \cite{Sullivan2008} analytical models for two Divinycell foam grades: (a) modulus and (b) yield strength; (c-d) probability density functions (PDF) of the cell mechanical anisotropy computed using the present analytical model~\eqref{eq:Kelvin_cell_compressive_anisotropy}. The analytical model predictions and the five-number summary statistics in (a-b) are shifted horizontally for better visibility. Experimental data collected from extensive literature are provided in (a-b) for reference. H100 data from \cite{Shafiq2015} and H200 data from \cite{Funari2021} are indicated by the black circles. H100 data from \cite{Zhang2012, Taher2012, Chen2013, Liu2020, Funari2021, HooFatt2022, Tang2022, Tong2022} and H200 data from \cite{Liu2020, Magliaro2023} are indicated by the black crosses.}\label{fig:eff_prop_ratios_foam_2}
\end{figure}
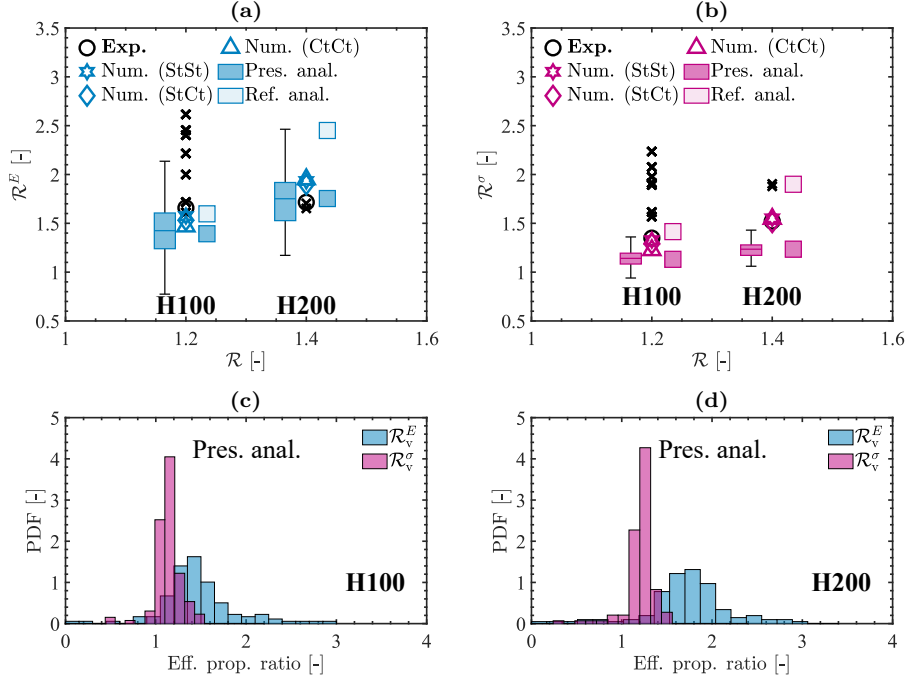

With the effective compressive properties in Figure~\ref{fig:eff_props_foam}, the mechanical anisotropy $\mathcal{R}^E$ and $\mathcal{R}^\sigma$ of different model sets are computed, and reported in Figures~\ref{fig:eff_prop_ratios_foam_2}(a) and (b), respectively. Three model sets (``StSt", ``StCt" and ``CtCt"), are found to deliver comparable predictions of both $\mathcal{R}^E$ and $\mathcal{R}^\sigma$, with the relative difference in between $<10\%$. This implies that the cell wall thickness and cell size stochastics only weakly affect the resulting mechanical anisotropy.

Again, $\mathcal{R}^E$ and $\mathcal{R}^\sigma$ computed using the experimental data in Figure~\ref{fig:eff_props_foam}, are supplemented to Figures~\ref{fig:eff_prop_ratios_foam_2}(a) and (b), respectively (black markers). It seems that the experimental mechanical anisotropy (black circles) from \cite{Shafiq2015, Funari2021}, especially strength anisotropy, can be well reproduced using any of the three model sets. Nevertheless, given the large inconsistency among the experimental data from different literature (see also Figure~\ref{fig:eff_props_foam}), it is hardly feasible to conduct any in-depth analysis regarding the accuracy of model predictions.

In order to identify the influence of cell shape anisotropy stochastics, an idealized foam mesostructural model is introduced as an array of periodically repeated Kelvin cells, which has been frequently used in the literature (see e.g.\ \cite{Song2010, Shi2018, Vengatachalam2019, Gahlen2022a}). $\mathcal{R}^E$ and $\mathcal{R}^\sigma$ of the idealized model for each Divinycell foam grade can be determined by substituting the overall cell shape anisotropy $\mathcal{R}$ (see Table~\ref{tab:tess_geo_parameters}) into the present analytical model~\eqref{eq:Kelvin_cell_compressive_anisotropy}, which has been validated against numerical simulations (see Figure~\ref{fig:eff_prop_ratios_cell}(b)) \textcolor{red}{and holds for a broad range of relative densities (see~\ref{sec:influence_plas})}. The results are reported in Figures~\ref{fig:eff_prop_ratios_foam_2}(a) and (b), respectively (dark cyan and pink squares). It can be seen that the idealized model underestimates $\mathcal{R}^E$ and $\mathcal{R}^\sigma$ with respect to the tessellation-based models. In particular, the predictive deviation on $\mathcal{R}^\sigma$ is $>30\%$. This may be attributed to the high sensitivity of compressive strength to mesostructural stochastics.

To unravel the idealized model predictive deviations in more detail, the individual cell mechanical anisotropy $\mathcal{R}_\stt{v}^E$ and $\mathcal{R}_\stt{v}^\sigma$ are computed by substituting the shape anisotropy $\mathcal{R}_\stt{v}$ (see Figures~\ref{fig:PDF_meso_foam}(a-b)) into eq.~\eqref{eq:Kelvin_cell_compressive_anisotropy}. The corresponding probability density functions (PDF) for H100 and H200 are reported in Figures~\ref{fig:eff_prop_ratios_foam_2}(c) and (d), respectively. Because of the stronger dependency on $\mathcal{R}_\stt{v}$ (see Figure~\ref{fig:eff_prop_ratios_cell}(b)), $\mathcal{R}_\stt{v}^E$ exhibits a larger spread compared with $\mathcal{R}_\stt{v}^\sigma$. Furthermore, the five-number summary statistics of $\mathcal{R}_\stt{v}^E$ and $\mathcal{R}_\stt{v}^\sigma$ are indicated in Figures~\ref{fig:eff_prop_ratios_foam_2}(a) and (b), respectively (dark cyan and pink windows). It can be seen that the upper quartile (top edge of the window) of $\mathcal{R}_\stt{v}^E$ closely represents $\mathcal{R}^E$ of the tessellation-based models, while the upper bound (top black edge over the window) of $\mathcal{R}_\stt{v}^\sigma$ corresponds well with $\mathcal{R}^\sigma$ of the tessellation-based models. These observations suggest that the mechanical anisotropy, especially strength anisotropy, of a foam mesostructure with random cell shape anisotropy, is dominated by the cells with a relatively large shape anisotropy and cannot be simply correlated to the overall cell shape anisotropy.

Based on the comparative study above, it is concluded that:
\begin{enumerate}
\item The cell shape anisotropy stochastics have strong impacts on the resulting  mechanical anisotropy, in particular strength anisotropy.
\item The cell size and cell wall thickness stochastics play a rather secondary role.
\end{enumerate}
A model without taking into account the cell shape anisotropy stochastics would apparently underestimate the mechanical anisotropy of realistic foams. \textcolor{red}{It has been confirmed that the present findings remain valid, even when larger stochastic variations of cell size and cell wall thickness than those in Table~\ref{tab:tess_geo_parameters} are considered}.

Besides, $\mathcal{R}^E$ and $\mathcal{R}^\sigma$ predicted by the reference analytical model, Sullivan model \cite{Sullivan2008}, are reported in Figures~\ref{fig:eff_prop_ratios_foam_2}(a) and (b), respectively (light cyan and pink squares). It can be seen that the Sullivan model overestimates both $\mathcal{R}^E$ and $\mathcal{R}^\sigma$ with respect to the tessellation-based models. Especially for H200, the predictive deviations are $>30\%$. Interestingly for H100, the Sullivan model demonstrates even better predictive capabilities than the present analytical model.

The Sullivan model seems able to deliver reasonably good predictions. This is, however, a consequence of two sources of deviations compensating for each other. On one hand, introducing inappropriate mechanistic assumptions leads to the mechanical anisotropy significantly overestimated (see Figure~\ref{fig:eff_prop_ratios_cell}(b)). On the other hand, disregarding the cell shape anisotropy stochastics results in the mechanical anisotropy apparently underestimated (see above discussions for the present analytical model). Depending on the competition between the two sources of deviations, the predictive capabilities of the Sullivan model would largely vary from one case to another. This has been recognized in many studies for realistic foams (see e.g.\ \cite{Sullivan2008, EspadasEscalante2015, Andersons2016, Doyle2019}), where the experimental uncertainties also play an important role.

\section{Conclusions}\label{sec:conclusions}
Focusing on closed-cell foams with a high cell face fraction and \textcolor{red}{low relative density ($<0.15$)}, a systematic numerical study has been performed to investigate the anisotropic compressive behaviour, which takes into account cell shape anisotropy, cell structure and different mesostructural stochastics. The main findings are summarized:
\begin{enumerate}
\item The anisotropic compressive properties of Divinycell foam H100 and H200 predicted by the tessellation-based models that include all the cell size, cell wall thickness and cell shape anisotropy stochastics, can closely represent the experimental results in \cite{Shafiq2015, Funari2021}. Yet, as a remark, there is a large inconsistency among the experimental data from different literature, calling for more attention to the experimental aspects.
\item The cell wall membrane deformation dominates the initial elastic region, irrespective of the loading direction. Compared with this predominant deformation mechanism, the cell wall bending contribution is small at the initial stage and becomes important only \textcolor{red}{after buckling, followed by membrane yielding}.
\item The anisotropy of compressive properties is related to cell shape anisotropy through three pathways, cell load-bearing area fraction, cell wall buckling strength and cell wall inclination angle. The inclination angle has crucial impacts on the resulting mechanical anisotropy, in particular modulus anisotropy. \textcolor{red}{The base material yield stress does not contribute to the strength anisotropy, despite being relevant for the compressive strengths in different global directions}.
\item The cell shape anisotropy stochastics strongly affect the anisotropy of compressive properties, in particular strength anisotropy. In contrast, the impacts of the cell size and cell wall thickness stochastics are much less important.
\item The mechanistic assumptions introduced in the two widely used analytical models \cite{Gibson1997_ch6, Sullivan2008} appear to be inappropriate at a high cell face fraction and low relative density. This becomes another key source of deviations besides different uncertainties in the real foam mesostructures and experiments, and explains why the predictive capabilities of these analytical models can sometimes be quite low (see e.g.\ \cite{EspadasEscalante2015, Doyle2019, Liu2020, Zhou2023b}).
\end{enumerate}

Through quantitative analysis of the cell wall deformation behaviour, this contribution confirms the dominant role of membrane deformation in the initial elastic region, as already suggested by other studies (see e.g.\ \cite{Simone1998a, Grenestedt2000, Shi2018}). The present findings on the impacts of cell shape anisotropy, cell structure and different mesostructural stochastics, provide deeper insights into how the anisotropic compressive properties are related to mesostructural features. The developed analytical models that describe the relationships between mechanical anisotropy and cell shape anisotropy, may provide new design guidelines for not only traditional foams, but also lattice structures consisting of regular cells (see e.g.\ \cite{Berger2017, TancogneDejean2018, Guo2024}).

\textcolor{red}{It has been assumed that the cell wall elastic buckling dominates the compressive failure, which is qualitatively supported by detailed experimental observations (see e.g.\ \cite{Kidd2012, Duan2019, Chai2020}). To mitigate uncertainty and achieve quantitative validation, an integrated experimental-numerical study remains necessary. This would require careful experimental characterization (including scale-consistent characterization of the base material properties) and development of high-fidelity numerical models. Many other mesostructural feastures in realistic foams remain to be accounted for, which may influence the cell wall deformation behaviour}. For instance, cell walls tend to be thicker around the edges and thinner towards the face centers (see e.g.\ \cite{Jang2015, Tang2022}), and may undergo apparent distortion/damage during the manufacturing process, \textcolor{red}{which potentially causes apparent initial curvature with corrugations and wriggles, and even the absence of several cell walls (see e.g.\ \cite{Andrews1999, Jeon2005, PerezTamarit2019})}. \textcolor{red}{Besides the stochastic variations of cell size, cell wall thickness and cell shape, there is likely a spread on individual cell elongation directions. Also, the spatial variations of different mesostructural features may depend on each other. Incorporating all the above would require the use of more general tessellation techniques (see e.g.\ \cite{Sonon2015, Ghazi2019, Ghazi2020b}), in combination with quantitative experimental characterization}. Moreover, the cell wall plasticity becomes important upon a large applied strain, in particular to capture the plateau region observed on the compressive response. In addition, the strain rate effects and orientation effects of base materials may influence the anisotropic compressive behaviour. These aspects have not been systematically investigated here and will be addressed in the next steps.

\section*{Acknowledgements}
This research was carried out under project numbers 2020-04526 and 2023-01937 in the framework of the Strategic Innovation Programme LIGHTer, funded by the Swedish Agency for Innovation Systems (Vinnova) and supported by the Swedish government.

\appendix

\color{red}
\section{Influence of the cell wall curvature}\label{sec:influence_cell_wall_curv}
Cell walls in realistic foams have more or less initial curvature, which is known to influence the foam compressive properties. This effect has been disregarded in the present study and is preliminarily investigated in the following.

In order to gain general insights, how the initial curvature influences the cell wall compressive behaviour in the long direction is focused on. Rectangular cell walls are modelled, which are parametrized by length $L_\stt{w}$, width $B_\stt{w}$ and thickness $t$ (see also Figure~\ref{fig:geo_model_cell_wall}). Assuming that the cell wall curvature is induced by the gas pressure difference between cells, it is fair to parametrize the initial curvature pattern by a bubble function \cite{Nakshatrala2007}:
\begin{equation}\label{eq:bubble_fun}
h=h_0\left(1-\left(\frac{2l_\stt{w}}{L_\stt{w}}\right)^2\right)\left(1-\left(\frac{2b_\stt{w}}{B_\stt{w}}\right)^2\right),
\end{equation}
where $l_\stt{w}$ and $b_\stt{w}$ denote the in-plane coordinates in the length and width directions, with the cell wall center as the origin point; $h$ characterizes the out-of-plane geometrical deviation of a curved cell wall with respect to a flat one, and reaches its maximum $h_0$ at the cell wall center. The initial curvature level can be indicated using $h_0/\sqrt{A_\stt{w}}$, with $A_\stt{w}$ being the in-plane area.

Cell wall models for the two aspect ratios $\mathcal{R}_{\rm w}=1.0$ and $\mathcal{R}_{\rm w}=2.0$ are taken as examples, due to their different critical buckling modes \cite{Gerard1957}. $L_\stt{w}$ and $B_\stt{w}$ are scaled according to the prescribed $\mathcal{R}_{\rm w}$, with $A_\stt{w}$ and $t$ preserved, i.e.\ $L_{\rm w}=\sqrt{A_\stt{w}}\mathcal{R}_{\rm w}^{\frac{1}{2}}$ and $B_{\rm w}=\sqrt{A_\stt{w}}\mathcal{R}_{\rm w}^{\stt{-}\frac{1}{2}}$. The reference geometrical parameters are taken as $L_\stt{w}=0.4$ [mm], $B_\stt{w}=0.4$ [mm] and $t=0.01$ [mm], at $\mathcal{R}_{\rm w}=1.0$. For each $\mathcal{R}_{\rm w}$, model configurations with different $h_0/\sqrt{A_\stt{w}}$ are considered. Notice that the resulting cell wall volume would slightly increase as $h_0/\sqrt{A_\stt{w}}$ increases. 

FE discretization strategy and material model along with parameters follow those in \textcolor{blue}{Section}~\ref{sec:num_model_method}. Motivated by the observations in \textcolor{blue}{Section}~\ref{sec:relation_anisotropy}, the cell wall edges are fully clamped and uniaxial compressive loading in the length direction is applied. A small perturbation force is imposed at the cell wall center to trigger buckling.

\begin{figure}[ht]
\centering
\begin{overpic}[draft=false,width=0.9\textwidth]{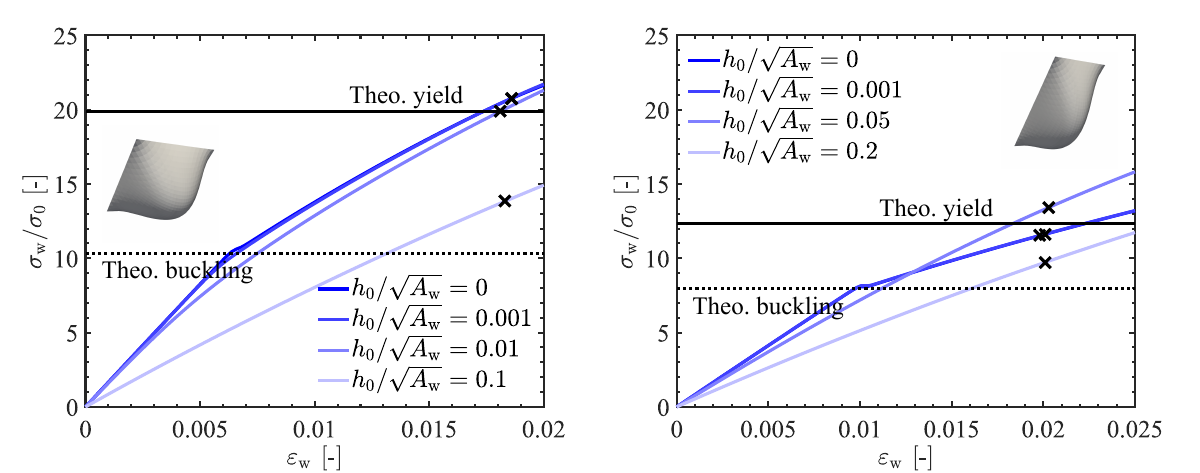}
\put(24.5,37.5){\color{black}\footnotesize \textbf{(a)}}
\put(73.625,37.5){\color{black}\footnotesize \textbf{(b)}}
\end{overpic}
\caption{\textcolor{red}{Compressive responses of the rectangular cell wall models with different initial curvature levels $h_0/\sqrt{A_{\rm w}}$, for the two aspect ratios: (a) $\mathcal{R}_{\rm w}=1.0$ and (b) $\mathcal{R}_{\rm w}=2.0$. The yield points are indicated by the black crosses}.}\label{fig:eff_response_cell_wall}
\end{figure}

The compressive stresses $\sigma_\stt{w}$ of different cell wall models, for the two aspect ratios $\mathcal{R}_{\rm w}=1.0$ and $\mathcal{R}_{\rm w}=2.0$, are plotted as functions of the applied strain $\varepsilon_\stt{w}$ in Figures~\ref{fig:eff_response_cell_wall}(a) and (b), respectively. Here, $\sigma_0=\sigma_\stt{c,w}/k_\stt{c}$ has been introduced, with $\sigma_\stt{c,w}$ being the theoretical buckling strength and $k_{\rm c}$ the theoretical buckling coefficient \cite{Gerard1957}. As expected, for each flat cell wall ($h_0/\sqrt{A_\stt{w}}=0$), the stress first increases linearly, followed by a sudden stiffness reduction, indicating the buckling event. The cell wall continues to carry additional load and eventually reaches its yield strength (black crosses), which is determined using the membrane yielding detector (\ref{eq:Y_detector}). The stresses at the stiffness reduction points, i.e.\ buckling points, agree well with the theoretical buckling strengths $\sigma_\stt{c,w}$ (with $k_\stt{c}=10.35$ at $\mathcal{R}_{\rm w}=1.0$ and $k_\stt{c}=7.95$ at $\mathcal{R}_{\rm w}=2.0$ \cite{Gerard1957}), and those at the yield points correspond well with the theoretical yield strengths $\sigma_\stt{y,w}=\sqrt{\sigma_\stt{c,w}\sigma_\stt{y}}$  \cite{Timoshenko1961}, with $\sigma_\stt{y}$ being the material yield stress.

The impacts of $h_0/\sqrt{A_\stt{w}}$ is then discussed. For $\mathcal{R}_{\rm w}=1.0$, as $h_0/\sqrt{A_\stt{w}}$ increases, the compressive response gets weakened, resulting in the compressive modulus and strength both reduced (see Figure~\ref{fig:eff_response_cell_wall}(a)). For $\mathcal{R}_{\rm w}=2.0$, the initial compressive response and compressive modulus remain a decreasing trend with increasing $h_0/\sqrt{A_\stt{w}}$ (see Figure~\ref{fig:eff_response_cell_wall}(b)). However, the later compressive response and compressive strength sometimes get strengthened instead ($h_0/\sqrt{A_\stt{w}}=0.05$). The yield strains for each $\mathcal{R}_{\rm w}$ seems independent on $h_0/\sqrt{A_\stt{w}}$. Notice that the buckling point is no more visible at a large $h_0/\sqrt{A_\stt{w}}$. 

\begin{figure}[ht]
\centering
\begin{overpic}[draft=false,width=0.9\textwidth]{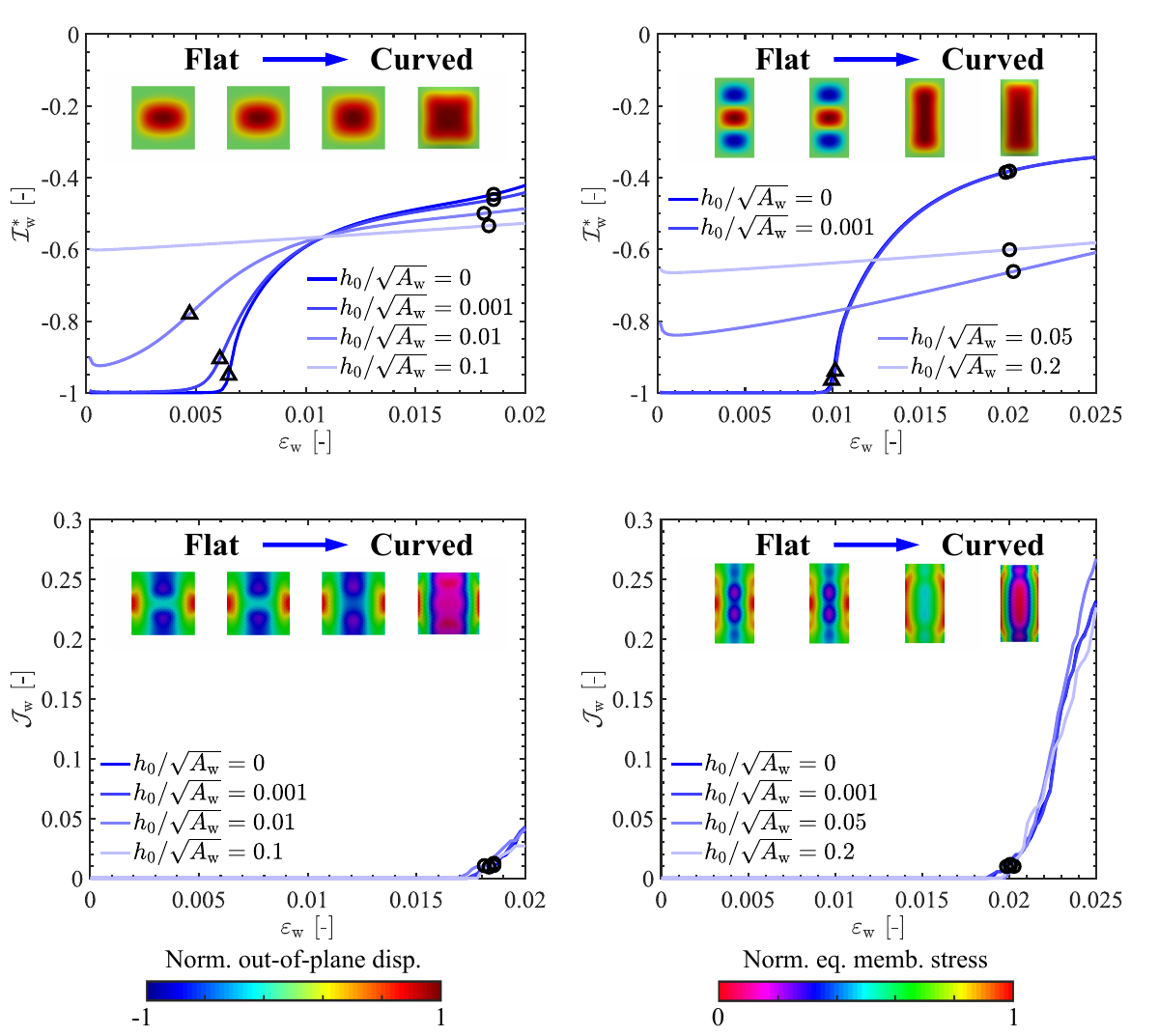}
\put(24.5,87.625){\color{black}\footnotesize \textbf{(a)}}
\put(73.625,87.625){\color{black}\footnotesize \textbf{(b)}}
\put(24.5,45.75){\color{black}\footnotesize \textbf{(c)}}
\put(73.625,45.75){\color{black}\footnotesize \textbf{(d)}}
\end{overpic}
\caption{\textcolor{red}{Cell wall strain energy partitioning indicators, membrane plasticity indicators, and bending patterns and membrane stress patterns at the yield points, of the rectangular cell wall models with different initial curvature levels $h_0/\sqrt{A_{\rm w}}$, for the two aspect ratios: (a, c) $\mathcal{R}_{\rm w}=1.0$ and (b, d) $\mathcal{R}_{\rm w}=2.0$. $(\bullet)^*$ in (a-b) indicates that the strain energy contribution by the perturbation force is removed. The buckling in (a-b) and yield points are indicated by the black triangles and circles, respectively}.}\label{fig:cell_wall}
\end{figure}

The strain energy partitioning indicators $\mathcal{I}_\stt{w}^*$ and membrane plasticity indicators $\mathcal{J}_\stt{w}$ are reported in Figure~\ref{fig:cell_wall}. It can be seen from Figures~\ref{fig:cell_wall}(a-b) that each flat cell wall deforms first by the membrane mode ($\mathcal{I}_\stt{w}^*\sim -1$) and then switches to a mixed membrane-bending mode after buckling (black triangles). For a given $\mathcal{R}_{\rm w}$, as $h_0/\sqrt{A_\stt{w}}$ increases, the bending contribution gets promoted (with a larger $\mathcal{I}_\stt{w}^*$) at the initial stage. This explains the compressive modulus reduction in Figure~\ref{fig:eff_response_cell_wall}. The deformation mode transition becomes less and less apparent, and it gets infeasible to detect buckling (without black triangles) for a large $h_0/\sqrt{A_\stt{w}}$. 

Figures~\ref{fig:cell_wall}(c-d) show that for a given $\mathcal{R}_{\rm w}$, the membrane plasticity evolution behaviours for different $h_0/\sqrt{A_\stt{w}}$ are quite similar.

In order to interpret the compressive strength change in Figure~\ref{fig:eff_response_cell_wall}, the bending patterns\footnote{\textcolor{red}{Since buckling becomes almost invisible for a large initial curvature, we adopt the term bending pattern for generality}.} and membrane stress patterns right after yielding are provided in Figure~\ref{fig:cell_wall}. For $\mathcal{R}_{\rm w}=1.0$, as $h_0/\sqrt{A_\stt{w}}$ increases, the bending pattern remains nearly unchanged as the critical buckling mode (one half wave) of a flat cell wall (see Figure~\ref{fig:cell_wall}(a)). However, the bending area increases, leading to a smaller load-carrying portion (see Figure~\ref{fig:cell_wall}(c)) and thus lower strength (see Figure~\ref{fig:eff_response_cell_wall}(a)). This is likely because the initial curvature pattern is compatible with the critical buckling mode. For $\mathcal{R}_{\rm w}=2.0$, as $h_0/\sqrt{A_\stt{w}}$ increases, the bending pattern first changes from the critical buckling mode (three half waves) of a flat cell wall to a higher-order mode (one half wave), which would require a higher compressive load to activate \cite{Gerard1957}, and then remains nearly unchanged, accompanied by increasing bending area (see Figure~\ref{fig:cell_wall}(b)). Attributed to the competition between the bending pattern change and bending area increase, the load-carrying portion increases first and then decreases (see Figure~\ref{fig:cell_wall}(d)), resulting in a complex  strength change (see Figure~\ref{fig:eff_response_cell_wall}(b)). This trend may be understood by that the initial curvature pattern is incompatible with the critical buckling mode but compatible with the higher-order mode. 

\begin{figure}[ht]
\centering
\begin{overpic}[draft=false,width=0.9\textwidth]{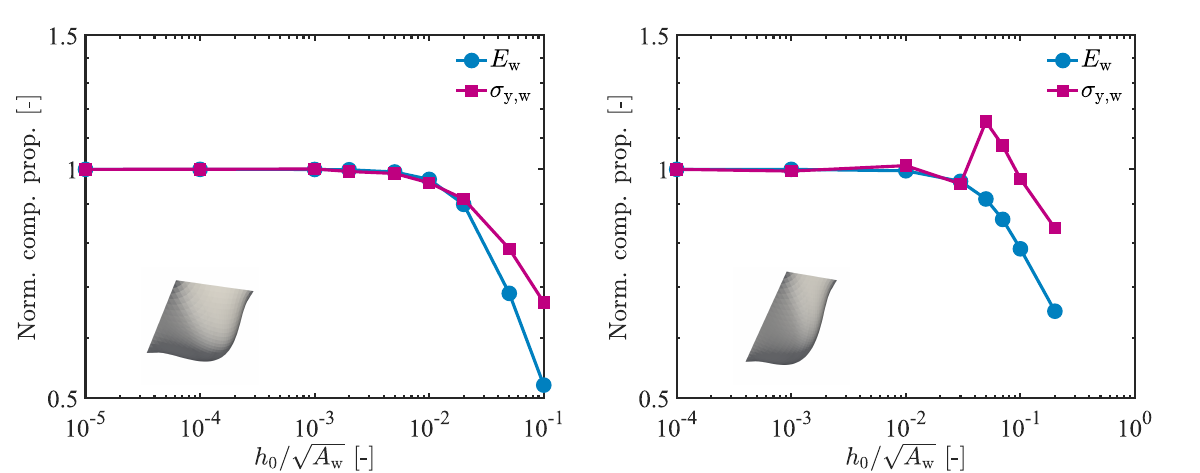}
\put(24.5,37.5){\color{black}\footnotesize \textbf{(a)}}
\put(73.625,37.5){\color{black}\footnotesize \textbf{(b)}}
\end{overpic}
\caption{\textcolor{red}{Compressive properties of the rectangular cell wall models with different initial curvature levels $h_0/\sqrt{A_{\rm w}}$, for the two aspect ratios: (a) $\mathcal{R}_{\rm w}=1.0$ and (b) $\mathcal{R}_{\rm w}=2.0$. Results have been normalized with respect to those at $h_0/\sqrt{A_{\rm w}}=0$}.}\label{fig:eff_props_cell_wall}
\end{figure}

Next, the compressive properties are extracted and plotted against varying $h_0/\sqrt{A_{\rm w}}$ in Figure~\ref{fig:eff_props_cell_wall}. For $\mathcal{R}_{\rm w}=1.0$, both the compressive modulus $E_\stt{w}$ and strength $\sigma_\stt{y,w}$ decrease with increasing $h_0/\sqrt{A_{\rm w}}$ (see Figure~\ref{fig:eff_props_cell_wall}(a)). Compared with $\sigma_\stt{y,w}$, $E_\stt{w}$ is more sensitive to $h_0/\sqrt{A_{\rm w}}$. For $\mathcal{R}_{\rm w}=2.0$, the compressive properties exhibit a rather complex trend against $h_0/\sqrt{A_{\rm w}}$ (see Figure~\ref{fig:eff_props_cell_wall}(b)). As $h_0/\sqrt{A_{\rm w}}$ increases, $E_\stt{w}$ decreases, while $\sigma_\stt{y,w}$ tends to increase first and then decrease. A critical $h_\stt{c}/\sqrt{A_{\rm w}}\sim0.01$ (corresponding to a normalized curvature $\sim0.1$ and see \cite{Andrews1999} for its definition) can be identified, below which the compressive properties are almost unaffected.

The cell wall curvature has been commonly regarded as a mesostructural feature which weakens the closed-foam compressive response (see e.g.\ \cite{Grenestedt1998, Simone1998b, Ghazi2019}). This preliminary study, however, has revealed that the initial curvature does not necessarily weaken but sometimes strengthens the cell wall compressive response (see also recent experimental evidence on lattice structures \cite{Ding2024}), depending on the cell wall aspect ratio and initial curvature level. The initial curvature pattern may play an important role as well. Therefore, it is believed that the impacts of cell wall curvature with corrugations and wriggles are more complex than what have been reported in the literature. Nevertheless, these impacts shall be negligible as long as the normalized curvature remains small ($< 0.1$, see also e.g.\ \cite{Grenestedt1998, Simone1998b}), which seems the case for most cell walls in many foams, especially polymer foams (see e.g.\ \cite{Andersons2016, Chai2020, Zhou2023a}). Accordingly, we have chosen to model each cell wall as a flat plate in the present study.
\color{black}

\section{Numerically realized mesostructural stochastics}\label{sec:meso_stochastics}
The stochastic variations of different mesostructural features for the tessellation-based models introduced in \textcolor{blue}{Section}~\ref{sec:num_model_method}, are elaborated in this appendix.

Detailed experimental characterization has been conducted in \cite{Zhou2023a} for Diab Divinycell foam H100 and H200. Both three-dimensional (3D) and 2D images obtained using X-ray CT scan and scanning electron microscope (SEM), respectively, are analysed. It is found that 3D and 2D measurements lead to similar distributions of cell equivalent diameter. However, cell wall thickness are largely overestimated with 3D measurements. Therefore, 2D measurements are adopted in the following. 

The cell equivalent diameters $d_\stt{v}$ follow a log-normal distribution \cite{Zhou2023a}:
\begin{equation}\label{eq:PDF_log_norm}
f(d_\stt{v})=\frac{1}{d_\stt{v}\sigma\sqrt{2\pi}}\exp\left(-\frac{(\ln{d_\stt{v}}-\mu)^2}{2\sigma^2}\right),
\end{equation}
where $\mu$ and $\sigma$ are the mean and standard deviation of $\ln d_\stt{v}$, respectively. $\mu$ and $\sigma$ can be related to the mean $\mu_d$ and standard deviation $\sigma_d$ of $d_\stt{v}$ through:
\begin{equation}\label{eq:PDF_log_norm_terms}
\mu=\ln\left(\frac{\mu_d^2}{\sqrt{\mu_d^2+\sigma_d^2}}\right),\,\,\sigma^2=\ln\left(1+\frac{\sigma_d^2}{\mu_d^2}\right).
\end{equation}

The cell wall thickness $t$ follows a gamma distribution \cite{Zhou2023a}:
\begin{equation}\label{eq:PDF_gamma}
f(t)=\frac{1}{\Gamma\theta^\alpha}t^{\alpha-1}\exp\left(-\frac{t}{\theta}\right),
\end{equation}
where $\alpha$ and $\theta$ denote the shape and scale parameters, respectively; $\Gamma$ is the gamma function, given by  $\Gamma(\alpha)=\int_0^\infty x^{\alpha-1}\exp(x)\, {\rm d}x$. $\alpha$ and $\theta$ can be related to the mean $\mu_t$ and standard deviation $\sigma_t$ of $t$ through:
\begin{equation}\label{eq:PDF_gamma_terms}
\alpha=\frac{\mu_t^2}{\sigma_t^2},\,\,\theta=\frac{\sigma_t^2}{\mu_t},
\end{equation}

Based on 2D measurements reported in \cite{Zhou2023a}, $(\mu_d,\sigma_d)$ and $(\mu_t,\sigma_t)$ for H100 and H200 are fitted, respectively, with the results listed in Table~\ref{tab:tess_geo_parameters}. No detailed measurements on individual cell shape anisotropy $\mathcal{R}_\stt{v}$ are provided in \cite{Zhou2023a}. Therefore, only the overall cell shape anisotropy $\mathcal{R}$ are given in Table~\ref{tab:tess_geo_parameters}.

\begin{figure}[t!]
\centering
\begin{overpic}[draft=false,width=0.9\textwidth]{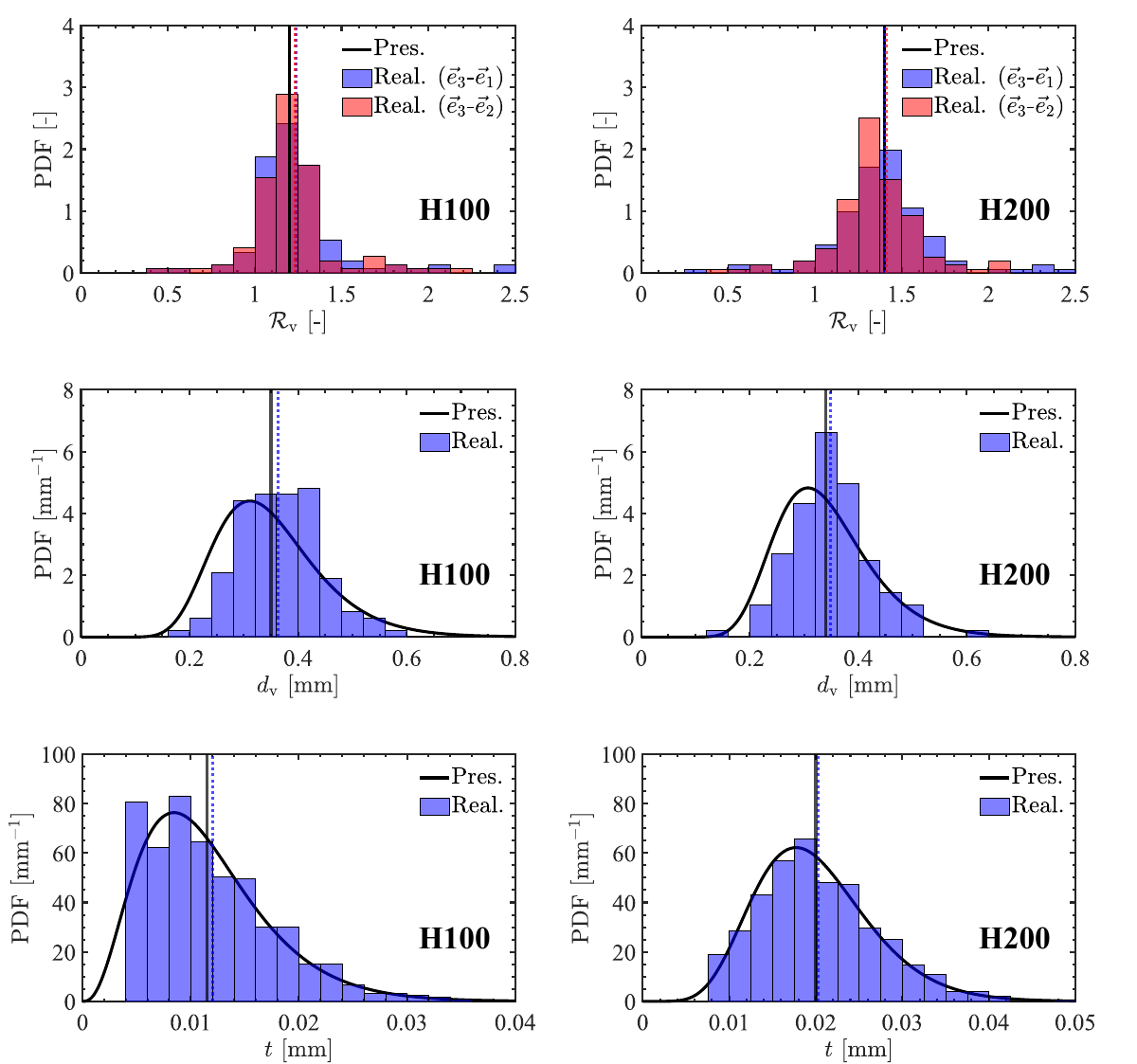}
\put(24.5,92.375){\color{black}\footnotesize \textbf{(a)}}
\put(73.625,92.375){\color{black}\footnotesize \textbf{(b)}}
\put(24.5,60.375){\color{black}\footnotesize \textbf{(c)}}
\put(73.625,60.375){\color{black}\footnotesize \textbf{(d)}}
\put(24.5,28.375){\color{black}\footnotesize \textbf{(e)}}
\put(73.625,28.375){\color{black}\footnotesize \textbf{(f)}}
\end{overpic}
\caption{Probability density functions (PDF) of different mesostructural features of the tessellation-based model set ``StSt" for two Divinycell foam grades: (a-b) cell shape anisotropy, (c-d) cell equivalent diameter and (e-f) cell wall thickness. Comparison between the prescribed and numerically realized distributions is shown. The means of the prescribed and numerically realized distributions are indicated by the solid and dashed vertical lines, respectively.}\label{fig:PDF_meso_foam}
\end{figure}
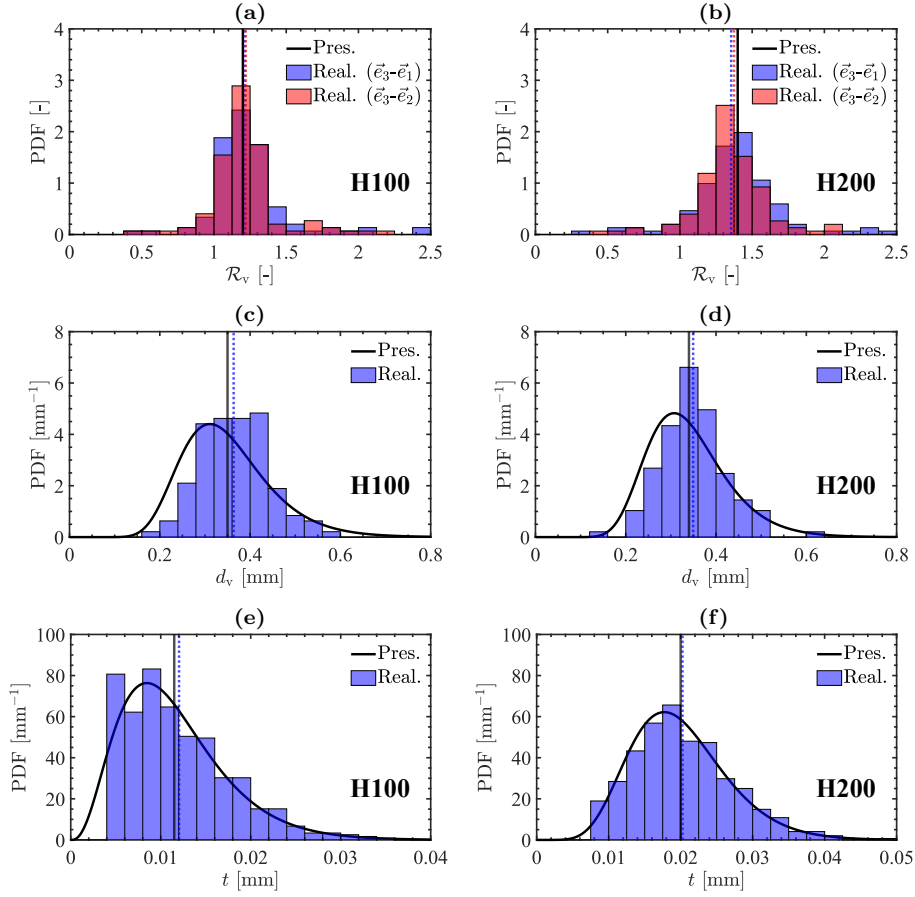

Probability density functions (PDF) of different mesostructural features of the generated model set ``StSt" for Divinycell foam H100 and H200, are compared to the prescribed ones in Figure~\ref{fig:PDF_meso_foam}. An excellent agreement can be observed between the prescribed and numerically realized distributions, indicating that ``StSt" can well approximate the real foam mesostructures. Notice that assigning an overall $\mathcal{R}$ would still cause varying $\mathcal{R}_\stt{v}$ of individual cells, because of the cell shape irregularity naturally induced by Laguerre tessellation. As expected, the cell aspect ratios $\mathcal{R}_\stt{v,31}$ and $\mathcal{R}_\stt{v,32}$ are quite comparable, and both approximately follow a normal distribution (see Figures~\ref{fig:PDF_meso_foam}(a-b)). \textcolor{red}{This trend is in qualitative agreement with more recent experimental measurements \cite{Skeens2024}}.

\section{Influence of the RVE size and random realization}\label{sec:influence_RVE_size_real}
The particular choice of RVE size and random realization may affect the macroscale effective responses shown in \textcolor{blue}{Section}~\ref{sec:analy_tess_meso_model}, and is thus examined in this appendix.

\begin{figure}[ht]
\centering
\begin{overpic}[draft=false,width=0.9\textwidth]{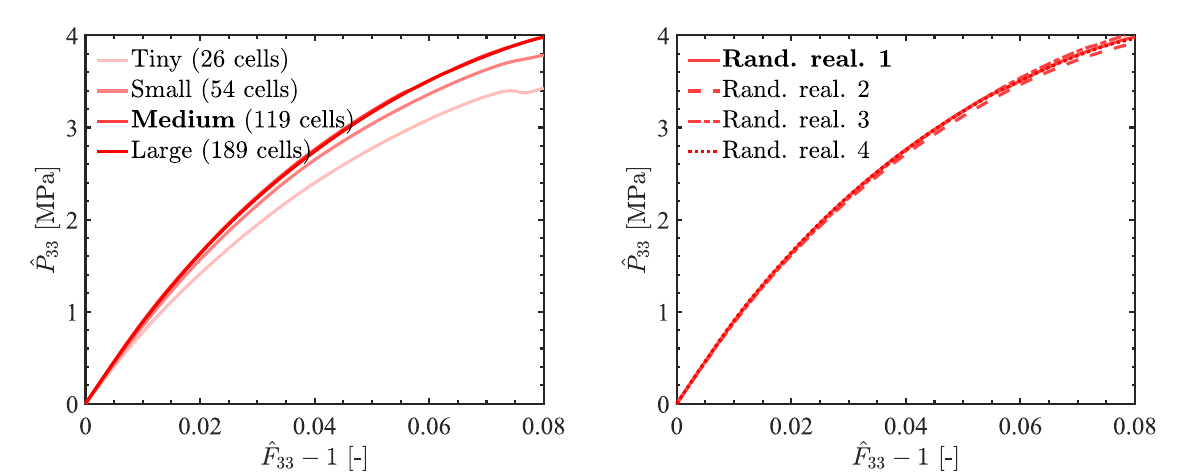}
\put(24.5,37.5){\color{black}\footnotesize \textbf{(a)}}
\put(73.625,37.5){\color{black}\footnotesize \textbf{(b)}}
\end{overpic}
\caption{(a) Effective stresses versus applied strain of the tessellation-based models ``StSt" H100 for (a) different RVE sizes and (b) different RVE random realizations, under uniaxial compression in the foam rise ($\vec{e}_3$) direction.
The curves for medium and large RVE sizes in (a) are nearly overlapping.}\label{fig:eff_props_foam_size_rand_real}
\end{figure}

The tessellation-based model ``StSt" for Divinycell foam H100 is focused on as one example. Four different choices of the RVE size $L_i$ are investigated, tiny 0.90 [mm], small 1.15 [mm], medium 1.50 [mm] and large 1.75 [mm]. Using the same random seed, RVE models consisting of 26, 54, 119 and 189 cells, respectively, are generated. The resulting overall relative densities $\rho/\rho_\stt{r}$ are 0.0682, 0.0742, 0.0806 and 0.0835, respectively. The effective stresses $\hat{\mathbf{P}}$ of the four RVE models under compression in the foam rise ($\vec{e}_3$) direction are reported in Figure~\ref{fig:eff_props_foam_size_rand_real}(a). It can be seen that as the RVE size increases, the effective stress response tends to increase (see also e.g.\ \cite{Vengatachalam2019, Zhou2023b, Ding2023}). This is likely because of the higher $\rho/\rho_\stt{r}$ associated with the larger RVE size. In addition, a sudden stress drop can be noticed for the ``tiny" size, which is likely caused by the high sensitivity to the presence of weak cell walls when the RVE size is too small. Nevertheless, the effective stress responses for the ``medium" and ``large" sizes are visually indistinguishable.

Next, the sensitivity to the RVE random realization is investigated. Four different random seeds with the ``medium" RVE size $L_i=1.50$ [mm] are considered. It has been verified that the resulting $\rho/\rho_\stt{r}$ are almost the same. The corresponding effective responses are reported in Figure~\ref{fig:eff_props_foam_size_rand_real}(b). It can be seen that different RVE random realizations deliver quite consistent results.

Similar findings have been confirmed for two transverse directions ($\vec{e}_1/\vec{e}_2$), which are thus not presented here. \textcolor{red}{The combination of ``medium" size and random realization ``1" has been adopted in the present study}.

\color{red}
\section{Influence of plasticity}\label{sec:influence_plas}

Focusing on low-density foams where the cell wall elastic buckling is the primary failure mode, the present study has disregarded plasticity in the material modelling. This simplification is examined in the following.

\begin{figure}[ht]
\centering
\begin{overpic}[draft=false,width=0.9\textwidth]{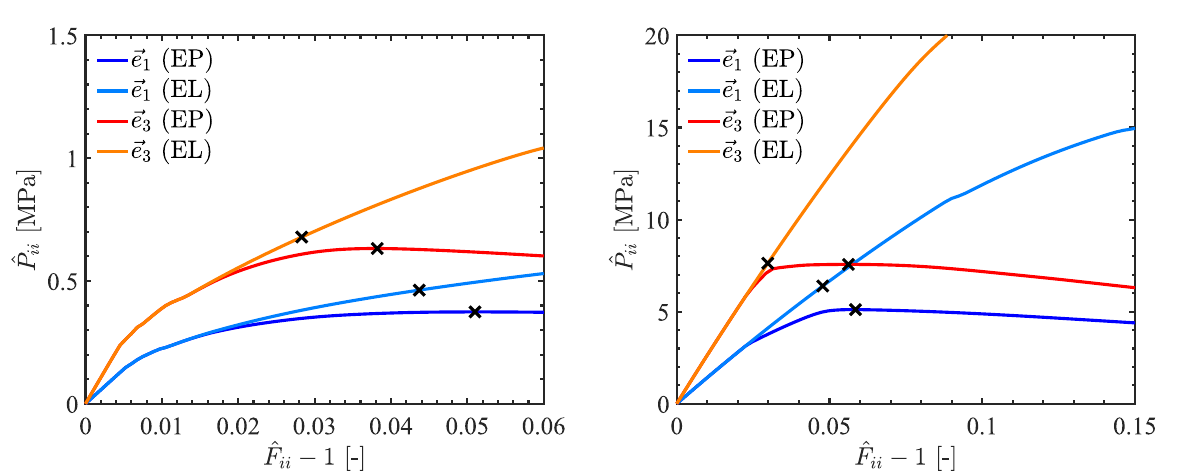}
\put(24.5,37.5){\color{black}\footnotesize \textbf{(a)}}
\put(73.625,37.5){\color{black}\footnotesize \textbf{(b)}}
\end{overpic}
\caption{\textcolor{red}{Effective responses of the Kelvin cell-based models, under uniaxial compression in the transverse ($\vec{e}_1$) and foam rise ($\vec{e}_3$) directions, for the two relative densities: (a) $\rho/\rho_\stt{r}=0.03$ and (c) $\rho/\rho_\stt{r}=0.15$. Comparison between the elasto-plastic (EP) and elastic (EL) models is shown. The yield points are indicated by the black crosses}.}\label{fig:eff_response_Kelvin_cell_05_T}
\end{figure}

The Kelvin cell-based model with cell shape anisotropy $\mathcal{R}=1.5$ (as described in \textcolor{blue}{Section}~\ref{sec:num_model_method}) is taken as an example, given its representativeness of many anisotropic foams. Model configurations with a broad range of cell wall thickness $t$ from 0.0036 [mm] to 0.018 [mm], are considered for the sake of comprehensive assessment. The resulting relative density $\rho/\rho_\stt{r}$ approximately varies from 0.03 to 0.15. The elasto-plastic numerical simulations are first conducted as the reference results. The effective stresses $\hat{\mathbf{P}}$ of the two relative densities $\rho/\rho_\stt{r}=0.03$ and $\rho/\rho_\stt{r}=0.15$ are plotted as functions of the applied strain in Figures~\ref{fig:eff_response_Kelvin_cell_05_T}(a) and (b), respectively. As expected, the typical plateau region observed in experiments are captured by an elasto-plastic model. 

The elastic numerical simulations are then conducted, with the results of $\rho/\rho_\stt{r}=0.03$ and $\rho/\rho_\stt{r}=0.15$ reported in Figures~\ref{fig:eff_response_Kelvin_cell_05_T}(a) and (b), respectively. Attributed to the absence of plasticity, the plateau regions are no more captured. A better agreement between the elasto-plastic and elastic models can be observed for $\rho/\rho_\stt{r}=0.03$ than $\rho/\rho_\stt{r}=0.15$. This can be understood by that as $\rho/\rho_\stt{r}$ decreases, the cell wall buckling plays a increasingly important role in governing the compressive response.

\begin{figure}[ht]
\centering
\begin{overpic}[draft=false,width=0.9\textwidth]{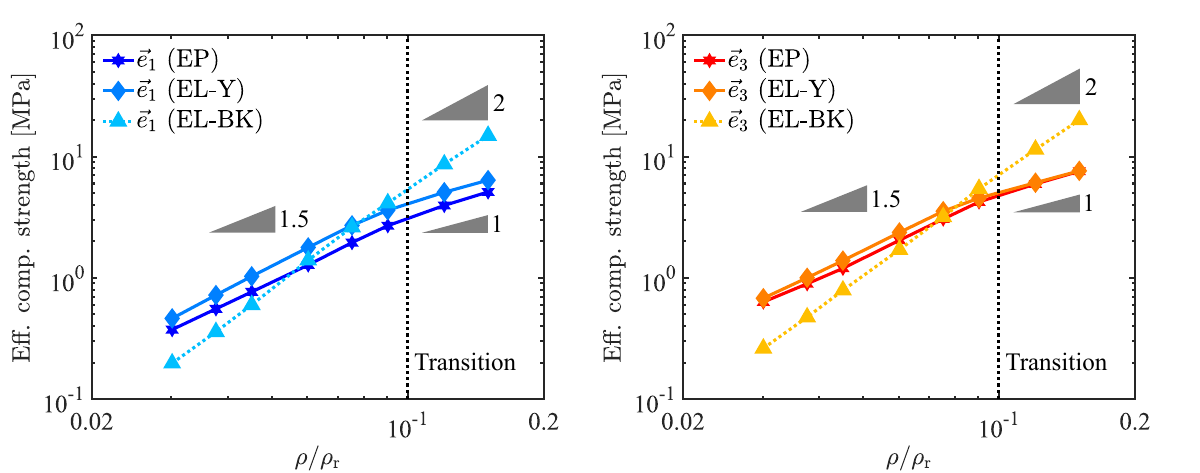}
\put(24.5,37.5){\color{black}\footnotesize \textbf{(a)}}
\put(73.625,37.5){\color{black}\footnotesize \textbf{(b)}}
\end{overpic}
\caption{\textcolor{red}{Effective compressive strengths of the Kelvin cell-based models for different relative densities $\rho/\rho_\stt{r}$, in the (a) transverse and (b) foam rise directions. Comparison between the elasto-plastic (EP) models, and elastic (EL) models in combination with the cell wall yielding (Y) and buckling (BK) criteria is shown. Several characteristic slopes in the logarithmic space are indicated by the grey triangles. The transition positions are indicated by the dashed vertical lines}.}\label{fig:eff_props_Kelvin_cell_wall_T}
\end{figure}

Next, the effective compressive strengths $\hat{\sigma}_{\rm y,11}$ and $\hat{\sigma}_{\rm y,33}$ in the transverse ($\vec{e}_1$) and foam rise ($\vec{e}_3$) directions, respectively, of the elasto-plastic and elastic models are extracted, and plotted against varying relative densities $\rho/\rho_\stt{r}$ in Figure~\ref{fig:eff_props_Kelvin_cell_wall_T}. Here, the yield strength of each elasto-plastic model is determined at the peak stress point (also indicated in Figure~\ref{fig:eff_response_Kelvin_cell_05_T}), as the reference result. The yield strength of each elastic model is determined at the first inclined cell wall yield point (as adopted in \textcolor{blue}{Section}~\ref{sec:analy_ideal_meso_model} and also indicated in Figure~\ref{fig:eff_response_Kelvin_cell_05_T}). The buckling strength is also provided, which is determined at the first inclined cell wall buckling point.

Figure~\ref{fig:eff_props_Kelvin_cell_wall_T} shows that for the elasto-plastic models, both $\hat{\sigma}_{\rm y,11}$ and $\hat{\sigma}_{\rm y,33}$ increase with increasing $\rho/\rho_\stt{r}$. For each loading direction, the slope is $\sim$ 1 at the higher-density regime, implying that the plastic collapse becomes the leading failure mode (see also \cite{Gibson1997_ch5}). At the lower-density regime, the slope is $\sim$ 1.5, implying that the elastic buckling followed by membrane yielding causes failure. A critical transition relative density of the failure mode can be identified as $\rho_\stt{c}/\rho_\stt{r}\sim 0.1$.

Both $\hat{\sigma}_{\rm y,11}$ and $\hat{\sigma}_{\rm y,33}$ of the elastic models against $\rho/\rho_\stt{r}$ demonstrate similar trends as those of the elasto-plastic models. $\hat{\sigma}_{\rm y,33}$ of the elastic models are almost identical to the elasto-plastic models, even when $\rho/\rho_\stt{r}>0.1$ (see Figure~\ref{fig:eff_props_Kelvin_cell_wall_T}(b)). $\hat{\sigma}_{\rm y,11}$ are overestimated by $\sim 20\%$ using the elastic models (see Figure~\ref{fig:eff_props_Kelvin_cell_wall_T}(a)), likely due to the larger inclined angle of primary load-bearing cell walls (see Figure~\ref{fig:geo_model_cell_wall}(b)) in the transverse direction, compared with the foam rise direction. Accordingly, the parallel cell walls, which tend to fail by early membrane yielding instead of buckling (see Figure~\ref{fig:cell_wall_Kelvin_cell_05}), may have more non-negligible impacts on $\hat{\sigma}_{\rm y,11}$ than $\hat{\sigma}_{\rm y,33}$. This effect has been disregarded in determining the compressive strengths of the elastic models.

The buckling strengths also increase as $\rho/\rho_\stt{r}$ increases. For each loading direction, the slope is nearly constant $\sim$ 2 (see also \cite{Gibson1997_ch5}). The buckling strength begins to exceed the yield strength approximately at $\rho/\rho_\stt{r}=0.09$, in good agreement with the identified $\rho_\stt{c}/\rho_\stt{r}\sim 0.1$ according to the slope change.

The above analyses confirm that for the Kelvin cells made from PVC, the critical transition relative density $\rho_\stt{c}/\rho_\stt{r}$ is $\sim$ 0.1, close to the value 0.11 suggested by experimental observations \cite{Kidd2012}. Notice that $\rho_\stt{c}/\rho_\stt{r}$ shall scale linearly with respect to the base material property ratio $\sqrt{\sigma_\stt{y}/E}$ (based on a rectangular plate analysis in \cite{Timoshenko1961}), and $\sqrt{\sigma_\stt{y}/E}=0.15$ for PVC. Taking PLA foams with $\sqrt{\sigma_\stt{y}/E}=0.19$ as another example, $\rho_\stt{c}/\rho_\stt{r}$ is estimated as $\sim$ 0.125, again in good agreement with the value 0.14 suggested by experimental observations \cite{Duan2019}. $\rho_\stt{c}/\rho_\stt{r}$ seems less well-established for metallic foams \cite{Michailidis2011, Kader2017}, likely because of more complicated mesostructural features, e.g.\ with more defects and larger stochastic variations.

In addition, the compressive strength ratios $\hat{\sigma}_{\stt{y},33}/\hat{\sigma}_{\stt{y},11}$ are discussed. It is found that $\hat{\sigma}_{\stt{y},33}/\hat{\sigma}_{\stt{y},11}$ tends to decrease with increasing $\rho/\rho_\stt{r}$ (see Figure~\ref{fig:eff_props_Kelvin_cell_wall_T}). Nevertheless, $\hat{\sigma}_{\stt{y},33}/\hat{\sigma}_{\stt{y},11}$ only varies $\sim \pm10\%$ with respect to that at $\rho/\rho_\stt{r}=0.075$, which has been chosen for the Kelvin cells in the present study. 

The effective compressive moduli $\hat{E}_{11}$ and $\hat{E}_{33}$ of the elasto-plastic and elastic models are identical to each other and thus not reported here. The compressive modulus ratios $\hat{E}_{33}/\hat{E}_{11}$ are nearly constant for all the considered $\rho/\rho_\stt{r}$. This is expected since the initial elastic region is dominated by the membrane deformation mode.
\color{black}

\section{Reference analytical models}\label{sec:ref_anal_models}
Two reference analytical models assessed in \textcolor{blue}{Sections}~\ref{sec:relation_anisotropy} and~\ref{sec:analy_tess_meso_model} for the predictions of foam mechanical anisotropy under compression, are detailed in this appendix.

The first one is the Gibson-Ashby model \cite{Gibson1997_ch6}, which is derived by adopting a rectangular parallelepiped cell structure. The cell edge bending deformation accompanied by the face tension along the direction perpendicular to the compressive loading, is assumed to govern the initial elastic region. The cell wall plastic collapse is assumed to be the leading failure mode. The mechanical anisotropy $\mathcal{R}^E$ and $\mathcal{R}^\sigma$ are expressed in terms of shape anisotropy $\mathcal{R}$ as:
\begin{subequations}\label{eq:eff_prop_ratios_GA}
	\begin{align}
	&\mathcal{R}^E=\phi\frac{2\mathcal{R}^2}{1+\mathcal{R}^\stt{-3}}+(1-\phi)\frac{2\mathcal{R}}{1+\mathcal{R}^\stt{-1}},\\
	&\mathcal{R}^\sigma=\frac{2\mathcal{R}}{1+\mathcal{R}^\stt{-1}},
	\end{align}
\end{subequations}
where $\phi$ denotes the cell edge fraction. Note, that the physical interpretation of $\phi$ has been relaxed, which is instead treated as a fitting parameter in the Gibson-Ashby model. Interestingly, the fitted $\phi$ are found to be much higher than the experimental data for some closed-cell foams (see e.g.\ \cite{Mills1999, Almanza2001, Mills2009}).

For consistency with numerical simulations, $\phi=0$ has been used in the present study. The results predicted using eq.~\eqref{eq:eff_prop_ratios_GA} are reported in Figure~\ref{fig:eff_prop_ratios_cell}(a).

The second one is the Sullivan model \cite{Sullivan2008}, which is derived by adopting a (open-cell) Kelvin cell structure. The cell edge axial and bending deformations are both taken into account. The cell edge axial failure after reaching the base material yield stress is assumed to be the foam failure mechanism. The resulting $\mathcal{R}^E$ and $\mathcal{R}^\sigma$ are expressed as:
\begin{subequations}\label{eq:eff_prop_ratios_S}
	\begin{align}
	&\mathcal{R}^E=\frac{\mathcal{R}^2}{4}\frac{\mathcal{A}_1\mathcal{A}_2}{\mathcal{A}_3}\\
	&\mathcal{R}^\sigma=\mathcal{R}\frac{\mathcal{B}_1\mathcal{B}_2}{\mathcal{B}_3},
	\end{align}
\end{subequations}
with the intermediate terms
\begin{subequations}\label{eq:eff_prop_ratios_S_terams_A}
	\begin{align}
	&\mathcal{A}_1=C_1\left(2\tilde{Q}^2\mathcal{R}^2+\frac{64Q^3}{\mathcal{R}_1^{\frac{1}{2}}}\right),\\
	&\mathcal{A}_2=\gamma C_2\frac{8\tilde{Q}^3(32+4Q\mathcal{R}_1^{\frac{1}{2}})\mathcal{R}}{\mathcal{R}_1\mathcal{R}_2},\\
	&\mathcal{A}_3=16C_1+\gamma C_2\frac{8\tilde{Q}^5\mathcal{R}^3}{\mathcal{R}_1\mathcal{R}_2},
	\end{align}
\end{subequations}
and 
\begin{subequations}\label{eq:eff_prop_ratios_S_terams_B}
	\begin{align}
    &\mathcal{B}_1=\sqrt{C_1}\tilde{Q}\mathcal{R},\\
    &\mathcal{B}_2=\sqrt{\gamma}C_3\frac{16\sqrt{2}\tilde{Q}^{\frac{3}{2}}\mathcal{R}^{\frac{1}{2}}}{\mathcal{R}_1^{\frac{1}{2}}\mathcal{R}_2^{\frac{1}{2}}},\\
    &\mathcal{B}_3=4\sqrt{C_1}+\sqrt{\gamma}C_3\frac{4\sqrt{2}\tilde{Q}^{\frac{5}{2}}\mathcal{R}^{\frac{3}{2}}}{\mathcal{R}_1^{\frac{1}{2}}\mathcal{R}_2^{\frac{1}{2}}},
	\end{align}
\end{subequations}
respectively, where two additional ratios $\mathcal{R}_1$ and $\mathcal{R}_2$ have been introduced:
\begin{equation}\label{eq:eff_prop_ratios_S_R}
	\mathcal{R}_1=16+\tilde{Q}^2\mathcal{R}^2,\,\,\mathcal{R}_2=4Q+2(16+\tilde{Q}^2\mathcal{R}^2)^{\frac{1}{2}}.
\end{equation}
Here, $Q$ denotes a Kelvin cell shape parameter and $\tilde{Q}=2+\sqrt{2}Q$; $\gamma=\rho/\rho_\stt{r}$ is the relative density; $C_1$, $C_2$ and $C_3$ are three constants which characterize the cell edge cross-section shape. 

Considering a standard Kelvin cell shape and three-cusp hypocycloid cross-section shape of cell edges, leads to $Q=\sqrt{2}$, and $C_1=\sqrt{3}-\frac{\pi}{2}$, $C_2=\frac{20\sqrt{3}-11\pi}{2\sqrt{3}-\pi}$ and $C_2=\frac{60-11\sqrt{3}\pi}{24(\sqrt{3}-\frac{\pi}{2})}$ \cite{Sullivan2008}. The results predicted using eq.~\eqref{eq:eff_prop_ratios_S} are reported in Figures~\ref{fig:eff_prop_ratios_cell}(b),~\ref{fig:eff_prop_ratios_foam_2}(a) and~\ref{fig:eff_prop_ratios_foam_2}(b).

\bibliographystyle{unsrt}
\bibliography{References_Lei_foam}

\begin{thebibliography}{100}

\bibitem{Smith2012}
B.H. Smith, S.~Szyniszewski, J.F. Hajjar, B.W. Schafer, and S.R. Arwade.
\newblock {Steel foam for structures: A review of applications, manufacturing
  and material properties}.
\newblock {\em Journal of Constructional Steel Research}, 71:1--10, 2012.

\bibitem{Sun2018}
Y.L. Sun and Q.M. Li.
\newblock {Dynamic compressive behaviour of cellular materials: A review of
  phenomenon, mechanism and modelling}.
\newblock {\em International Journal of Impact Engineering}, 112:74--115, 2018.

\bibitem{Rahimidehgolan2023}
F.~Rahimidehgolan and W.~Altenhof.
\newblock {Compressive behavior and deformation mechanisms of rigid polymeric
  foams: A review}.
\newblock {\em Composites Part B: Engineering}, 253:110513, 2023.

\bibitem{Mu2010}
Y.~Mu, G.~Yao, and H.~Luo.
\newblock {Effect of cell shape anisotropy on the compressive behavior of
  closed-cell aluminum foams}.
\newblock {\em Materials {\&} Design}, 31(3):1567--1569, 2010.

\bibitem{Zhou2023a}
Y.~Zhou, B.~Xue, W.~Zhang, and R.~Wang.
\newblock {Prediction of bulk mechanical properties of PVC foam based on
  microscopic model: Part I-Microstructure characterization and generation
  algorithm}.
\newblock {\em Polymer Testing}, 117:107872, 2023.

\bibitem{Jang2015}
W.Y. Jang, W.Y. Hsieh, C.C. Miao, and Y.C. Yen.
\newblock {Microstructure and mechanical properties of ALPORAS closed-cell
  aluminium foam}.
\newblock {\em Materials Characterization}, 107:228--238, 2015.

\bibitem{Tang2022}
Y.~Tang, W.~Zhang, X.~Jiang, J.~Zhao, W.~Xie, and T.~Chen.
\newblock {Experimental investigations on phenomenological constitutive model
  of closed-cell PVC foam considering the effects of density, strain rate and
  anisotropy}.
\newblock {\em Composites Part B: Engineering}, 238:109885, 2022.

\bibitem{Gibson1997_ch5}
L.J. Gibson and M.F. Ashby.
\newblock {The mechanics of foams: basic results}.
\newblock In {\em Cellular Solids}, pages 175--234. Cambridge University Press,
  Cambridge, 1997.

\bibitem{Andrews1999}
E.~Andrews, W.~Sanders, and L.J. Gibson.
\newblock {Compressive and tensile behaviour of aluminum foams}.
\newblock {\em Materials Science and Engineering: A}, 270(2):113--124, 1999.

\bibitem{Jeon2005}
I.~Jeon and T.~Asahina.
\newblock {The effect of structural defects on the compressive behavior of
  closed-cell Al foam}.
\newblock {\em Acta Materialia}, 53(12):3415--3423, 2005.

\bibitem{PerezTamarit2019}
S.~P{\'{e}}rez-Tamarit, E.~Sol{\'{o}}rzano, A.~Hilger, I.~Manke, and M.A.
  Rodr{\'{i}}guez-P{\'{e}}rez.
\newblock {Effect of solid phase corrugation on the thermo-mechanical
  properties of low density flexible cellular polymers}.
\newblock {\em Materials {\&} Design}, 161:106--113, 2019.

\bibitem{Ghazi2020a}
A.~Ghazi, P.~Berke, C.~Tiago, and T.J. Massart.
\newblock {Computed tomography based modelling of the behaviour of closed cell
  metallic foams using a shell approximation}.
\newblock {\em Materials {\&} Design}, 194:108866, 2020.

\bibitem{Andersons2016}
J.~Andersons, M.~Kirpluks, L.~Stiebra, and U.~Cabulis.
\newblock {Anisotropy of the stiffness and strength of rigid low-density
  closed-cell polyisocyanurate foams}.
\newblock {\em Materials {\&} Design}, 92:836--845, 2016.

\bibitem{Tan2005}
P.J. Tan, S.R. Reid, J.J. Harrigan, Z.~Zou, and S.~Li.
\newblock {Dynamic compressive strength properties of aluminium foams. Part
  I-experimental data and observations}.
\newblock {\em Journal of the Mechanics and Physics of Solids},
  53(10):2174--2205, 2005.

\bibitem{Zenkert2009}
D.~Zenkert and M.~Burman.
\newblock {Tension, compression and shear fatigue of a closed cell polymer
  foam}.
\newblock {\em Composites Science and Technology}, 69(6):785--792, 2009.

\bibitem{Michailidis2011}
N.~Michailidis, F.~Stergioudi, and A.~Tsouknidas.
\newblock {Deformation and energy absorption properties of powder-metallurgy
  produced Al foams}.
\newblock {\em Materials Science and Engineering: A}, 528(24):7222--7227, 2011.

\bibitem{Kidd2012}
T.H. Kidd, S.~Zhuang, and G.~Ravichandran.
\newblock {In situ mechanical characterization during deformation of PVC
  polymeric foams using ultrasonics and digital image correlation}.
\newblock {\em Mechanics of Materials}, 55:82--88, 2012.

\bibitem{Koohbor2018}
B.~Koohbor, S.~Ravindran, and A.~Kidane.
\newblock {Effects of cell-wall instability and local failure on the response
  of closed-cell polymeric foams subjected to dynamic loading}.
\newblock {\em Mechanics of Materials}, 116:67--76, 2018.

\bibitem{Duan2019}
Y.~Duan, B.~Du, X.~Shi, B.~Hou, and Y.~Li.
\newblock {Quasi-static and dynamic compressive properties and deformation
  mechanisms of 3D printed polymeric cellular structures with Kelvin cells}.
\newblock {\em International Journal of Impact Engineering}, 132:103303, 2019.

\bibitem{Kader2017}
M.A. Kader, M.A. Islam, M.~Saadatfar, P.J. Hazell, A.D. Brown, S.~Ahmed, and
  J.P. Escobedo.
\newblock {Macro and micro collapse mechanisms of closed-cell aluminium foams
  during quasi-static compression}.
\newblock {\em Materials {\&} Design}, 118:11--21, 2017.

\bibitem{Benouali2005}
A.H. Benouali, L.~Froyen, T.~Dillard, S.~Forest, and F.~N'guyen.
\newblock {Investigation on the influence of cell shape anisotropy on the
  mechanical performance of closed cell aluminium foams using micro-computed
  tomography}.
\newblock {\em Journal of Materials Science}, 40(22):5801--5811, 2005.

\bibitem{Bafti2013}
H.~Bafti and A.~Habibolahzadeh.
\newblock {Compressive properties of aluminum foam produced by powder-Carbamide
  spacer route}.
\newblock {\em Materials {\&} Design}, 52:404--411, 2013.

\bibitem{Cheng2018}
Y.~Cheng, Y.~Li, X.~Chen, X.~Zhou, and N.~Wang.
\newblock {Compressive properties and energy absorption of aluminum foams with
  a wide range of relative densities}.
\newblock {\em Journal of Materials Engineering and Performance},
  27(8):4016--4024, 2018.

\bibitem{Saha2005}
M.C. Saha, H.~Mahfuz, U.K. Chakravarty, M.~Uddin, Md.~E. Kabir, and S.~Jeelani.
\newblock {Effect of density, microstructure, and strain rate on compression
  behavior of polymeric foams}.
\newblock {\em Materials Science and Engineering: A}, 406(1-2):328--336, 2005.

\bibitem{Zhou2023b}
Y.~Zhou, B.~Xue, W.~Zhang, and R.~Wang.
\newblock {Prediction of bulk mechanical properties of PVC foam based on
  microscopic model: Part II-Material characterization and analytical
  formulae}.
\newblock {\em Polymer Testing}, 117:107846, 2023.

\bibitem{Mondal2009}
D.P. Mondal, M.D. Goel, and S.~Das.
\newblock {Compressive deformation and energy absorption characteristics of
  closed cell aluminum-fly ash particle composite foam}.
\newblock {\em Materials Science and Engineering: A}, 507(1-2):102--109, 2009.

\bibitem{Guo2015}
C.~Guo, T.~Zou, C.~Shi, X.~Yang, N.~Zhao, E.~Liu, and C.~He.
\newblock {Compressive properties and energy absorption of aluminum composite
  foams reinforced by in-situ generated MgAl 2 O 4 whiskers}.
\newblock {\em Materials Science and Engineering: A}, 645:1--7, 2015.

\bibitem{Kim2005}
Y.W. Kim, Y.J. Jin, Y.S. Chun, I.H. Song, and H.D. Kim.
\newblock {A simple pressing route to closed-cell microcellular ceramics}.
\newblock {\em Scripta Materialia}, 53(8):921--925, 2005.

\bibitem{Celzard2010}
A.~Celzard, W.~Zhao, A.~Pizzi, and V.~Fierro.
\newblock {Mechanical properties of tannin-based rigid foams undergoing
  compression}.
\newblock {\em Materials Science and Engineering: A}, 527(16-17):4438--4446,
  2010.

\bibitem{Deshpande2001}
V.S. Deshpande and N.A. Fleck.
\newblock {Multi-axial yield behaviour of polymer foams}.
\newblock {\em Acta Materialia}, 49(10):1859--1866, 2001.

\bibitem{EdwinRaj2009}
R.~{Edwin Raj} and B.S.S. Daniel.
\newblock {Structural and compressive property correlation of closed-cell
  aluminum foam}.
\newblock {\em Journal of Alloys and Compounds}, 467(1-2):550--556, 2009.

\bibitem{Linul2018}
E.~Linul, N.~Movahedi, and L.~Marsavina.
\newblock {The temperature and anisotropy effect on compressive behavior of
  cylindrical closed-cell aluminum-alloy foams}.
\newblock {\em Journal of Alloys and Compounds}, 740:1172--1179, 2018.

\bibitem{Liu2020}
Y.~Liu, F.~Rahimidehgolan, and W.~Altenhof.
\newblock {Anisotropic compressive behavior of rigid PVC foam at strain rates
  up to 200 s$^{-1}$}.
\newblock {\em Polymer Testing}, 91:106836, 2020.

\bibitem{Hamilton2013}
A.R. Hamilton, O.T. Thomsen, L.A.O. Madaleno, L.R. Jensen, J.M. Rauhe, and
  R.~Pyrz.
\newblock {Evaluation of the anisotropic mechanical properties of reinforced
  polyurethane foams}.
\newblock {\em Composites Science and Technology}, 87:210--217, 2013.

\bibitem{MarviMashhadi2018}
M.~Marvi-Mashhadi, C.S. Lopes, and J.~LLorca.
\newblock {Effect of anisotropy on the mechanical properties of polyurethane
  foams: An experimental and numerical study}.
\newblock {\em Mechanics of Materials}, 124:143--154, 2018.

\bibitem{Li2019}
P.~Li, Y.B. Guo, M.W. Zhou, and V.P.W. Shim.
\newblock {Response of anisotropic polyurethane foam to compression at
  different loading angles and strain rates}.
\newblock {\em International Journal of Impact Engineering}, 127:154--168,
  2019.

\bibitem{Linul2013}
E.~Linul, L.~Marsavina, T.~Voiconi, and T.~Sadowski.
\newblock {Study of factors influencing the mechanical properties of
  polyurethane foams under dynamic compression}.
\newblock {\em Journal of Physics: Conference Series}, 451:012002, 2013.

\bibitem{Gibson1997_ch6}
L.J. Gibson and M.F. Ashby.
\newblock {The mechanics of foams: refinements}.
\newblock In {\em Cellular Solids}, pages 235--282. Cambridge University Press,
  Cambridge, 1997.

\bibitem{Gong2005}
L.~Gong, S.~Kyriakides, and W.Y. Jang.
\newblock {Compressive response of open-cell foams. Part I: Morphology and
  elastic properties}.
\newblock {\em International Journal of Solids and Structures},
  42(5-6):1355--1379, 2005.

\bibitem{Sullivan2008}
R.M. Sullivan, L.J. Ghosn, and B.A. Lerch.
\newblock {A general tetrakaidecahedron model for open-celled foams}.
\newblock {\em International Journal of Solids and Structures},
  45(6):1754--1765, 2008.

\bibitem{EspadasEscalante2015}
J.J. Espadas-Escalante and F.~Avil{\'{e}}s.
\newblock {Anisotropic compressive properties of multiwall carbon
  nanotube/polyurethane foams}.
\newblock {\em Mechanics of Materials}, 91:167--176, 2015.

\bibitem{Doyle2019}
L.~Doyle, I.~Weidlich, and M.~Illguth.
\newblock {Anisotropy in polyurethane pre-insulated pipes}.
\newblock {\em Polymers}, 11(12):2074, 2019.

\bibitem{Poapongsakorn2011}
P.~Poapongsakorn and C.~Kanchanomai.
\newblock {Time-dependent deformation of closed-cell PVC foam}.
\newblock {\em Journal of Cellular Plastics}, 47(4):323--336, 2011.

\bibitem{Luong2013}
D.D. Luong, D.~Pinisetty, and N.~Gupta.
\newblock {Compressive properties of closed-cell polyvinyl chloride foams at
  low and high strain rates: Experimental investigation and critical review of
  state of the art}.
\newblock {\em Composites Part B: Engineering}, 44(1):403--416, 2013.

\bibitem{Concas2019}
F.~Concas, S.~Diebels, and A.~Jung.
\newblock {Multiaxial failure surface of PVC foams and monitoring of
  deformation bands by three-dimensional digital image correlation}.
\newblock {\em Journal of the Mechanics and Physics of Solids}, 130:195--215,
  2019.

\bibitem{Bolintineanu2021}
D.S. Bolintineanu, R.~Waymel, H.~Collis, K.N. Long, E.C. Quintana, and S.L.B.
  Kramer.
\newblock {Anisotropy evolution of elastomeric foams during uniaxial
  compression measured via in-situ X-ray computed tomography}.
\newblock {\em Materialia}, 18:101112, 2021.

\bibitem{Chai2020}
H.W. Chai, Z.L. Xie, X.H. Xiao, H.L. Xie, J.Y. Huang, and S.N. Luo.
\newblock {Microstructural characterization and constitutive modeling of
  deformation of closed-cell foams based on in situ x-ray tomography}.
\newblock {\em International Journal of Plasticity}, 131:102730, 2020.

\bibitem{Santosa1998}
S.~Santosa and T.~Wierzbicki.
\newblock {On the modeling of crush behavior of a closed-cell aluminum foam
  structure}.
\newblock {\em Journal of the Mechanics and Physics of Solids}, 46(4):645--669,
  1998.

\bibitem{Simone1998a}
A.E. Simone and L.J. Gibson.
\newblock {Effects of solid distribution on the stiffness and strength of
  metallic foams}.
\newblock {\em Acta Materialia}, 46(6):2139--2150, 1998.

\bibitem{Simone1998b}
A.E. Simone and L.J. Gibson.
\newblock {The effects of cell face curvature and corrugations on the stiffness
  and strength of metallic foams}.
\newblock {\em Acta Materialia}, 46(11):3929--3935, 1998.

\bibitem{Grenestedt1998}
J.L. Grenestedt.
\newblock {Influence of wavy imperfections in cell walls on elastic stiffness
  of cellular solids}.
\newblock {\em Journal of the Mechanics and Physics of Solids}, 46(1):29--50,
  1998.

\bibitem{Grenestedt2000}
J.L. Grenestedt and F.~Bassinet.
\newblock {Influence of cell wall thickness variations on elastic stiffness of
  closed-cell cellular solids}.
\newblock {\em International Journal of Mechanical Sciences}, 42(7):1327--1338,
  2000.

\bibitem{DeGiorgi2010}
M.~{De Giorgi}, A.~Carofalo, V.~Dattoma, R.~Nobile, and F.~Palano.
\newblock {Aluminium foams structural modelling}.
\newblock {\em Computers {\&} Structures}, 88(1-2):25--35, 2010.

\bibitem{Sadek2013}
E.~Sadek and N.A. Fouad.
\newblock {Finite element modeling of compression behavior of extruded
  polystyrene foam using X-ray tomography}.
\newblock {\em Journal of Cellular Plastics}, 49(2):161--191, 2013.

\bibitem{Chen2018}
Y.~Chen, R.~Das, M.~Battley, and Z.~Xu.
\newblock {Compressive and shear strengths of the ductile closed-cell Kelvin
  and Weaire-Phelan foams along the lattice direction [100]}.
\newblock {\em Thin-Walled Structures}, 132:237--249, 2018.

\bibitem{Shakibanezhad2022}
R.~Shakibanezhad, M.~Sadighi, and R.~Hedayati.
\newblock {Numerical and experimental study of quasi-static loading of aluminum
  closed-cell foams using weaire–phelan and kelvin tessellations}.
\newblock {\em Transport in Porous Media}, 142(1-2):229--248, 2022.

\bibitem{Caty2008}
O.~Caty, E.~Maire, S.~Youssef, and R.~Bouchet.
\newblock {Modeling the properties of closed-cell cellular materials from
  tomography images using finite shell elements}.
\newblock {\em Acta Materialia}, 56(19):5524--5534, 2008.

\bibitem{Jeon2010}
I.~Jeon, T.~Asahina, K.J. Kang, S.~Im, and T.~Lu.
\newblock {Finite element simulation of the plastic collapse of closed-cell
  aluminum foams with X-ray computed tomography}.
\newblock {\em Mechanics of Materials}, 42(3):227--236, 2010.

\bibitem{Sulong2015}
M.A. Sulong, M.~Taherishargh, I.V. Belova, G.E. Murch, and T.~Fiedler.
\newblock {On the mechanical anisotropy of the compressive properties of
  aluminium perlite syntactic foam}.
\newblock {\em Computational Materials Science}, 109:258--265, 2015.

\bibitem{Natesaiyer2015}
K.~Natesaiyer, C.~Chan, S.~Sinha-Ray, D.~Song, C.L. Lin, J.D. Miller, E.J.
  Garboczi, and A.M. Forster.
\newblock {X-ray CT imaging and finite element computations of the elastic
  properties of a rigid organic foam compared to experimental measurements:
  insights into foam variability}.
\newblock {\em Journal of Materials Science}, 50(11):4012--4024, 2015.

\bibitem{Sun2017}
Y.~Sun, X.~Zhang, Z.~Shao, and Q.M. Li.
\newblock {Image-based correlation between the meso-scale structure and
  deformation of closed-cell foam}.
\newblock {\em Materials Science and Engineering: A}, 688:27--39, 2017.

\bibitem{Chen2017b}
Y.~Chen, R.~Das, and M.~Battley.
\newblock {Finite element analysis of the compressive and shear responses of
  structural foams using computed tomography}.
\newblock {\em Composite Structures}, 159:784--799, 2017.

\bibitem{Talebi2019}
S.~Talebi, M.~Sadighi, and M.M. Aghdam.
\newblock {Numerical and experimental analysis of the closed-cell aluminium
  foam under low velocity impact using computerized tomography technique}.
\newblock {\em Acta Mechanica Sinica}, 35(1):144--155, 2019.

\bibitem{Song2010}
Y.~Song, Z.~Wang, L.~Zhao, and J.~Luo.
\newblock {Dynamic crushing behavior of 3D closed-cell foams based on Voronoi
  random model}.
\newblock {\em Materials {\&} Design}, 31(9):4281--4289, 2010.

\bibitem{Shi2018}
X.~Shi, S.~Liu, H.~Nie, G.~Lu, and Y.~Li.
\newblock {Study of cell irregularity effects on the compression of closed-cell
  foams}.
\newblock {\em International Journal of Mechanical Sciences}, 135:215--225,
  2018.

\bibitem{Vengatachalam2019}
B.~Vengatachalam, L.H. Poh, Z.S. Liu, Q.H. Qin, and S.~Swaddiwudhipong.
\newblock {Three dimensional modelling of closed-cell aluminium foams with
  predictive macroscopic behaviour}.
\newblock {\em Mechanics of Materials}, 136:103067, 2019.

\bibitem{Roberts2001}
A.P. Roberts and E.J. Garboczi.
\newblock {Elastic moduli of model random three-dimensional closed-cell
  cellular solids}.
\newblock {\em Acta Materialia}, 49(2):189--197, 2001.

\bibitem{Koll2016}
J.~K{\"{o}}ll and S.~Hallstr{\"{o}}m.
\newblock {Elastic properties of equilibrium foams}.
\newblock {\em Acta Materialia}, 113:11--18, 2016.

\bibitem{Chen2015}
Y.~Chen, R.~Das, and M.~Battley.
\newblock {Effects of cell size and cell wall thickness variations on the
  stiffness of closed-cell foams}.
\newblock {\em International Journal of Solids and Structures}, 52:150--164,
  2015.

\bibitem{Chen2017a}
Y.~Chen, R.~Das, and M.~Battley.
\newblock {Effects of cell size and cell wall thickness variations on the
  strength of closed-cell foams}.
\newblock {\em International Journal of Engineering Science}, 120:220--240,
  2017.

\bibitem{MarviMashhadi2020}
M.~Marvi-Mashhadi, C.S. Lopes, and J.~LLorca.
\newblock {High fidelity simulation of the mechanical behavior of closed-cell
  polyurethane foams}.
\newblock {\em Journal of the Mechanics and Physics of Solids}, 135:103814,
  2020.

\bibitem{Sun2015}
Y.~Sun and Q.M. Li.
\newblock {Effect of entrapped gas on the dynamic compressive behaviour of
  cellular solids}.
\newblock {\em International Journal of Solids and Structures}, 63:50--67,
  2015.

\bibitem{Zhang2015}
J.~Zhang, G.P. Zhao, T.J. Lu, and S.Y. He.
\newblock {Strain rate behavior of closed-cell Al-Si-Ti foams: experiment and
  numerical modeling}.
\newblock {\em Mechanics of Advanced Materials and Structures}, 22(7):556--563,
  2015.

\bibitem{Su2018}
B.Y. Su, C.M. Huang, H.~Sheng, and W.Y. Jang.
\newblock {The effect of cell-size dispersity on the mechanical properties of
  closed-cell aluminum foam}.
\newblock {\em Materials Characterization}, 135:203--213, 2018.

\bibitem{Gebhart2019}
T.M.J. Gebhart, Di. Jehnichen, R.~Koschichow, M.~M{\"{u}}ller, M.~G{\"{o}}bel,
  V.~Geske, M.~Stegelmann, and M.~Gude.
\newblock {Multi-scale modelling approach to homogenise the mechanical
  properties of polymeric closed-cell bead foams}.
\newblock {\em International Journal of Engineering Science}, 145:103168, 2019.

\bibitem{Su2022}
B.Y. Su and W.Y. Jang.
\newblock {The microstructure characterization and elastic properties of
  closed-cell foams}.
\newblock {\em International Journal of Solids and Structures}, 257:111700,
  2022.

\bibitem{HossingerKalteis2022}
A.~H{\"{o}}ssinger-Kalteis, M.~Reiter, M.~Jerabek, and Z.~Major.
\newblock {Application of computed tomography data-based modelling technique
  for polymeric low density foams, Part A: Model development}.
\newblock {\em Journal of Cellular Plastics}, 58(3):429--448, 2022.

\bibitem{Ding2023}
Y.~Ding, X.~Zhou, J.~Wang, Y.~Feng, J.~Tang, N.~Shang, S.~Xin, Xi. Jian,
  M.~Gude, and J.~Xu.
\newblock {A sophisticated periodic micro-model for closed-cell foam based on
  centroidal constraint and capacity constraint}.
\newblock {\em Composite Structures}, 303:116175, 2023.

\bibitem{Gahlen2022a}
P.~Gahlen and M.~Stommel.
\newblock {Modeling of the local anisotropic mechanical foam properties in
  polyisocyanurate metal panels using mesoscale FEM simulations}.
\newblock {\em International Journal of Solids and Structures}, 244-245:111595,
  2022.

\bibitem{Gahlen2022b}
P.~Gahlen and M.~Stommel.
\newblock {Multiscale approach to determine the anisotropic mechanical
  properties of polyisocyanurate metal panels using FEM simulations}.
\newblock {\em Mechanics of Materials}, 174:104475, 2022.

\bibitem{Sonon2015}
B.~Sonon, B.~Fran{\c{c}}ois, and T.J. Massart.
\newblock {An advanced approach for the generation of complex cellular material
  representative volume elements using distance fields and level sets}.
\newblock {\em Computational Mechanics}, 56(2):221--242, 2015.

\bibitem{Ghazi2019}
A.~Ghazi, P.~Berke, K.~{Ehab Moustafa Kamel}, B.~Sonon, C.~Tiago, and T.J.
  Massart.
\newblock {Multiscale computational modelling of closed cell metallic foams
  with detailed microstructural morphological control}.
\newblock {\em International Journal of Engineering Science}, 143:92--114,
  2019.

\bibitem{Ghazi2020b}
A.~Ghazi, C.~Tiago, B.~Sonon, P.~Berke, and T.J. Massart.
\newblock {Efficient computational modelling of closed cell metallic foams
  using a morphologically controlled shell geometry}.
\newblock {\em International Journal of Mechanical Sciences}, 168:105298, 2020.

\bibitem{Huo2022}
X.~Huo, Z.~Jiang, Q.~Luo, Q.~Li, and G.~Sun.
\newblock {Mechanical characterization and numerical modeling on the yield and
  fracture behaviors of polymethacrylimide (PMI) foam materials}.
\newblock {\em International Journal of Mechanical Sciences}, 218:107033, 2022.

\bibitem{OstojaStarzewski2006}
M.~Ostoja-Starzewski.
\newblock {Material spatial randomness: From statistical to representative
  volume element}.
\newblock {\em Probabilistic Engineering Mechanics}, 21(2):112--132, 2006.

\bibitem{Reissner1961}
E.~Reissner and Y.~Stavsky.
\newblock {Bending and stretching of certain types of heterogeneous aeolotropic
  elastic plates}.
\newblock {\em Journal of Applied Mechanics}, 28(3):402--408, 1961.

\bibitem{Campello2003}
E.M.B. Campello, P.M. Pimenta, and P.~Wriggers.
\newblock {A triangular finite shell element based on a fully nonlinear shell
  formulation}.
\newblock {\em Computational Mechanics}, 31(6):505--518, 2003.

\bibitem{vanDijk2016}
N.P. van Dijk.
\newblock {Formulation and implementation of stress‐driven and/or
  strain‐driven computational homogenization for finite strain}.
\newblock {\em International Journal for Numerical Methods in Engineering},
  107(12):1009--1028, 2016.

\bibitem{Saadat2023}
M.A. Saadat and D.~Durville.
\newblock {A mixed stress-strain driven computational homogenization of spiral
  strands}.
\newblock {\em Computers {\&} Structures}, 279:106981, 2023.

\bibitem{Larsson2023}
C.~Larsson, F.~Larsson, J.~Xu, K.~Runesson, and L.E. Asp.
\newblock {Effects of lithium insertion induced swelling of a structural
  battery negative electrode}.
\newblock {\em Composites Science and Technology}, 244:110299, 2023.

\bibitem{Coenen2010}
E.W.C. Coenen, V.G. Kouznetsova, and M.G.D. Geers.
\newblock {Computational homogenization for heterogeneous thin sheets}.
\newblock {\em International Journal for Numerical Methods in Engineering},
  83(8-9):1180--1205, 2010.

\bibitem{Liu2021}
L.~Liu, A.~Sridhar, M.G.D. Geers, and V.G. Kouznetsova.
\newblock {Computational homogenization of locally resonant acoustic
  metamaterial panels towards enriched continuum beam/shell structures}.
\newblock {\em Computer Methods in Applied Mechanics and Engineering},
  387:114161, 2021.

\bibitem{Kouznetsova2001}
V.G. Kouznetsova, W.A.M. Brekelmans, and F.P.T. Baaijens.
\newblock {An approach to micro-macro modeling of heterogeneous materials}.
\newblock {\em Computational Mechanics}, 27(1):37--48, 2001.

\bibitem{Miehe2002}
C.~Miehe.
\newblock {Strain-driven homogenization of inelastic microstructures and
  composites based on an incremental variational formulation}.
\newblock {\em International Journal for Numerical Methods in Engineering},
  55(11):1285--1322, 2002.

\bibitem{Hill1963}
R.~Hill.
\newblock {Elastic properties of reinforced solids: Some theoretical
  principles}.
\newblock {\em Journal of the Mechanics and Physics of Solids}, 11(5):357--372,
  1963.

\bibitem{Alkhader2009}
M.~Alkhader and M.~Vural.
\newblock {The partition of elastic strain energy in solid foams and lattice
  structures}.
\newblock {\em Acta Materialia}, 57(8):2429--2439, 2009.

\bibitem{Ding2024}
J.~Ding, Q.~Ma, X.~Li, L.~Zhang, H.~Yang, S.~Qu, M.Y. Wang, W.~Zhai, H.~Gao,
  and X.~Song.
\newblock {Imperfection‐enabled strengthening of ultra‐lightweight lattice
  materials}.
\newblock {\em Advanced Science}, 2024.

\bibitem{Timoshenko1961}
S.P. Timoshenko and J.M. Gere.
\newblock {Buckling of thin plates}.
\newblock In {\em Theory of Elastic Stability}, pages 348--439. McGraw-Hill
  Book Co. Inc., New York, 1961.

\bibitem{Skeens2024}
J.W. Skeens and S.~Kyriakides.
\newblock {Crushing of a closed-cell polymeric foam under triaxial loading}.
\newblock {\em International Journal of Solids and Structures}, 291:112686,
  2024.

\bibitem{Quey2018}
R.~Quey, A.~Villani, and C.~Maurice.
\newblock {Nearly uniform sampling of crystal orientations}.
\newblock {\em Journal of Applied Crystallography}, 51(4):1162--1173, 2018.

\bibitem{Diab2023}
{Technical Data Divinycell H}.
\newblock Technical report, Diab Group, 2023.

\bibitem{Geuzaine2009}
C.~Geuzaine and J.‐F. Remacle.
\newblock {Gmsh: A 3‐D finite element mesh generator with built‐in pre‐
  and post‐processing facilities}.
\newblock {\em International Journal for Numerical Methods in Engineering},
  79(11):1309--1331, 2009.

\bibitem{Helfer2020}
T.~Helfer, J.~Bleyer, T.~Frondelius, I.~Yashchuk, T.~Nagel, and D.~Naumov.
\newblock {The MFrontGenericInterfaceSupport project}.
\newblock {\em Journal of Open Source Software}, 5(48):2003, 2020.

\bibitem{Logg2012}
A.~Logg, K.-A. Mardal, and G.~Wells.
\newblock {\em {Automated Solution of Differential Equations by the Finite
  Element Method: The FEniCS book}}.
\newblock Springer Science {\&} Business Media, 2012.

\bibitem{Bleyer2018}
J.~Bleyer.
\newblock {\em {Numerical Tours of Computational Mechanics with FEniCS}}.
\newblock Zenodo, 2018.

\bibitem{Gerard1957}
G.~Gerard and H.~Becker.
\newblock {Handbook of structural stability Part I: Buckling of flat plates}.
\newblock Technical report, New York University, 1957.

\bibitem{Shafiq2015}
M.~Shafiq, R.S. Ayyagari, M.~Ehaab, and M.~Vural.
\newblock {Multiaxial yield surface of transversely isotropic foams: Part
  II—Experimental}.
\newblock {\em Journal of the Mechanics and Physics of Solids}, 76:224--236,
  2015.

\bibitem{Funari2021}
M.F. Funari, S.~Spadea, P.~Lonetti, and P.B. Louren{\c{c}}o.
\newblock {On the elastic and mixed-mode fracture properties of PVC foam}.
\newblock {\em Theoretical and Applied Fracture Mechanics}, 112:102924, 2021.

\bibitem{Zhang2012}
S.~Zhang, J.M. Dulieu-Barton, R.K. Fruehmann, and O.T. Thomsen.
\newblock {A methodology for obtaining material properties of polymeric foam at
  elevated temperatures}.
\newblock {\em Experimental Mechanics}, 52(1):3--15, 2012.

\bibitem{Taher2012}
S.T. Taher, O.T. Thomsen, J.M. Dulieu-Barton, and S.~Zhang.
\newblock {Determination of mechanical properties of PVC foam using a modified
  Arcan fixture}.
\newblock {\em Composites Part A: Applied Science and Manufacturing},
  43(10):1698--1708, 2012.

\bibitem{Chen2013}
L.~Chen and M.S. {Hoo Fatt}.
\newblock {Transversely isotropic mechanical properties of PVC foam under
  cyclic loading}.
\newblock {\em Journal of Materials Science}, 48(19):6786--6796, 2013.

\bibitem{HooFatt2022}
M.S. {Hoo Fatt} and A.R. Vedire.
\newblock {Mechanical properties of marine polymer foams in the arctic
  environment}.
\newblock {\em Marine Structures}, 86:103308, 2022.

\bibitem{Tong2022}
X.~Tong, M.S. {Hoo Fatt}, and A.R. Vedire.
\newblock {A new crushable foam model for polymer-foam core sandwich
  structures}.
\newblock {\em International Journal of Crashworthiness}, 27(5):1460--1480,
  2022.

\bibitem{Magliaro2023}
J.~Magliaro, F.~Rahimidehgolan, P.~Mohammadkhani, W.~Altenhof, and A.T. Alpas.
\newblock {Modular energy absorbing capabilities achieved with compounded
  deformation mechanisms in composite AA6061-T6/PVC foam structures}.
\newblock {\em Acta Mechanica}, 234(9):4217--4258, 2023.

\bibitem{Berger2017}
J.B. Berger, H.N.G. Wadley, and R.M. McMeeking.
\newblock {Mechanical metamaterials at the theoretical limit of isotropic
  elastic stiffness}.
\newblock {\em Nature}, 543(7646):533--537, 2017.

\bibitem{TancogneDejean2018}
T.~Tancogne‐Dejean, M.~Diamantopoulou, M.B. Gorji, C.~Bonatti, and D.~Mohr.
\newblock {3D plate‐lattices: an emerging class of low‐density metamaterial
  exhibiting optimal isotropic stiffness}.
\newblock {\em Advanced Materials}, 30(45), 2018.

\bibitem{Guo2024}
H.~Guo and J.~Zhang.
\newblock {Performance‐oriented and deformation‐constrained dual‐topology
  metamaterial with high‐stress uniformity and extraordinary plastic
  property}.
\newblock {\em Advanced Materials}, 2024.

\bibitem{Nakshatrala2007}
K.B. Nakshatrala, A.~Masud, and K.D. Hjelmstad.
\newblock {On finite element formulations for nearly incompressible linear
  elasticity}.
\newblock {\em Computational Mechanics}, 41(4):547--561, 2007.

\bibitem{Mills1999}
N.J. Mills and H.X. Zhu.
\newblock {The high strain compression of closed-cell polymer foams}.
\newblock {\em Journal of the Mechanics and Physics of Solids}, 47(3):669--695,
  1999.

\bibitem{Almanza2001}
O.~Almanza, M.A. Rodriguez-Perez, and J.A. de~Saja.
\newblock {The microestructure of polyethylene foams produced by a nitrogen
  solution process}.
\newblock {\em Polymer}, 42(16):7117--7126, 2001.

\bibitem{Mills2009}
N.J. Mills, R.~St{\"{a}}mpfli, F.~Marone, and P.A. Br{\"{u}}hwiler.
\newblock {Finite element micromechanics model of impact compression of
  closed-cell polymer foams}.
\newblock {\em International Journal of Solids and Structures},
  46(3-4):677--697, 2009.

\end{thebibliography}

\end{document}